\documentclass[11pt,onecolumn,final,reqno]{article}

\usepackage[margin=1in]{geometry}
\usepackage{graphicx}
\usepackage{subcaption}  
\usepackage{enumerate}
\usepackage{paralist}
\usepackage{xcolor} 
\usepackage[hypcap]{caption}
\usepackage{hyperref}

\usepackage[small,compact]{titlesec}
\titleformat{\subsubsection}[runin]
  {\normalfont\normalsize\bfseries}{\thesubsubsection}{1em}{}
  
\usepackage{amsmath}
\usepackage{amsthm}
\usepackage{amsfonts}
\usepackage{amssymb}
\usepackage{mathrsfs}  
\usepackage{mathtools}
\usepackage{bm} 
\usepackage{bbm} 

\newcommand{\mc}{\mathcal}
\newcommand{\msf}{\mathsf}
\newcommand{\mbb}{\mathbb}  	
\newcommand{\mf}{\mathfrak}

\newcommand{\N}{ \ensuremath{\mathbb{N}}}  
\newcommand{\Z}{ \ensuremath{\mathbb{Z}}}  

\newcommand{\R}{ \ensuremath{\mathbb{R}}}  
\newcommand{\C}{ \ensuremath{\mathbb{C}}}  

\DeclarePairedDelimiterX\set[1]{\lbrace}{\rbrace}{  #1 }  

\DeclarePairedDelimiterX\norm[1]{\lVert}{\rVert}{#1}  			
\DeclarePairedDelimiterX\inner[2]{\langle}{\rangle}{#1 \,,\, #2}  	
\newcommand{\transp}{\mathsf{T}}  								
\DeclareMathOperator{\linspan}{span} 							
\newcommand{\mat}[1]{\mathbf{#1}}

\DeclarePairedDelimiterX\abs[1]{\lvert}{\rvert}{#1} 

\newcommand{\ind}[1]{\mathbb{I}_{#1}}							
\DeclareMathOperator*{\argmin}{arg\,min}

\DeclarePairedDelimiterX\floor[1]{\lfloor}{\rfloor}{#1}

\usepackage[linesnumbered,ruled,vlined]{algorithm2e}

\theoremstyle{plain}

\theoremstyle{remark}

\theoremstyle{definition}

\usepackage[T1]{fontenc}

\newcommand{\blue}[1]{{#1}}

\newcommand{\koop}{\mathcal{K}}

\begin{document}
\title{Koopman Reduced Order Modeling with Confidence Bounds}
\author{Ryan Mohr, Maria Fonoberova, Igor Mezi\'{c}}
\date{}

\maketitle

\begin{abstract}
This paper introduces a reduced order modeling technique based on Koopman operator theory that gives confidence bounds on the model's predictions. It is based on a data-driven spectral decomposition of the Koopman operator. The reduced order model is constructed using a finite number of Koopman eigenvalues and modes, while the rest of spectrum is treated as a noise process. This noise process is used to extract the confidence bounds. Additionally, we propose a heuristic algorithm to choose the number of deterministic modes to keep in the model.
\end{abstract}

\section{Introduction}

Koopman operator theory and Dynamic Mode Decomposition (DMD) have become ubiquitous in the data-driven analysis of complex and often high-dimensional systems \cite{mezic:2005, rowleyetal:2009,budisicetal:2012,brunton2021modern}. Given a vector-valued observable, the Koopman operator can be decomposed into eigenvalues, eigenfunctions, and Koopman modes plus the operator's continuous spectrum \cite{mezic:2005}. Dynamic Mode Decomposition is an algorithm that can compute approximations to the Koopman modes. When the continuous spectrum exists, DMD tries to approximate it with modes and eigenvalues, when it would be better approximated by a stochastic process. We take this idea and try to find an optimal number of modes and model the rest as a stochastic process.
    
The result of this representation, first proposed in \cite{mezic:2005}, is that instead of only being able to make point predictions, we can make predictions in a distributional sense using the coherent part and the statistics of the residual process \blue{to get confidence bounds on the prediction}. 



\paragraph{Related literature.}
\blue{
The method in this paper constructs a reduced order model by projecting onto a low-dimensional subspace spanned by a finite number of Koopman modes and represents the unresolved dynamics with a stochastic process. Similar methods have been proposed in turbulence and climate studies \cite{Trefethen1993HydrodynamicSW, StochasticityandSpatialResonanceinInterdecadalClimateFluctuations, Reynolds_Hussain_1972, doi:10.1073/pnas.96.26.14687, FargeSneider:2001, Butler1992ThreedimensionalOP, 10.1063/1.1398044}.
}

\blue{
A related concept in Koopman operator theory is that of the stochastic Koopman operator associated with a random dynamical system \cite{mezic:2005}. Given a probability space $(\Omega, \mathfrak F, P)$ and a family of random maps $\{T_\omega\}$ parameterized by $\omega \in \Omega$ the stochastic Koopman operator computes the expectation of an observable over realizations of the random maps by $\koop f = \mbb E_P[f \circ T_\omega(x)] = \int_\Omega f \circ T_\omega(x) dP$. The difference is that the stochastic Koopman operator starts from a random dynamical system and computes expectations of the evolution of the observable versus the method in this paper that starts from a dynamical system and represents it as a deterministic part plus a stochastic process.
}

\blue{Our method discards a subset of the computed DMD modes and eigenvalues in favor of representing those modes with a stochastic process. While our goal is to get confidence bounds on the predictions, one could argue that the discarded modes do not adequately represent the dynamics and a stochastic process better represents the unresolved dynamics. Previous methods, \cite{drmac2018data,colbrook2024rigorous,colbrook2023residual} truncate modes and eigenvalues by estimating the residual $\norm{Kg - \lambda g}$ and discarding those modes and eigenvalues whose residual is greater than a specified accuracy. The result of these methods is that they discard potentially useful, if noisy, parts of the dynamics.

}

In \cite{lu2020prediction}, the authors investigate the prediction accuracy of Extended Dynamic Mode Decomposition (EDMD) on parabolic PDEs. They show that the local truncation error between the true solution and EDMD can be bound by a term that only depends on the number of snapshots used for EDMD and which goes to zero as the number of snapshots goes to zero.

The problem of approximating the Koopman operator in a reproducing kernel Hilbert space (RKHS) was taken up in \cite{philipp2023error}. A cross-covariance operator acting on the RKHS is defined via as an integral operator involving the Koopman operator and the reproducing kernel. Bounds between the cross-covariance operator and its empirical approximation are derived and, using these results, a bound on a certain norm between the Koopman operator and its empirical approximation are found.

In \cite{lian2020gaussian}, the authors define a Koopman operator over a Gaussian process $\msf{GP}(\mu, k)$. The Gaussian process is an infinite-dimensional distribution over the space of real-valued functions $f :\R^n \to \R$. The Gaussian process is characterized by specifying, a priori, the mean $\mu : \R^N \to \R$ and the covariance functions $k : \R^{N\times N} \to \R$. Assuming deterministic dynamics, the authors show that if an observable is distributed according to the Gaussian process, $f \sim \msf{GP}(\mu, k)$, then the Koopman operator applied to $f$ is also a Gaussian process such that $U^t f = \msf{GP}(U^t \mu, U^t k)$. The limitations of this approach are that it requires the dynamics to be deterministic and the mean $\mu$ and covariances $k$ must be specified a priori rather than being  learned from data.

\paragraph{Contributions of this paper.}
We build off reduced order modeling techniques using the Koopman operator \cite{Koopman:1931,LasotaandMackey:1994,SinghandManhas:1993,MezicandBanaszuk:2004,mezic:2005}. In many cases, this linear operator has a spectral expansion which decomposes the operator into the point spectrum (eigenvalues) and their associated (Koopman) modes and the continuous spectrum \cite{MezicandBanaszuk:2004,mezic:2005,budisicetal:2012}. Much recent work has been done on algorithms that compute an approximation of the point spectrum of the Koopman operator \cite{Schmid:2010,rowleyetal:2009,tu2013dynamic,bagheri2014effects,williamsetal:2015,hemati2017biasing,drmac2018data,drmac2019data,drmac2020least}. These methods are purely data-driven and are based off computing a spectral decomposition of the operator without requiring \blue{an analytic representation of the infinite-dimensional operator}. However, it is known that these algorithms can return spurious eigenvalues due to numerical issues or trying to capture the continuous part of the spectrum with eigenvalues.



    In this paper, we develop a method, based on Koopman operator theory, that constructs a model and also gives confidence bounds on its predictions. \blue{For measure-preserving dynamical systems, the Koopman operator $\koop$ splits into a singular component, $\koop_s$, and a regular component, $\koop_r$ \cite{mezic:2005}. The singular part contains only the point spectrum while the regular part contains only the continuous spectrum. Low-dimensional projects can be obtained by retaining a finite number of modes and introducing stochastic terms to account for the neglected modes \cite{mezic:2005}. The algorithm in this paper implements the above idea by computing a Koopman Mode Decomposition on the training data; a finite number of the modes are selected as part of the ROM; the rest of the dynamics are modeled as a noise process, splitting the part of the noise in the subspace spanned by the Koopman modes used for the ROM (called ``in-plane'' for short) and the part orthogonal to those modes. The in-plane noisy dynamics are used to estimate confidence bounds on the predictions. In this paper, we are concerned in situations where we only have a single trajectory of the system and we wish to forecast future states.}

The contributions of this paper can be summarized by the following points:
\begin{compactenum}[(1)]
\item A method of providing confidence bounds on the predictions of a reduced order Koopman model based on the estimated noise distribution of the residual parts of the dynamics.
\item A heuristic method for determining the minimum number of modes to be used for the reduced order model. The heuristic is based on the Shapiro-Wilk hypothesis test comparing how well the empirical distribution is represented by a Gaussian distribution. \blue{More sophisticated methods could use Kernel Density Estimation or Gaussian mixture models.}
\item Application of the method to stochastic systems with changing network topology.
\end{compactenum}

The rest of the paper is structured as follows. In section \ref{sec:krom}, we develop the mathematical framework for the reduced order modeling framework with confidence bounds. Section \ref{sec:examples} applies the methodology to a sequence of examples, two of which have dynamically changing network topology. Section \ref{sec:heuristic} gives a heuristic algorithm for selecting how many modes to keep for the nominal dynamics so one is not modeling the continuous spectrum with the point spectrum. We conclude the paper in section \ref{sec:conclusions}.

\section{Koopman Reduced Order Modeling}\label{sec:krom}

Consider a dynamical system
	\begin{equation}\label{eq:stoch-ds}
	\mat x(t+1) = T(\mat x(t)), \quad t\in \N
	\end{equation}
where $\mat x(t) \in \Omega \subset \R^d$ or $\C^d$, and $T: \Omega \to \Omega$ is a (possibly) nonlinear map. Let $\mc F$ be a linear space of real- or complex-valued functions on $\Omega$. The induced Koopman operator $\koop : \mc F \to \mc F$ is defined by the composition operation
	\begin{equation}
	\koop f(\mat x) = f(T(\mat x))
	\end{equation}
If $f_i \in \mc F$ for $i=1,\dots, n$, we define $\mat f = (f_1,\dots, f_n)^{\mathsf{T}}$ and define the action of $\koop$ on $\mat f$ as $\koop \mat f = (\koop f_1,\dots, \koop f_n)^{\mathsf{T}}$.

In this paper, we assume there is a spectral decomposition for $\koop$ (see \cite{mezic:2005} which states precise conditions for this to hold) that can be written as 
	\begin{equation}
	\koop^t f_i(\mat x) = \sum_j c_{j} \lambda_{j}^t \phi_{j}(\mat x) + \int_\C z^t dE_{c}(z)f_i(\mat x),
	\end{equation}
where $c_{j} \in \C$, $\lambda_{j}\in \C$ is an eigenvalue, $\phi_{j} \in \mc F$ is an eigenfunction, and $dE_{c}$ is a projection-valued measure corresponding to the continuous part of the spectrum. For the vector-valued observable case, we have
	\begin{equation}
	\koop^t \mat f(\mat x) = \sum_j \mat m_{j} \lambda_{j}^t \phi_{j}(\mat x) + \int_\C z^t dE_{c}(z) \mat f(\mat x),
	\end{equation}
where $\mat m_{j}$ is the $j$-th Koopman mode.

\subsection{Model reduction}

With finite data, we cannot compute the infinite expansion. Instead we compute a finite-dimensional approximation of the point spectrum and model the rest of the evolution of the observable as a stochastic process. That is
	\begin{equation}
	\koop^t \mat f(\mat x) = \sum_{j=1}^{J} \mat m_j \lambda_j^t \phi_j(\mat x) + \mat r(t),
	\end{equation}
where here $\mat r(t)$ is the residual part of the dynamics that we cannot resolve \blue{with a finite number of modes and contains the unresolved point spectrum and the continuous spectrum. The finite sequence formed by the residual dynamics forms empirical} distribution $\mf R$ (e.g. a Gaussian distribution). Additionally, we order the modes $\mat m_j$ according to their norm so that $\norm{\mat m_j} \geq \norm{\mat m_{j+1}}$. \blue{In case the modes are normalized, the eigenvalues can be ranked by their associated reconstruction coefficient. We choose an order $J$ for the ROM and take the highest ranked modes and eigenvalues.} This implies that the modes and eigenvalues of a lower-order model are subsets of those of higher-order models. To make the notation more compact, let us denote the finite-mode expansion as
	\begin{equation}\label{eq:nominal-dyn}
	D^t \mat f(\mat x) = \sum_{j=1}^{J} \mat m_j \lambda_j^t \phi_j(\mat x).
	\end{equation}
	

So far we have split the evolution of the observable into a finite-mode expansion and the residual dynamics. Now we look into splitting the residual dynamics into a  modal component and an innovation component
Let $M = \linspan\set{\mat m_1, \dots, \mat m_J}$, $\mat M$ the matrix having the Koopman modes as its  columns, and $P_M$ the orthogonal projection onto $M$. Then,
	\begin{equation}\label{eq:kmd-with-noise}
	\koop^t \mat f(\mat x) = \underbrace{\sum_{j=1}^{J} \mat m_j \lambda_j^t \phi_j(\mat x) + P_M \mat r(t)}_{\in M} + \underbrace{(I - P_M)\mat r(t)}_{\text{orthogonal to $M$}}.
	\end{equation}
Consider the terms lying in $M$:
	\begin{align}
	\sum_{j=1}^{J} \mat m_j \lambda_j^t \phi_j(\mat x) + P_M \mat r(t) 
	= \mat M \begin{bmatrix} \lambda_1^t \phi_1(\mat x) \\ \vdots\\ \lambda_J^t \phi_J(\mat x)\end{bmatrix} + \mat M \mat M^+ \mat r(t),
	\end{align}
where we have used the fact that $P_M$ can be computed as $P_M = \mat M \mat M^+$, where $\mat M^+$ is the Moore-Penrose pseudoinverse of $\mat M$.
We define the \textbf{modal noise} $\rho(t) = (\rho_1(t), \dots, \rho_J(t))^{\msf T}$ as
	\begin{equation}
	\rho(t) = \mat M \mat M^+ \mat r(t).
	\end{equation}
The \textbf{innovation noise} is then defined as 
	\begin{equation}
	\eta(t) = (I - P_M)\mat r(t) = \mat r(t) - P_M \mat r(t) = \mat r(t) - \rho(t).
	\end{equation}
The innovation noise $\eta(t)$ is orthogonal to $M$. With these definitions, \eqref{eq:kmd-with-noise} becomes
	\begin{equation}\label{eq:stoch-rom}
	\koop^t \mat f(\mat x) =  \left( \mat M\begin{bmatrix} \lambda_1^t \phi_1(\mat x) \\ \vdots\\ \lambda_J^t \phi_J(\mat x)\end{bmatrix} + \rho(t) \right) + \eta(t).
	\end{equation}

\subsection{Decomposition into deterministic dynamics and stochastic dynamics}
%

Algorithm \ref{alg:compute-rom} summarizes the basic computations to compute a reduced-order model. \blue{The novelty here is that the reduced order model has a deterministic part and a stochastic part --- thus modeling a (possibly deterministic) dynamical system as a sum of a ''coherent'' and ''stochastic'' parts. In theory the ''stochastic'' part is related to the continuous part of the Koopman spectrum.} The algorithm returns the deterministic part of the model, the residual noise sequence, the modal noise sequence, and the innovation noise sequence. \blue{Since the modes are normalized, the mode weights in step 2 are defined as
\begin{align}
    \mat V &= [\mat m_1, \dots, \mat m_M] \\
    c_j &= (\mat V^{+} \mat f(x))_j \label{eq:coefficient-of-projection} \\
    w_j &= \frac{\abs{c_j}}{\max\{\abs{c_j}\}} \label{eq:mode-weights}
\end{align}
Note that $c_j$ is the $j^{th}$ coefficient of projection of the first snapshot $\mat f(x)$ onto all of the Koopman modes. Once the number of modes $J$ to retain for the ROM is chosen, the projection coefficients are recomputed using \eqref{eq:model-reconstruction-coeff}.
}


\begin{algorithm}[H]\label{alg:compute-rom}
\DontPrintSemicolon
  
  \KwData{$\set{\koop^t \mat f(\mat x)}_{t=0}^{T}$ evolution of observable}
  \KwResult{ \\ Modal model: $D = \set{(c_j, \lambda_j, \mat m_j) : j=1,\dots J}$ \\ Noise sequences: $\set{(\mat r(t), \rho(t), \eta(t)) : t=0,\dots , T-1 }$}  
  
  Compute KMD of data $\to (\set{\mat m_j}_{j=1}^{M}, \set{\lambda_j}_{j=1}^{M} )$.
  
  Order the Koopman modes $\set{\mat m_j}_{1}^M$ by their their weights $\{w_j\}$ so that $w_{j} \geq w_{j+1}$, where $w_j$ is given by \eqref{eq:mode-weights}.
  
  Once ordered, normalize each mode to have norm 1.
  
  Choose truncation number $J > 0$ according to the method in section \ref{sec:heuristic}.
  
  Compute modal reconstruction coefficients: 
  \begin{equation}\label{eq:model-reconstruction-coeff}
  \mat c = (c_1, \dots, c_J)^\transp = \argmin_{\alpha = (\alpha_1, \dots, \alpha_J) \in \C^J} \sum_{t=0}^{T-1} \norm{\koop^t \mat f(\mat x) - \sum_{j=1}^J \alpha_j \lambda_j^t \mat m_j}_2^2.
  \end{equation}
  
  Construct deterministic reconstruction model
	\begin{equation}
	D^t \mat f(\mat x) = \sum_{j=1}^J c_j \lambda_j^t \mat m_j.
	\end{equation}
	
  Compute residual sequence
	\begin{equation}
	\mat r(t) = \koop^t \mat f(\mat x) - D^t \mat f(\mat x).
	\end{equation}

  Compute modal noise sequence
	\begin{equation}
	\rho(t) = \mat M \mat M^+ \mat r(t).
	\end{equation}
	
  Compute innovation sequence
	\begin{equation}
	\eta(t) = \mat r(t) - \rho(t).
	\end{equation}

\caption{Koopman Reduced Order Model with Confidence Bounds (KROM-CB)}
\end{algorithm}

\subsection{Forecasting with Confidence Bounds}

\blue{
Once the noise sequences $\mat r(t), \rho(t), \eta(t)$ are computed, empirical probability density models need to be constructed for which drawing samples is easy, for example by fitting a parameterized density function like the Gaussian distribution or using other methods such as kernel density estimation techniques taking the form $p_{\rho}(x) = \frac{1}{nh}\sum_{i=1}^{n} K(\frac{x - \rho(t_i)}{h})$, where $K$ is the kernel function (such as a Gaussian), $h$ is the bandwidth, and $\rho(t_i)$ are the samples of the residual sequence. Assume that the densities for $\rho$ and $\eta$ are computed.

Recall that $M$ is the subspace spanned by the Koopman modes, $M = \linspan\set{\mat m_1,\dots, \mat m_J}$, and $P_M$ is the orthogonal projection onto $M$. By \eqref{eq:stoch-rom} and since the innovation sequence $\eta$ is orthogonal to $M$ we have
	\begin{equation}\label{eq:projected-ROM-woth-modal-noise}
	P_M \koop^t \mat f(\mat x) = D^t \mat f(\mat x) + \rho(t)
	\end{equation}
as the form of the reduced order model. Note that the reduced-order model is split into a deterministic part and a stochastic part.
    
In our experiments, we assume that the resulting statistics of the modal sequence form a Gaussian distribution (which can be checked; see Section \ref{sec:heuristic}) and construct a prediction interval at each time $t$. We use the finite-mode expansion to give the mean of the prediction and use $\pm$ 2 standard deviations of the estimated Gaussian to give the prediction interval at each point. Thus, the prediction interval defined by the reduced order model at each time is 
    \begin{equation}\label{eq:ROM-prediction-interval}
        I(t) = [D^t \mat f(\mat x) - 2\msf{std}(\rho), D^t \mat f(\mat x) + 2\msf{std}(\rho)]
    \end{equation}

To measure the quality of the ROM, we compute both a reconstruction error and the percent time that the real trajectory were within 2 standard deviations of the deterministic model's prediction. This residence time can be computed for every coordinate in $\mat f$.
}

\blue{
In \eqref{eq:ROM-prediction-interval} and the experiments below, we plot $(K^t \mat f)_i$ and $(D^t \mat f)_i \pm 2\sigma_{\rho, i}$, where $\sigma_{\rho, i}$ is the estimated standard deviation of the $i^{th}$ coordinate of the modal noise $\{\rho(t) \in \R^n : t=0, \dots, T\}$. The true signal $\koop^t \mat f$ generally lies within these bounds, despite the true signal containing contributions from both the modal and innovation noise: $\koop^t \mat f = D^t \mat f + \rho(t) + \eta(t)$. The below computation shows why we can only consider the modal noise.

Let $\koop^t \mat f \in \R^n$ ($t\in \N_0$) be the true evolution and $D^t \mat f$ be the deterministic model in the modal subspace given by 
\begin{equation}
D^t \mat f = M \Lambda^t c_{\mat f}
\end{equation}
where $c_{\mat f} = M^+ \mat f \in \C^p$, $M \in \C^{n \times p}$ is the matrix of Koopman modes and $\Lambda \in \C^{p\times p}$ is the diagonal matrix of Koopman eigenvalues.

The residual and modal noise are given, respectively, by
\begin{align}
    r(t) &= \koop^t \mat f - D^t \mat f \\
    \rho(t) &= P_{M} r(t) = P_M (\koop^t \mat f - D^t \mat f) = P_M\koop^t \mat f - D^t \mat f,
\end{align}
 where $P_M = M M^+$ is the orthogonal projection onto $\linspan M$, $M^+$ is the pseudoinverse, and $P_M D^t \mat f = D^t \mat f$ since $D^t \mat f \in \linspan M$.

\paragraph{Modal noise standard deviation.}
We estimate the standard deviation of the modal noise from the sequence $\{e_i^* \rho(t) : t=0,\dots, T\}$, where $e_i$ is the $i^{th}$ canonical basis vector in $\R^n$. From the definitions,
\begin{equation}
    e_i^*\rho(t) = e_i^* P_M \koop^t \mat f - e_i^* D^t \mat f.
\end{equation}
Assuming $\mbb E e_i^*\rho(t) = 0$, the unbiased estimator for the variance of the $i^{th}$ coordinate is
\begin{align}
    \hat\sigma_{\rho, i}^2 &= \frac{1}{T} \sum_{t=0}^T \abs{e_i^* \rho(t)}^2 = \frac{1}{T} \sum_{t=0}^T \abs{e_i^* P_M \koop^t \mat f - e_i^* D^t \mat f}^2.
\end{align}

We plot $e_i^* \koop^t \mat f$ and $e_i^* D^t \mat f$ for some coordinate $i$ and see if they are within $\pm 2\hat\sigma_{\rho, i}$. Therefore, we need to bound the error between the true signal and the deterministic part of the ROM in terms of the standard deviations of the modal noise.

\begin{align}
    \abs{e_i^* \koop^t \mat f - e_i^* D^t \mat f} 
    &= \abs{e_i^* (I - P_M) \koop^t \mat f + e_i^* P_M \koop^t \mat f - e_i^* D^t \mat f} \nonumber \\
    &\leq \abs{e_i^* (I - P_M) \koop^t \mat f} + \abs{e_i^* P_M \koop^t \mat f - e_i^* D^t \mat f} \nonumber \\
    &= \abs{e_i^* (I - P_M) \koop^t \mat f} + \left(\abs{e_i^* P_M \koop^t \mat f - e_i^* D^t \mat f}^2 \right)^{1/2}.
\end{align}
Taking the expectation, we have
\begin{align}
    \mbb E\left[\abs{e_i^* \koop^t \mat f - e_i^* D^t \mat f}\right] 
    &:= \frac{1}{T+1} \sum_{t=0}^T \abs{e_i^* \koop^t \mat f - e_i^* D^t \mat f} \nonumber \\
    &\leq \frac{1}{T+1} \sum_{t=0}^T \abs{e_i^* (I - P_M) \koop^t \mat f} + \frac{1}{T+1} \sum_{t=0}^T \left(\abs{e_i^* P_M \koop^t \mat f - e_i^* D^t \mat f}^2 \right)^{1/2}
\end{align}
Applying Jensen's inequality to the last term on the right gives
\begin{align}
    \mbb E\left[\abs{e_i^* \koop^t \mat f - e_i^* D^t \mat f}\right] 
    &\leq \frac{1}{T+1} \sum_{t=0}^T \abs{e_i^* (I - P_M) \koop^t \mat f} +  \left(\frac{1}{T+1}\sum_{t=0}^T \abs{e_i^* P_M \koop^t \mat f - e_i^* D^t \mat f}^2 \right)^{1/2} \nonumber \\
    &\leq \frac{1}{T+1} \sum_{t=0}^T \abs{e_i^* (I - P_M) \koop^t \mat f} +  \hat\sigma_{\rho, i}.
\end{align}
Therefore, if 
\begin{equation}
    \frac{1}{T+1} \sum_{t=0}^T \abs{e_i^* (I - P_M) \koop^t \mat f} \leq \hat\sigma_{\rho, i},
\end{equation}
then $e_i^* \koop^t \mat f$ and $e_i^* D^t \mat f$ are within $\pm 2\hat\sigma_{\rho, i}$ in expectation.

From the definition $\eta(t) = r(t) - \rho(t)$, we can show that $\eta(t) = (I - P_M)\koop^t \mat f$. Using the same assumption as above ($\mbb E e_i^* \eta(t) = 0$) and Jensen's inequality, it can be shown that
\begin{equation}
    \frac{1}{T+1} \sum_{t=0}^T \abs{e_i^* (I - P_M) \koop^t \mat f} \leq \hat\sigma_{\eta, i}.
\end{equation}
Therefore,
\begin{align}
    \mbb E\left[\abs{e_i^* \koop^t \mat f - e_i^* D^t \mat f}\right] 
    &\leq \hat\sigma_{\eta, i} + \hat\sigma_{\rho, i}.
\end{align}
If $\mat f$ is in the invariant subspace of Koopman eigenfunctions, then $(I-P_M)\koop^t \mat f = 0$ and $\hat\sigma_{\eta, i} = 0$. If $\mat f$ is nearly invariant with respect to the subspace, then $\norm{(I-P_M)\koop^t \mat F}$ is small on the training set and we can check numerically that $\hat \sigma_{\eta, i} \leq \hat \sigma_{\rho, i}$. In either case, we have
\begin{equation}\label{eq:2std-modal-noise-bound}
    \mbb E\left[\abs{e_i^* \koop^t \mat f - e_i^* D^t \mat f}\right] \leq 2\hat\sigma_{\rho, i}.
\end{equation}
}

\section{Performance and Comparisons}\label{sec:examples}

\subsection{Linear Modal System}

We first test the algorithm for a system with known modes and eigenvalues. While in practice, it would not make sense to represent a known linear finite-dimensional model with the Koopman operator, this example acts as a sanity check. Additionally, in practice we only have trajectory data and may not know the underlying dynamics is linear. This example also shows that when lifting to a higher-dimensional space, spurious eigenvalues can be introduced. While everything was formulated with respect to deterministic dynamical systems, we will test both with a small amount of additive Gaussian noise to the dynamics. The general model of the system is given by
	\begin{equation}\label{eq:linear-modal-model}
	\mat x(t) = \left(\sum_{j=1}^{J} \mat m_j \lambda_j^t \right) + \xi(t) \in \C^n
	\end{equation}
where $\mat m_j \in \R^n$, $\lambda_j \in \C$, $t\in \N$, and $\xi(t) \in \R^n$ is a noise sequence with each coordinate independent and identically distributed with a normal distribution. We note that the evolution of this system is in $\C^n$ rather than $\R^n$ due to the choice of complex $\lambda_j$'s; therefore, \blue{the eigenvalues do not need to include the complex conjugates}.

For the following numerical results, the system parameters were chosen as $J = 20$, $n=40$, and $\xi_i(t) \sim N(0, 0.25^2)$. The values of the coordinates of each mode $\mat m_j$ were chosen from a uniform distribution between -1 and 1 and afterwards the mode was normalized to have norm 1. The associated eigenvalues were chosen uniformly at random in  the box $[-1, 1] \times [-1i , 1i]$ in the complex plane, then scaled to be on the unit circle.

The trajectory was simulated for 401 time steps resulting in a simulation data matrix in $\C^{40\times 401}$. To construct the observable $\mat f : X \to \R^N$ used for the Koopman model, we take a delay-embedding of the trajectory $\mat x(t)$ with a delay of 300. This gives a data matrix of size $\mat D \in \C^{12,000\times 101}$ as input for the reduced-order modeling algorithm. Therefore, there are a total of 100 Koopman modes in $\C^{12,000}$ for this model, even though $\mat x(t)$ is constructed using only 20 modes in $\C^{40}$. \blue{For all examples, the Cauchy-Vandermonde version of DMD was used \cite{drmac2019data}.}

\begin{figure}[ht!]
\centering
    \begin{subfigure}{0.49\textwidth}
        \includegraphics[width=\textwidth]{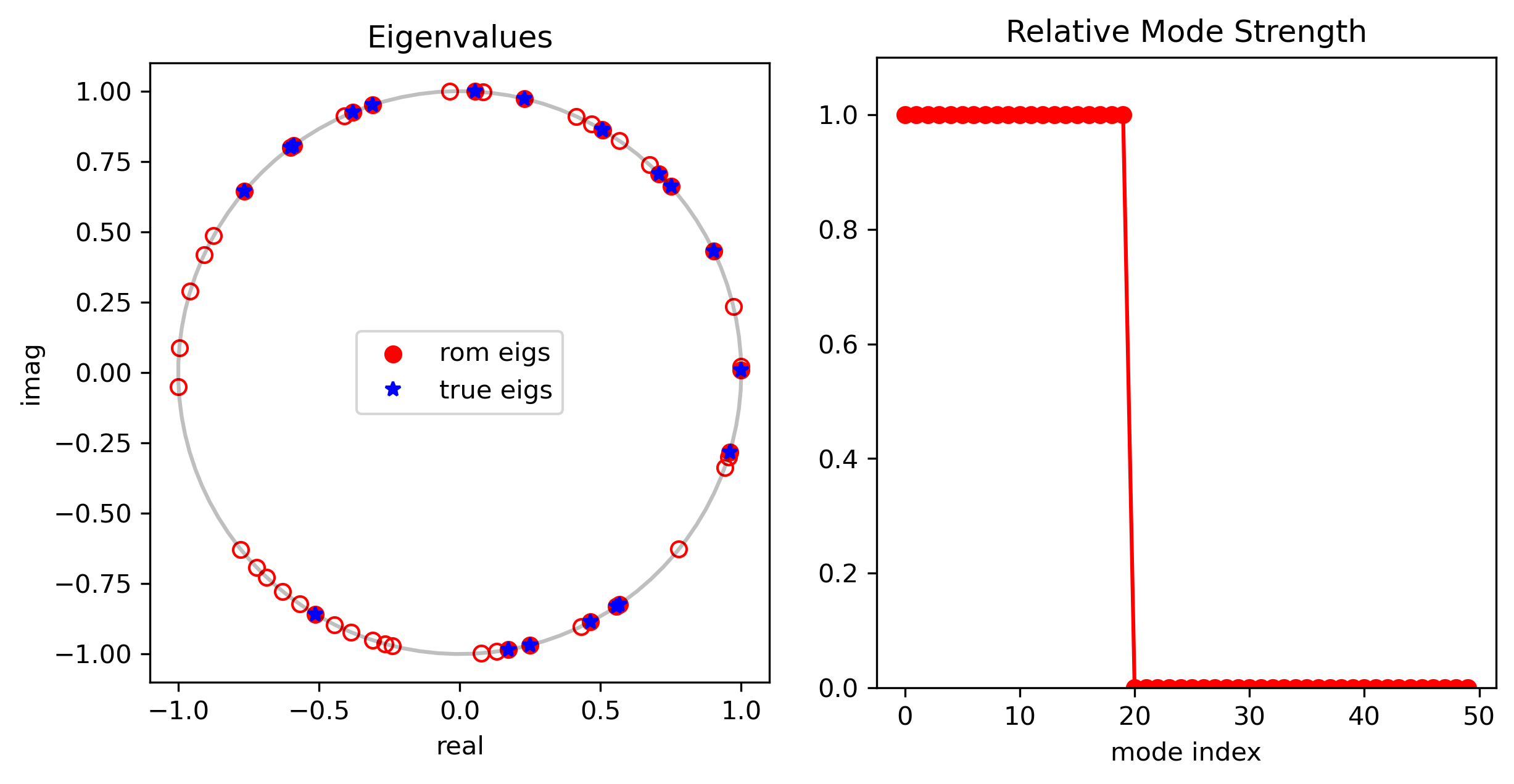}
        \caption{Noise $\xi(t) = 0$.}
        \label{fig:linear-modal-model-eigenvalues-no-noise}
    \end{subfigure}
    \hfill
    \begin{subfigure}{0.49\textwidth}
        \includegraphics[width=\textwidth]{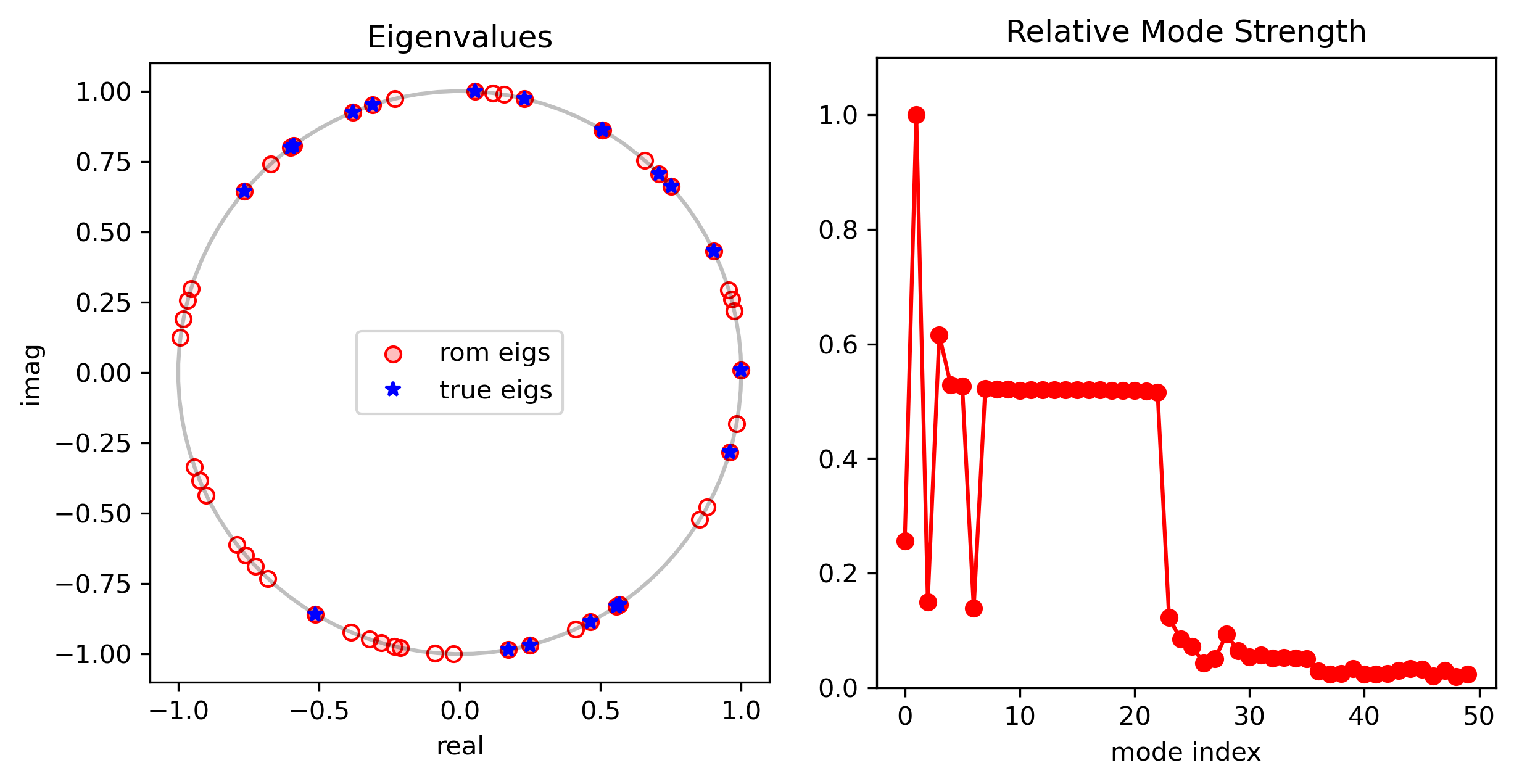}
        \caption{Standard deviation of $\xi(t) = 0.05$.}
        \label{fig:linear-modal-model-eigenvalues-0.05-noise}
    \end{subfigure}
    \hfill
    \begin{subfigure}{0.49\textwidth}
        \includegraphics[width=\textwidth]{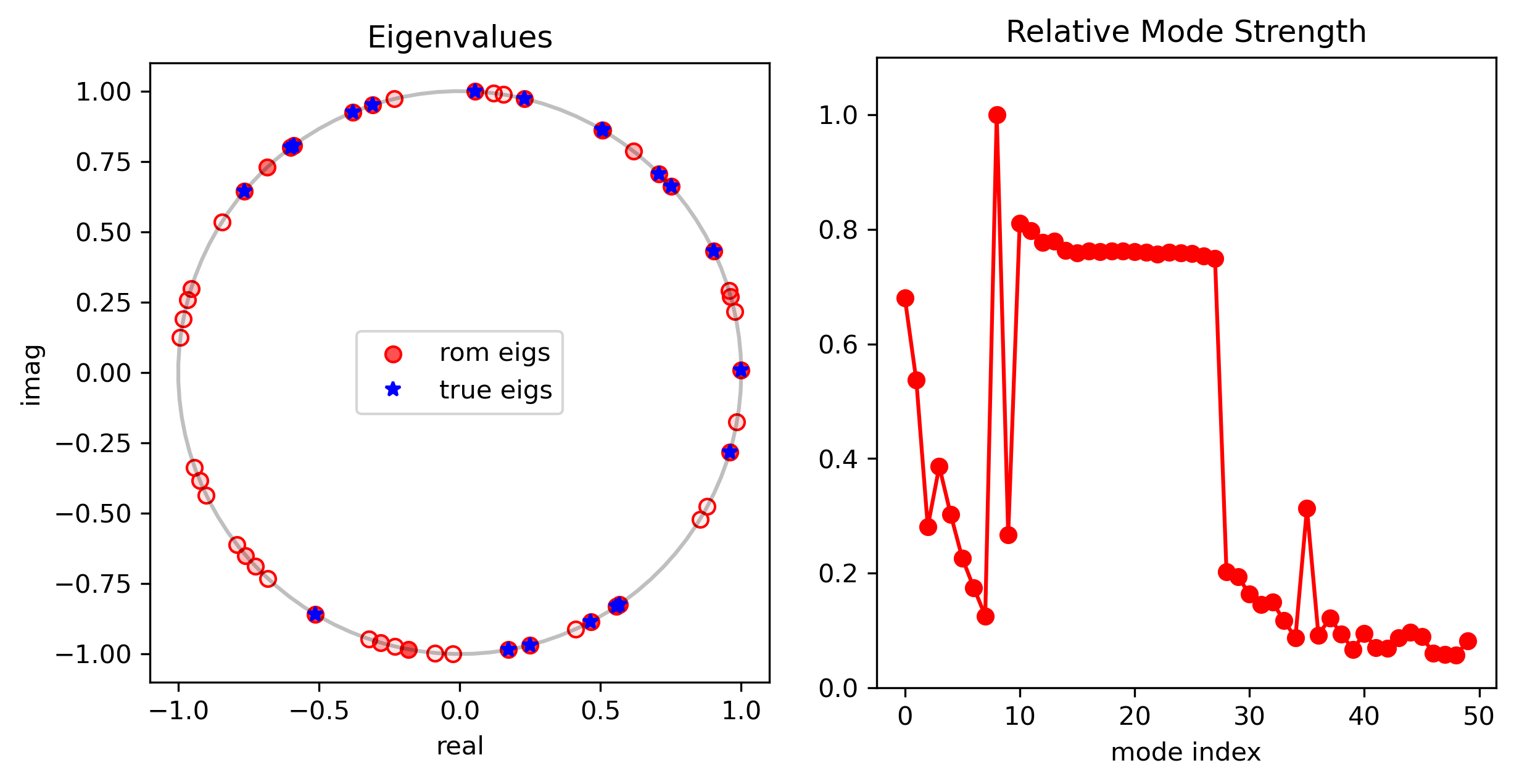}
        \caption{Standard deviation of $\xi(t) = 0.1$.}
        \label{fig:linear-modal-model-eigenvalues-0.1-noise}
    \end{subfigure}
    \hfill
    \begin{subfigure}{0.49\textwidth}
        \includegraphics[width=\textwidth]{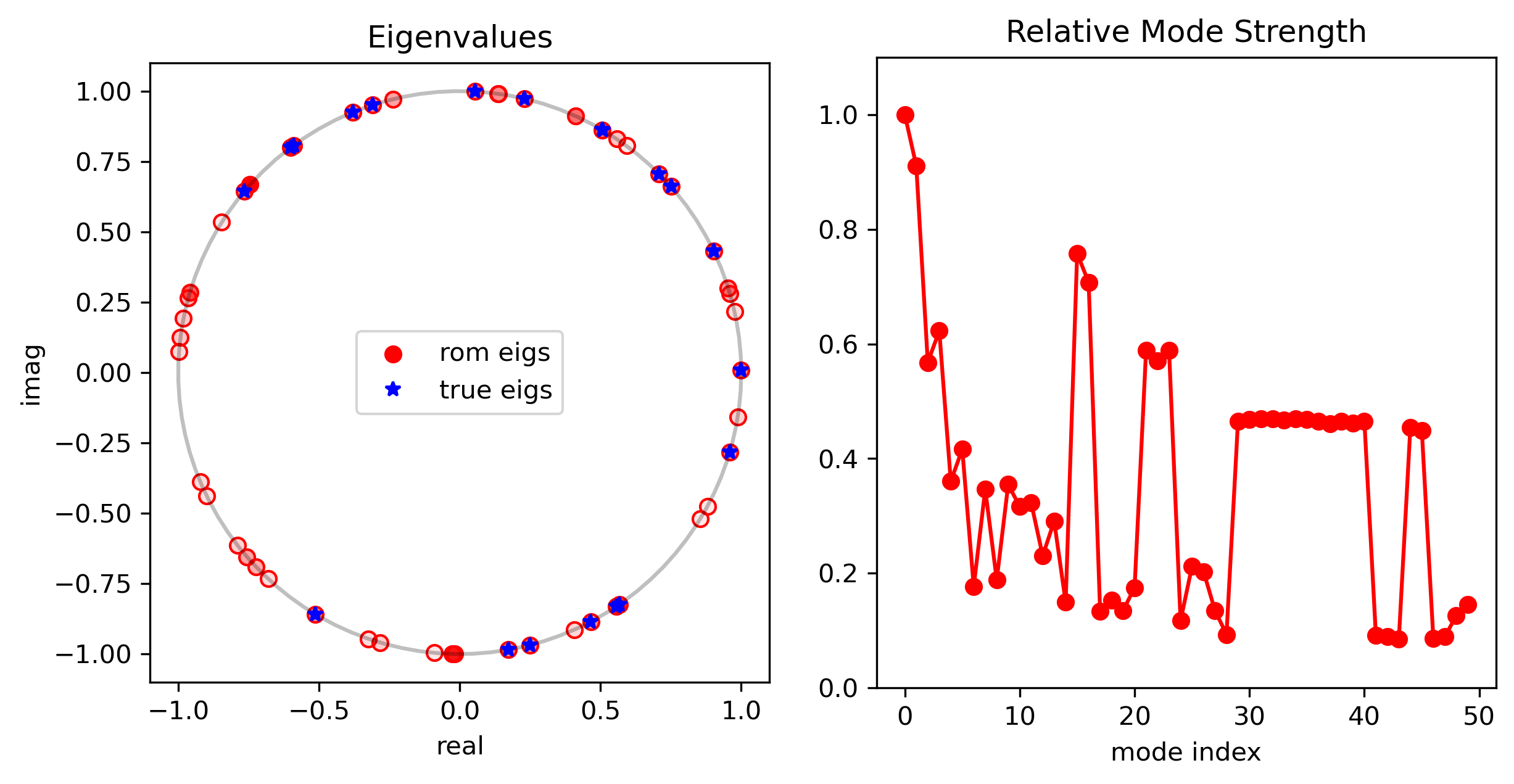}
        \caption{Standard deviation of $\xi(t)$ = 0.25.}
        \label{fig:linear-modal-model-eigenvalues-0_25-std-noise}
    \end{subfigure}
    \caption{\textbf{Linear Modal Model Eigenvalues}: \blue{Eigenvalues and relative mode weights for 50-mode ROM with various levels of noise. Modes and eigenvalues are ordered by their weights computed via \eqref{eq:mode-weights}. The relative strength of a mode in the ROM is determined by \eqref{eq:model-reconstruction-coeff} and is shown in the plots. ROM eigenvalues' (red) opacity is determined by the relative strength of the corresponding mode. The observable was a delay-embedding of the state vector with a delay of 300. When the noise is zero (a), the true eigenvalues are captured and the extra eigenvalues have (numerically) zero weight. With increasing the levels of noise, the true eigenvalues are still captured, but some of the extra eigenvalues have corresponding modes with non-trivial weights. As the noise increases, (b), (c), (d), the number of modes with non-trivial weight increases to capture the stochasticity in the underlying model.}}
    \label{fig:linear-modal-model-eigenvalues}
\end{figure}

\blue{
Figure \ref{fig:linear-modal-model-eigenvalues} compares the true eigenvalues of the model that generated the data and the eigenvalues computed when constructing the reduced order model. Recall that we have lifted the space from $\C^{40}$ to $\C^{12,000}$ using a delay embedding. Due to this, the true model's eigenvalues are a subset of the \blue{100} Koopman eigenvalues which can be seen in Figure \ref{fig:linear-modal-model-eigenvalues}. Figure \ref{fig:linear-modal-model-eigenvalues-no-noise} shows the computed Koopman eigenvalues for a 50-mode ROM when noise $\xi(t) = 0$. Note that the true eigenvalues are computed exactly and the weights of the other eigenvalues are zero, where the weights are computed using \eqref{eq:mode-weights} with the $c_j$'s computed in \eqref{eq:model-reconstruction-coeff}. Figure \ref{fig:linear-modal-model-eigenvalues-0_25-std-noise} shows the eigenvalues of a 50-mode ROM when the standard deviation of $\xi(t)$ is 0.25. Note that the true eigenvalues are still captured in this noisy case; however, there are extra eigenvalues in the Koopman ROM, not corresponding to the underlying system's eigenvalues, that have non-trivial weights due to the added noise. The number of modes with non-trivial weights increases with increasing levels of the noise due to the need to capture the extra stochasticity.
}

\begin{figure}[ht!]
\centering
    \begin{subfigure}{0.45\textwidth}
        \includegraphics[width=\textwidth]{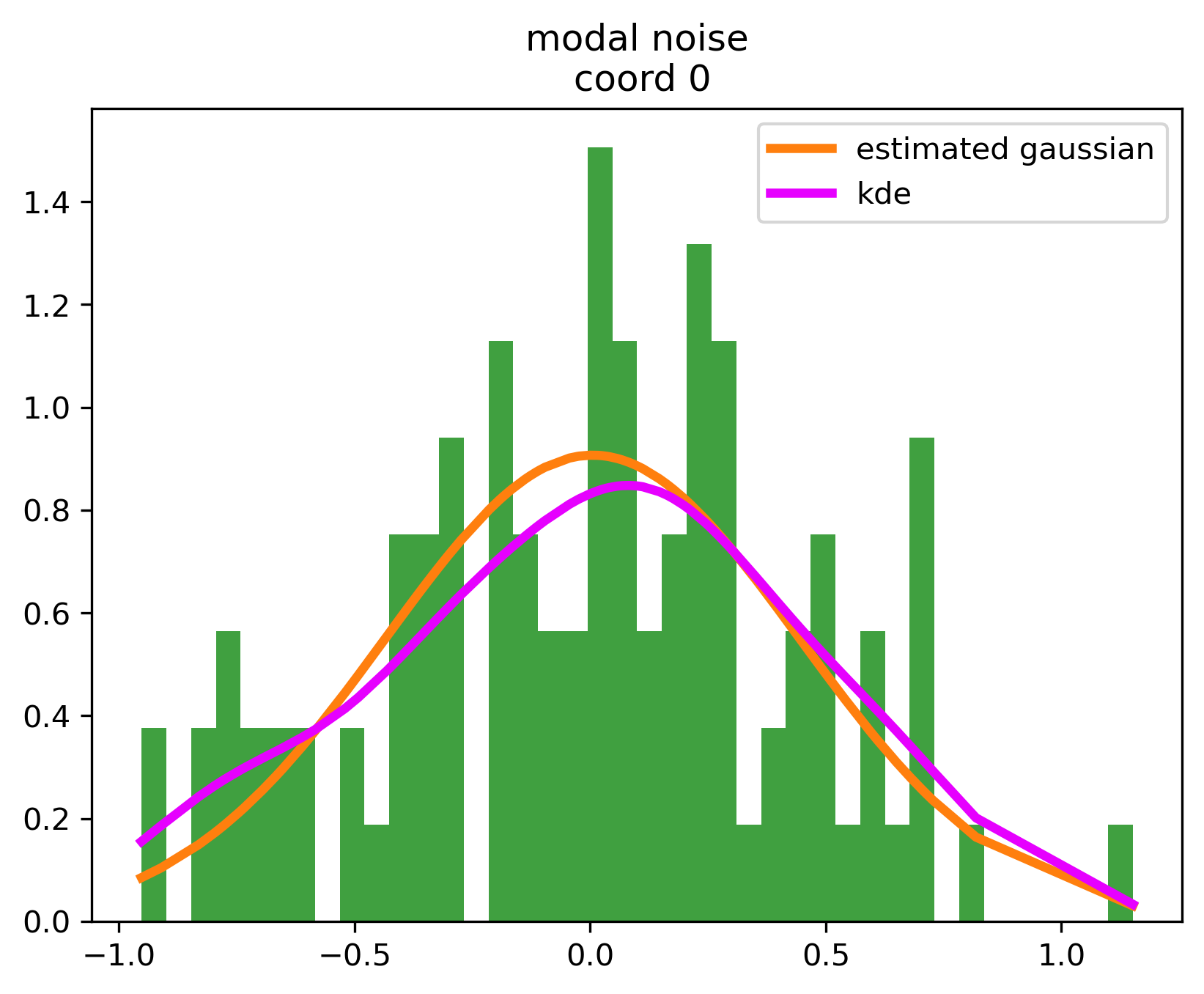}
        \caption{Modal distribution}
        \label{fig:linear-modal-model-modal-noise-coord-0} 
    \end{subfigure}
    \hfill
    \begin{subfigure}{0.45\textwidth}
        \includegraphics[width=\textwidth]{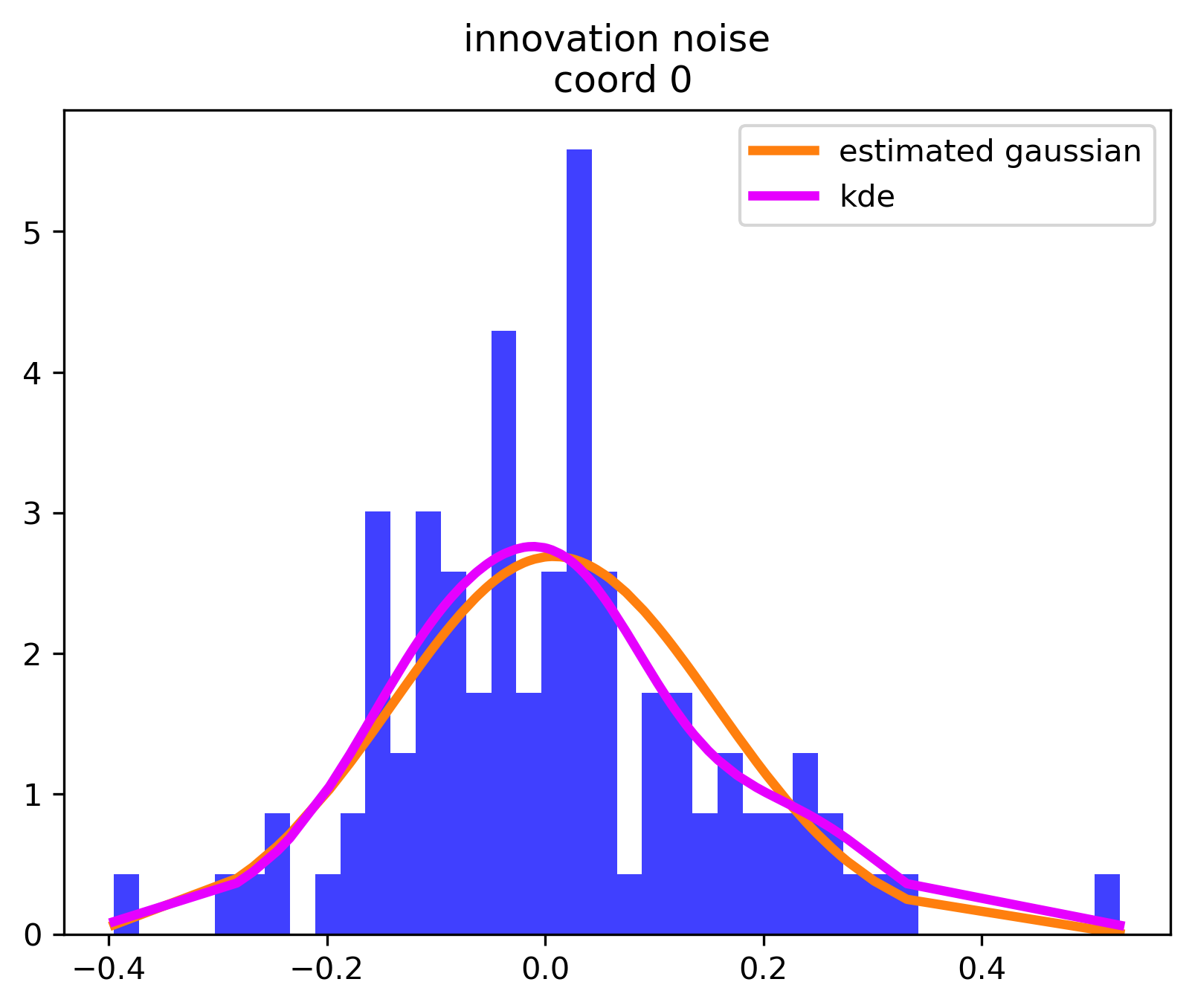}
        \caption{Innovation distribution}
\label{fig:linear-modal-model-innovation-noise-coord-0} 
    \end{subfigure}
    \caption{\textbf{Linear Modal Model Noise Distributions}: Noise distributions in the first coordinate for a 50-mode ROM. Standard deviation of $\xi$ in \eqref{eq:linear-modal-model} is 0.25. \blue{The observable was a delay-embedding of the state vector with a delay of 300.} (a) Modal noise. Green is the histogram, orange is a fitted Gaussian distribution, and pink is a kernel density approximation. (b) Innovation noise. Blue is the histogram, orange is a fitted Gaussian distribution, and pink is a kernel density approximation.}
    \label{fig:linear-modal-noise-distributions}
\end{figure}

\begin{figure}[ht!]
\begin{center}
\includegraphics[width=0.9\textwidth]{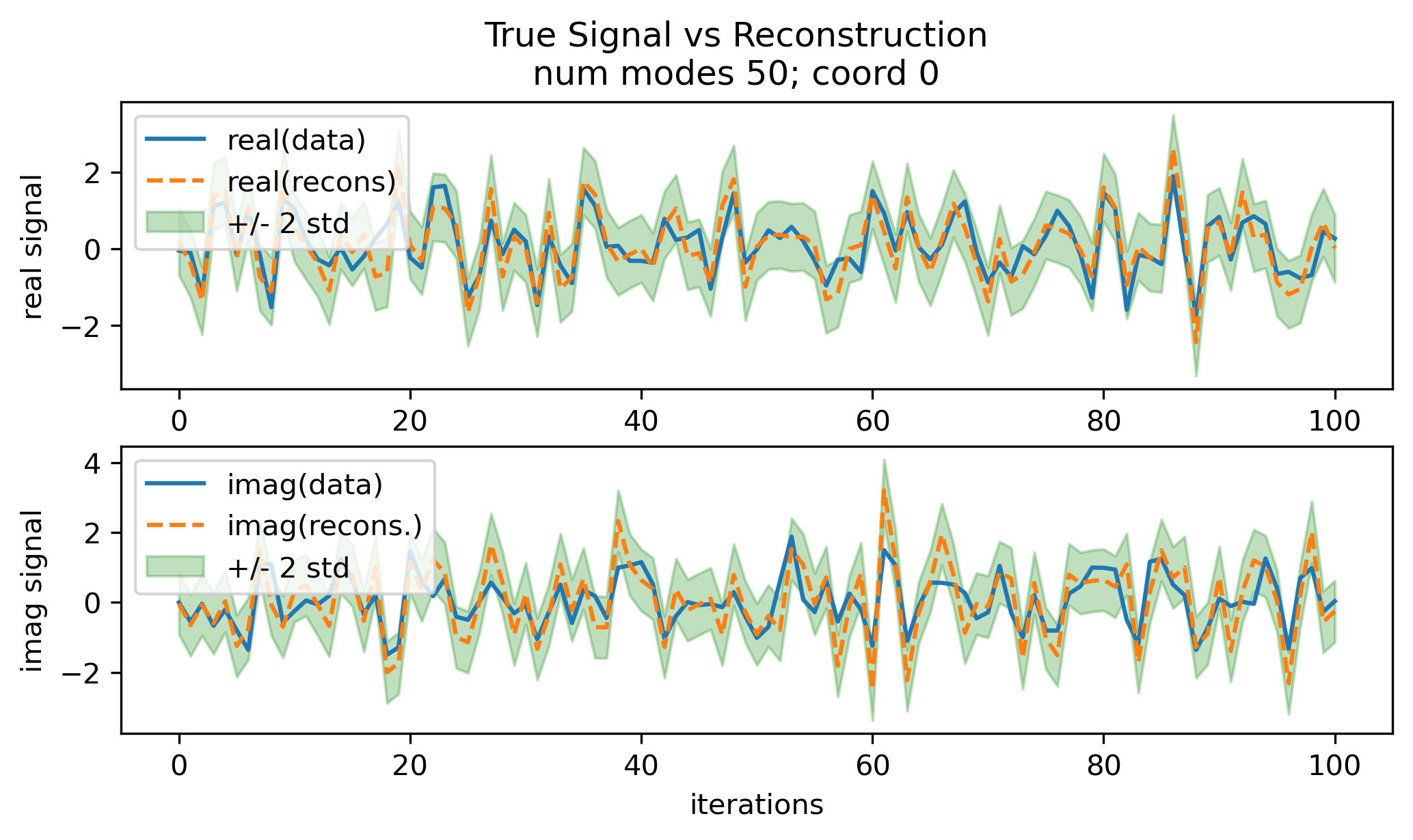}
\caption{\textbf{Linear Modal Model Reconstruction}: \blue{Comparison of the 50-mode projected ROM and confidence bounds (eq.'s \eqref{eq:projected-ROM-woth-modal-noise} and \eqref{eq:ROM-prediction-interval}) versus the true signal for the first coordinate. Standard deviation of $\xi$ in \eqref{eq:linear-modal-model} is 0.25. The true signal is given by the blue trace, the ROM reconstruction is the orange trace, and the green band is $\pm$2 standard deviations of the modal noise. The green band represent the 95\% confidence interval for the ROM reconstruction. The top trace is the real part of the signal. The bottom trace is the imaginary part.}}
\label{fig:linear-modal-model-rom-reconstruction-coord-0}
\end{center}
\end{figure}

\begin{figure}[ht!]
\centering
    \begin{subfigure}{0.45\textwidth}
        \includegraphics[width=\textwidth]{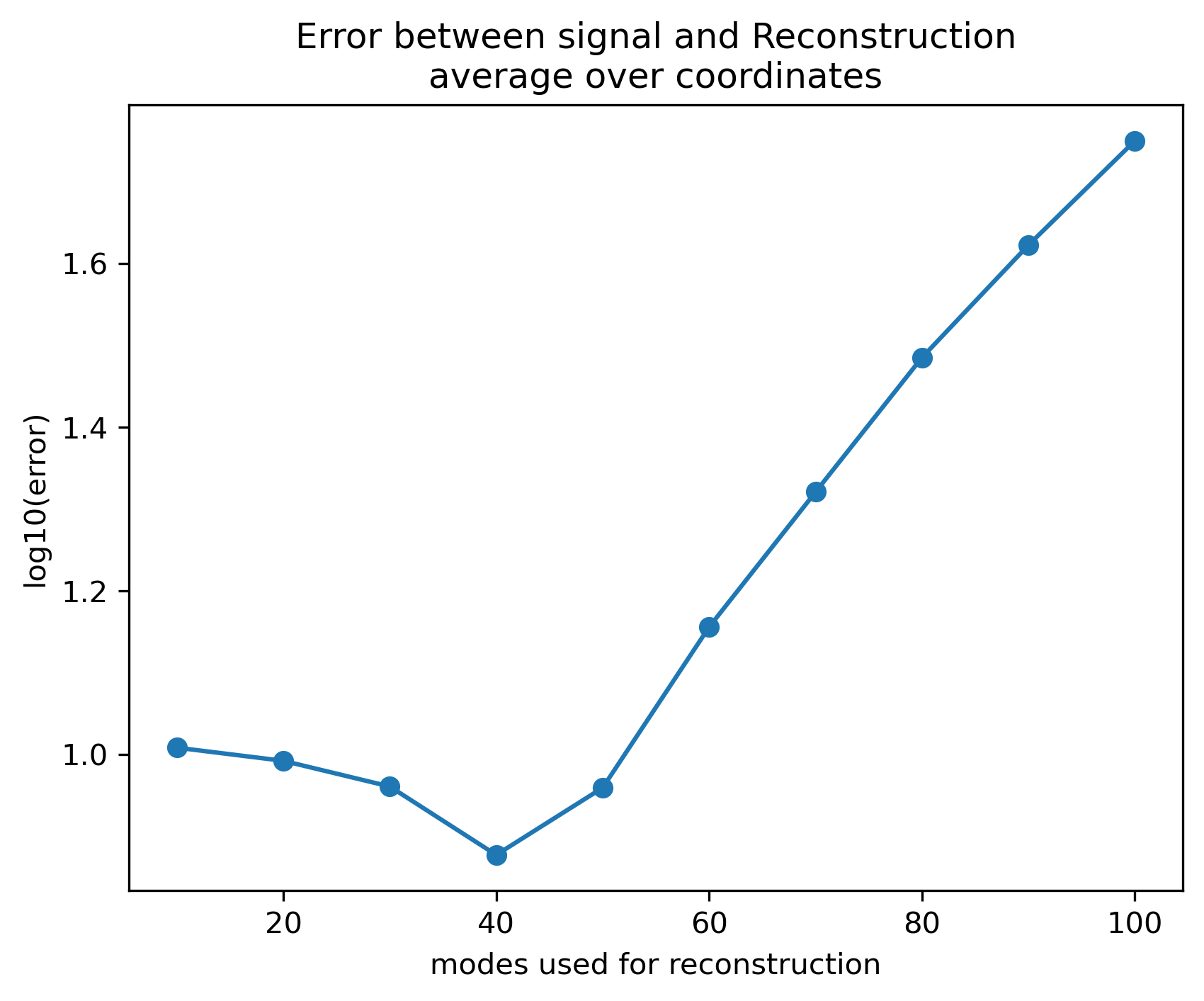}
        \caption{Log Error}
    \end{subfigure}
    \hfill
    \begin{subfigure}{0.45\textwidth}
        \includegraphics[width=\textwidth]{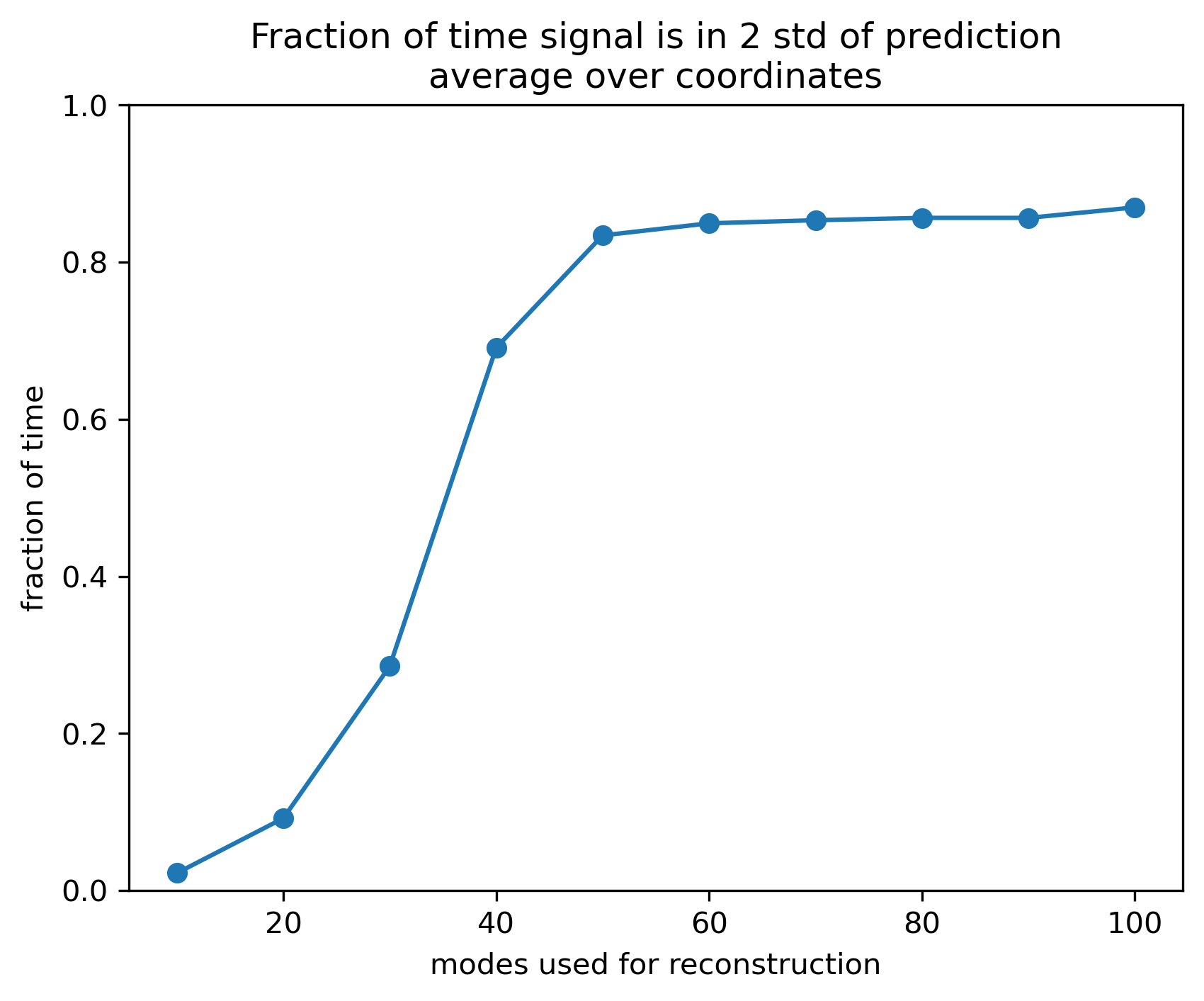}
        \caption{Residence times}
    \end{subfigure}
\caption{\textbf{Linear Modal Model Error and Residence Times}: \blue{Standard deviation of $\xi$ in \eqref{eq:linear-modal-model} is 0.25. For both the error and the residence times, the metrics were computed for each coordinate and then averaged over coordinates. (a) Error vs. number of modes in the ROM. (b) Fraction of time the true signal is in the confidence bounds vs.\ number of modes in the ROM. The error never drops too low (minimum at 40 modes) suggesting that the point prediction of the ROM is never that accurate. However, the residence times show a sharp transition and plateau around 50 modes suggesting that 50 modes capture the ``coherent'' part of the dynamics well and the rest of the dynamics are well-modeled by a stochastic term.}}
\label{fig:linear-modal-model-rom-error-residence-time-coord-0}
\end{figure}

Figure \ref{fig:linear-modal-noise-distributions} shows the corresponding modal noise distribution (noise in plane to the ROM modes) and the innovation noise distribution (orthogonal to the ROM modes). Figure \ref{fig:linear-modal-model-rom-reconstruction-coord-0} shows a comparison of the \blue{projected ROM, \eqref{eq:projected-ROM-woth-modal-noise}, and its confidence bounds, \eqref{eq:ROM-prediction-interval}, with the true signal} for just the first coordinate of the state vector. Note that the distributions for both the modal and innovation noises are approximately Gaussian. Additionally, the ROM signal matches the general trend of the true signal, with the true signal lying within the 95\% confidence interval of the ROM model. Here, we have shown only a single coordinate of the ROM and the true signal, but these results are typical for the other coordinates. We note that while the true model has additive Gaussian noise (in $\C^{40}$), the Koopman model's modal and innovation noise are in $\C^{12,000}$. Figure \ref{fig:linear-modal-model-rom-error-residence-time-coord-0} shows the ROM error with respect to the true signal and the fraction of time that the true signal is within $\pm 2$ standard deviations of the ROM's prediction. \blue{The error never drops too low (minimum at 40 modes) suggesting that the point prediction of the ROM is never that accurate. However, the residence times show a sharp transition and plateau around 50 modes suggesting that 50 modes capture the ``coherent'' part of the dynamics well and the rest of the dynamics are well-modeled by a stochastic term.}


\clearpage
\subsection{Switched Oscillators}

We apply the methodology to a switched anharmonic oscillator example. A single anharmonic oscillator can be written in action-angle coordinates as
	\begin{align}
		&\dot{I}(t) = 0 \\
		&\dot{\theta}(t) = f(I(t))
	\end{align}
where $I(t)$ is the action and $\theta(t)$ is the angle. Note that the rotational frequency $\dot{\theta}(t)$ depends on the action variable. For $f(I(t)) = I(t)$, the higher the action, the faster the rotation. The anharmonic oscillator was chosen due to its lack of point spectra other than the eigenvalue at 1 corresponding to the constant eigenfunctions. It can be shown \cite{mezic2020spectrum} that the rest of the spectrum is purely continuous with ``eigendistributions'' $\phi_j(I, \theta) = e^{ij\theta}\delta(I - c)$ supported on level-sets $\set{(c, \theta) : \theta \in S^1}$ of the action variable having ``eigenvalues'' $e^{ict}$. 

For this example, we use a collection of anharmonic oscillators and randomly wire pairs of them together at each integer time point. At each integer time point, energy is instantaneously transferred from the oscillator with higher action to the one with lower action. Between integer time points, the oscillators are decoupled.  Additionally, we incorporate an additive Gaussian noise term to the action variables. We investigate our prediction algorithm on such a system. Below, we give the construction of the model.

Let there be $j=0, \dots, n-1$ oscillators with action/angle dynamics given by 
	\begin{align}
	\dot I_j(t)  &= 0, \quad  (t \notin \N) \label{eq:I_dot} \\
	I_j(t) &= \sum_{i=0}^{n-1} \ind{E(t)}((i,j))\, c_{j,i}( I_i(k) - I_j(k) ) + \eta_j(t), \quad( t\in \N)  \label{eq:randomly-coupled-anharmonic-oscillators-1} \\ 
	\dot \theta_j(t) &= f( I_j(t) ) \label{eq:randomly-coupled-anharmonic-oscillators-2} 
	\end{align}
where $f(0) = 0$ and $f$ is a nondecreasing function, $0 < c_{j,i} < 1$ and $c_{i,j} = c_{j,i}$, $E(t)$ is the set of undirected edges denoting the connections between pairs of oscillators at time $t \in \N$, and $\ind{E(t)}$ is the indicator function for the set $E(t)$, and $\eta_j(t) \sim N(0, \sigma^2)$ is a normally distributed noise sequence. We use $f(I) = I$ in our examples. In everything that follows, \textbf{all edges are undirected edges}.

In our examples, the edge set dynamics are given by permutation dynamics. Let $\mat y(t), \mat z(t)  \in \Theta(n)$ where $\Theta(n)$ is the set of all permutations of the vector $(0,\dots, n-1)$. Write $\mat y(t) = (y_0(t),\dots, y_{n-1}(t))$ (similarly for $\mat z(t)$). Let $\mc P(n) = \{ \mat P:\Theta(n) \to \Theta(n)\}$ be the set of all permutation matrices on $\Theta(n)$. Fix $\mat Y, \mat Z \in \mc P(n)$. Then, the edge set dynamics are given by
	\begin{align}
	\mat y(t) &= \mat Y^{\floor{t}}\mat y \label{eq:p-permutation-dynamics-y} \\
	\mat z(t) &= \mat Z^{\floor{t}}\mat z \label{eq:q-permutation-dynamics-z} \\
	E(t) &=  \set{ (y_0(t), z_0(t)) }  \label{eq:edge-set-t}
	\end{align}
Here, the edge set changes at integer multiples in time. The edge set consists of a single, undirected edge which connects oscillators $y_0(t)$ and $z_0(t)$ at time $t$. \blue{If two oscillators are connected at any point in time, there future connections will be periodic due to the permutation dynamics.} The state space of the full model is
	\begin{equation}\label{eq:state-space}
	X = \set{ (\vec I, \vec \theta) \times (\mat y, \mat z) : \vec I \in (\R^+)^n, \vec\theta \in (\R/\Z)^n, \mat y, \mat z \in \Theta(n)},
	\end{equation}
where $\vec I = (I_0, \dots, I_{n-1})$ and $\vec \theta = (\theta_0, \dots, \theta_{n-1})$.

\paragraph{Experiments.}
The action \blue{variables' initial conditions $I_j(0)$} will be independent and identically distributed according to the exponential distribution defined as
	\begin{equation}\label{eq:exponential-distribution}
	h(t; \lambda) = \begin{cases}
	\lambda e^{-\lambda t}, & t \geq 0 \\
	0, & t < 0.
	\end{cases}
	\end{equation}
The mean of this exponential distribution is $1/\lambda$. For $\vec I = (I_1,\dots, I_n) \in (\R^+)^n$, we will take the $n$-product of exponential distributions
	\begin{equation}
	H(\vec I; \lambda) = \prod_{j=0}^{n-1} h(I_j; \lambda)
	\end{equation}
The distribution for the initial conditions will be
	\begin{equation}
	\mbb P = H(\vec I; 1) \otimes dist(\vec\theta) \otimes dist(\mat y) \otimes dist(\mat z).
	\end{equation}
where $dist(\vec\theta)$, $dist(\mat y)$, and $dist(\mat z)$ are the uniform distributions on their respective spaces.

Fix the \blue{vector-valued} observable $\mat f \in L^2(X; \R^{2n})$, which maps $X$ to \blue{$\R^{2n}$}, as 
	\begin{equation}\label{eq:example-observable}
	\mat f(\vec I, \vec \theta, \mat y, \mat z) = ( I_0, \dots, I_{n-1}, \theta_0, \dots, \theta_{n-1})^{\msf T} \in (\R^+)^n \times (\R/\Z)^n \subset \R^{2n}.
	\end{equation}
$\mat f$ does not directly observe the permutation variables $\mat y, \mat z$. Table \ref{tab:anharmonic-fixed-parameters-krom} specifies the parameters used in the simulation. When applying the Cauchy-Vandermonde DMD algorithm \cite{drmac2019data}, we use a delay embedding of the observable $\mat f$ . For all simulations, we use a delay embedding of 300 time steps that results in a data matrix of size $6,000\times 101$ for the DMD algorithm to act on. We therefore will be able to compute a maximum of 100 modes.

\begin{table}[h]
\caption{Fixed parameters for switch anharmonic simulation }
\begin{center}
\begin{tabular}{|ll|}
\hline
number of oscillators, $n$ & 10 \\ 
\hline
numerical time step, $dt$ & 0.05 sec \\
\hline
simulation time, $t$ & 20 sec \\
\hline
$c_{i,j}$, for all $i,j$ (see eq.\ \eqref{eq:randomly-coupled-anharmonic-oscillators-1}) & 0.5 \\
\hline
observable, $\vec f$ & eq.\ \eqref{eq:example-observable} \\ 
\hline
Hankel delay embedding & 300 \\
\hline
$\lambda$ & 1 \\
\hline
$\sigma$ & 0.05 \\
\hline
\end{tabular}
\end{center}
\label{tab:anharmonic-fixed-parameters-krom}
\end{table}%

\begin{figure}[ht]
	\centering
	\begin{subfigure}[b]{0.45\textwidth}
         \centering
         \includegraphics[width=\textwidth]{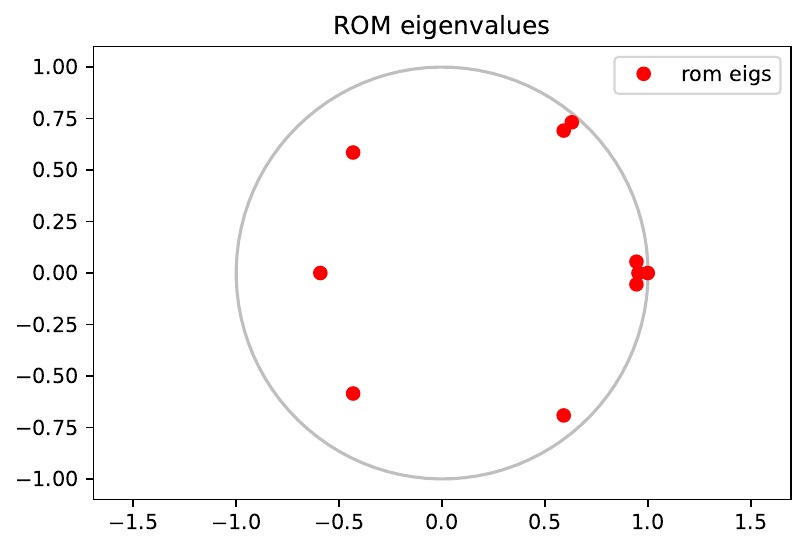}
         \caption{10 modes}
     \end{subfigure}
     \hfill
     \begin{subfigure}[b]{0.45\textwidth}
         \centering
         \includegraphics[width=\textwidth]{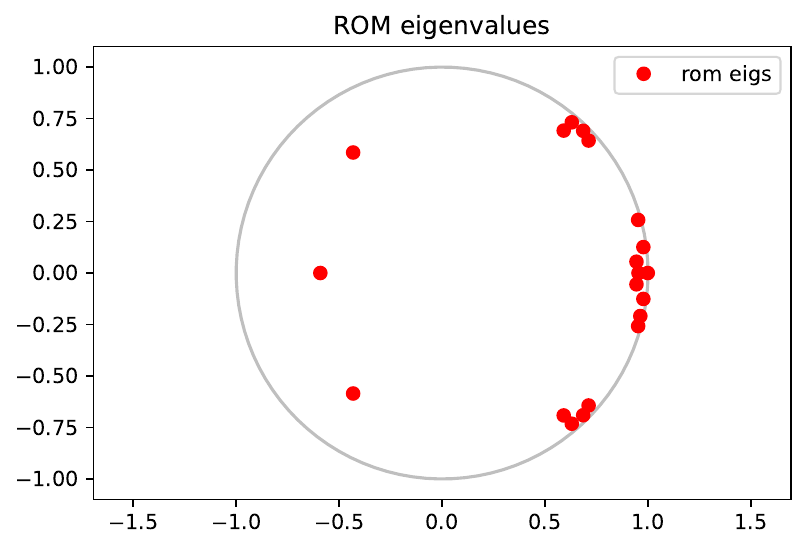}
         \caption{20 modes}
     \end{subfigure} 
     \hfill
     \begin{subfigure}[b]{0.45\textwidth}
         \centering
         \includegraphics[width=\textwidth]{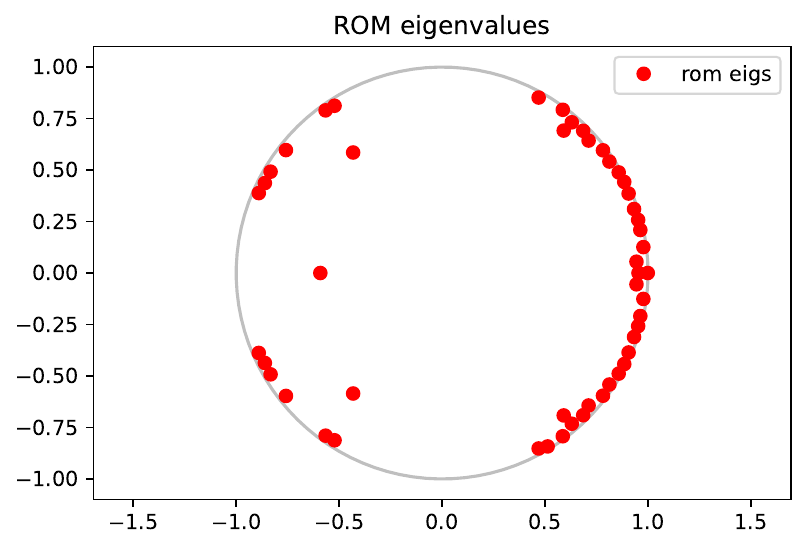}
         \caption{50 modes}
     \end{subfigure}
     \hfill
     \begin{subfigure}[b]{0.45\textwidth}
         \centering
         \includegraphics[width=\textwidth]{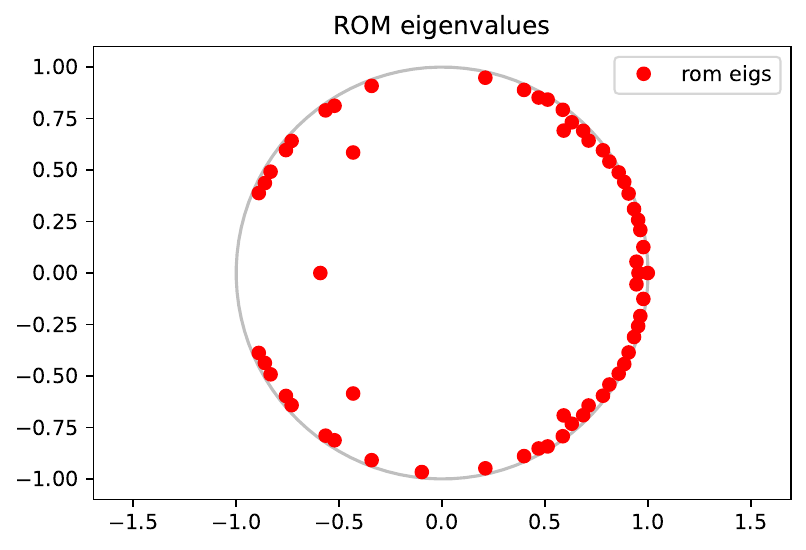}
         \caption{60 modes}
     \end{subfigure}
     \hfill
     \begin{subfigure}[b]{0.45\textwidth}
         \centering
         \includegraphics[width=\textwidth]{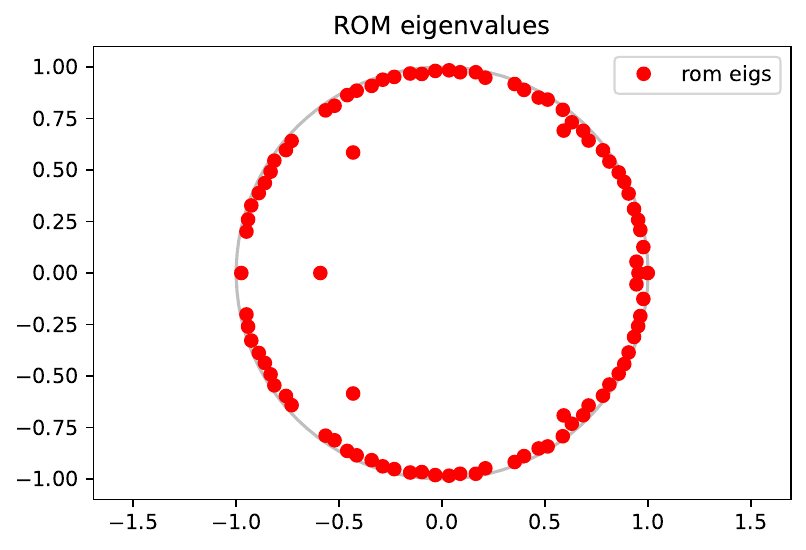}
         \caption{90 modes}
     \end{subfigure}
     \hfill
     \begin{subfigure}[b]{0.45\textwidth}
         \centering
         \includegraphics[width=\textwidth]{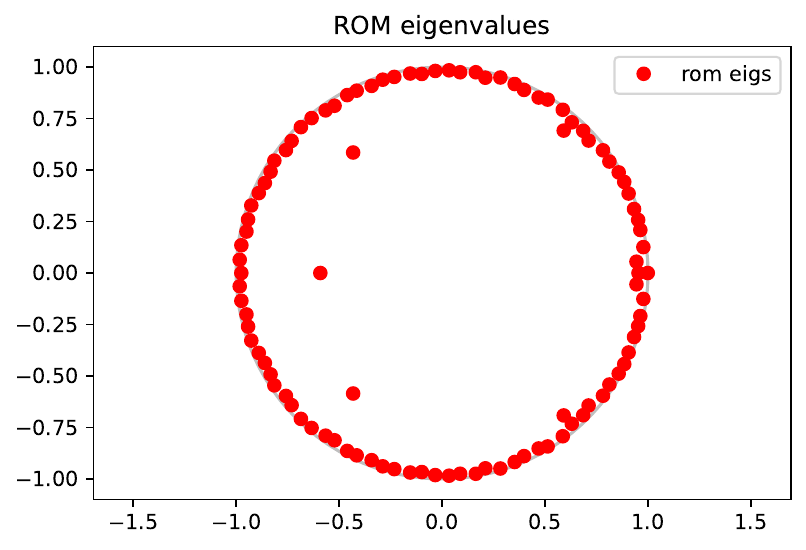}
         \caption{100 modes}
     \end{subfigure}
     \hfill
     \caption{\textbf{Anharmonic Oscillators}. Computed Koopman eigenvalues for the switched anharmonic oscillator model with added added Gaussian noise $(\sim N(0, 0.05^2)$. The eigenvalues of lower-order models are always a subset of those of higher-order models.}
\label{fig:anharmonic-eigenvalues}
\end{figure}

\begin{figure}[ht]
	\centering
	\begin{subfigure}[b]{0.45\textwidth}
         \centering
         \includegraphics[width=\textwidth]{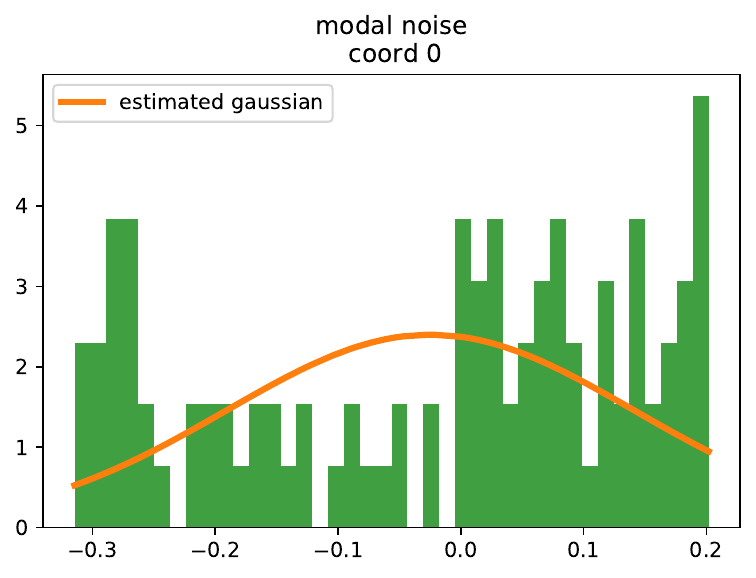}
         \caption{10 modes}
     \end{subfigure}
     \hfill
     \begin{subfigure}[b]{0.45\textwidth}
         \centering
         \includegraphics[width=\textwidth]{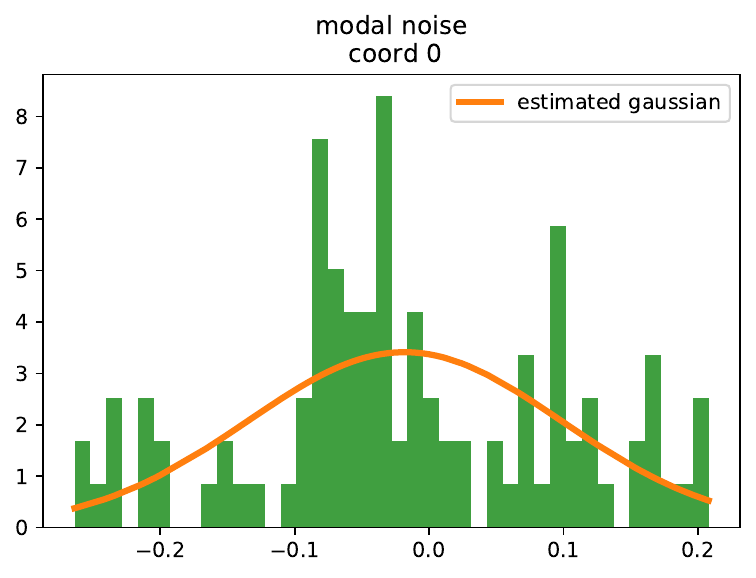}
         \caption{20 modes}
     \end{subfigure} 
     \hfill
     \begin{subfigure}[b]{0.45\textwidth}
         \centering
         \includegraphics[width=\textwidth]{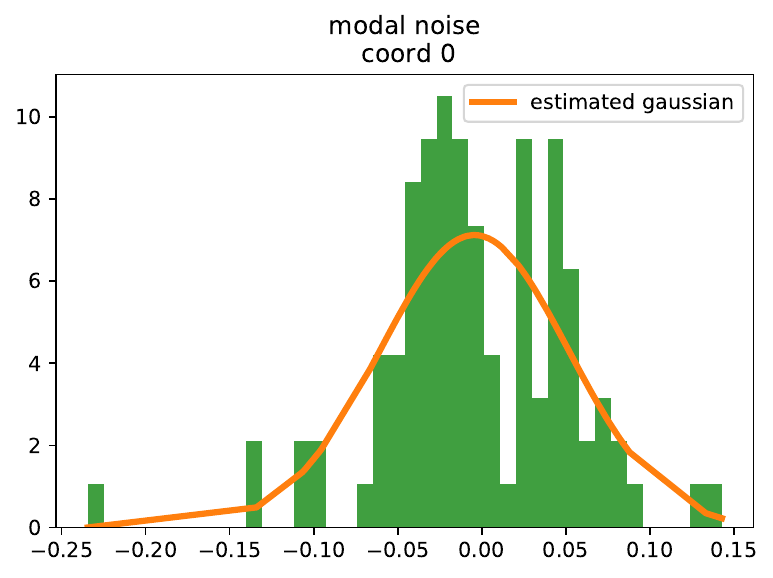}
         \caption{50 modes}
     \end{subfigure}
     \hfill
     \begin{subfigure}[b]{0.45\textwidth}
         \centering
         \includegraphics[width=\textwidth]{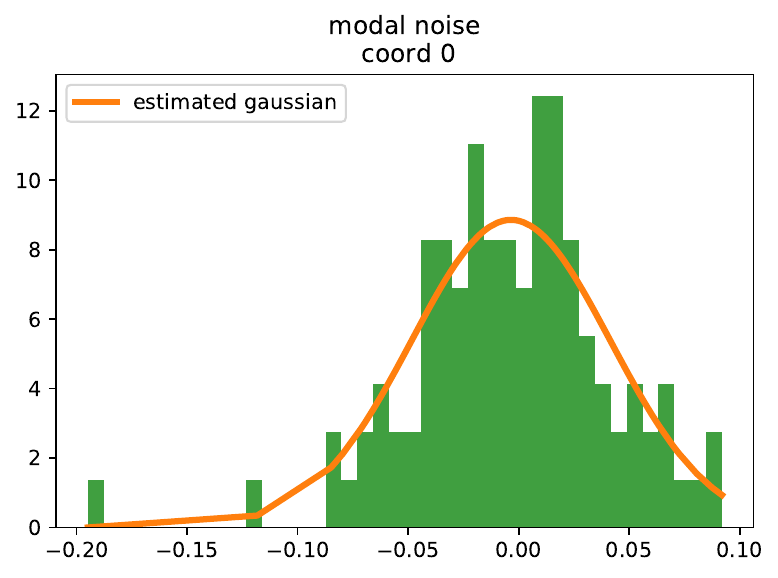}
         \caption{60 modes}
     \end{subfigure}
     \hfill
     \begin{subfigure}[b]{0.45\textwidth}
         \centering
         \includegraphics[width=\textwidth]{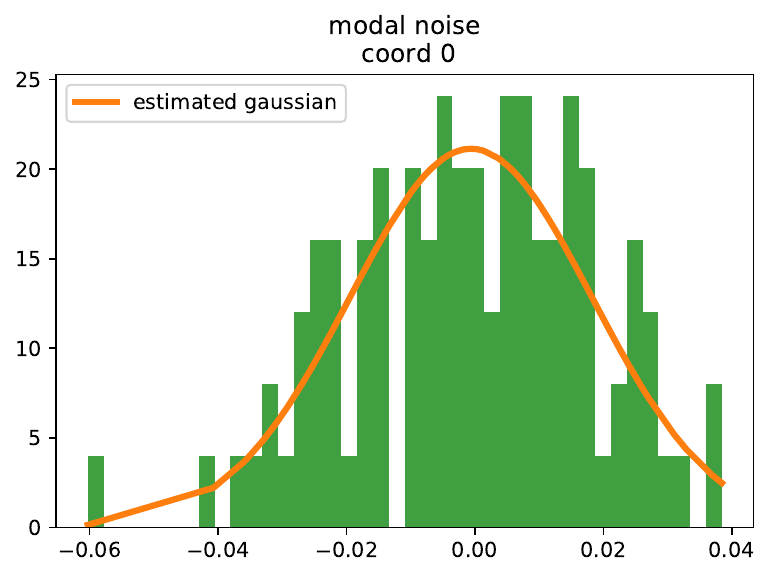}
         \caption{90 modes}
     \end{subfigure}
     \hfill
     \begin{subfigure}[b]{0.45\textwidth}
         \centering
         \includegraphics[width=\textwidth]{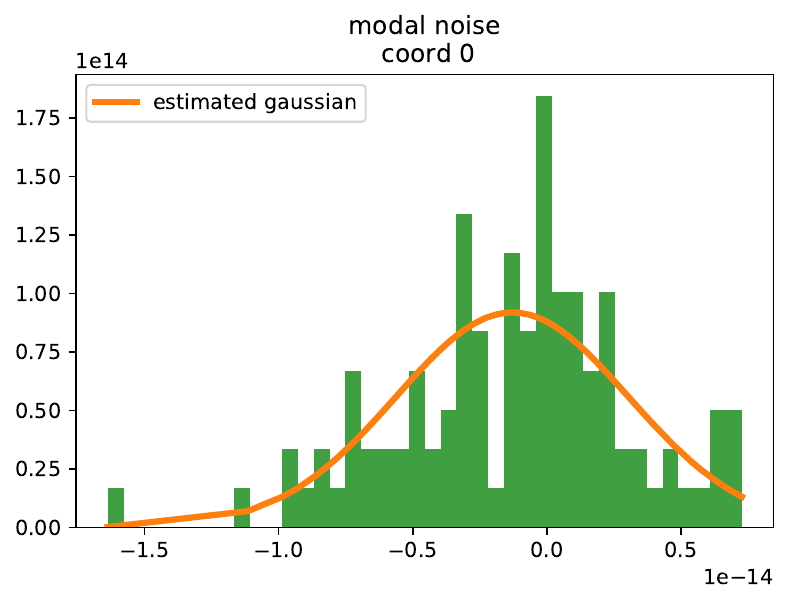}
         \caption{100 modes}
     \end{subfigure}
     \hfill
     \caption{\textbf{Anharmonic Oscillators}. The in-plane modal noise for oscillator 0. We estimate the Gaussian by computing the standard deviation of the modal noise distribution. The estimated standard deviation will be used to compute the confidence bounds for the reduced order models signal. Note the decreasing range of the scales on the x-axis as the number of modes increases. These results are typical for the other oscillators as well.}
\label{fig:anharmonic-modal-noise}
\end{figure}

\begin{figure}[ht]
	\centering
	\begin{subfigure}[b]{0.45\textwidth}
         \centering
         \includegraphics[width=\textwidth]{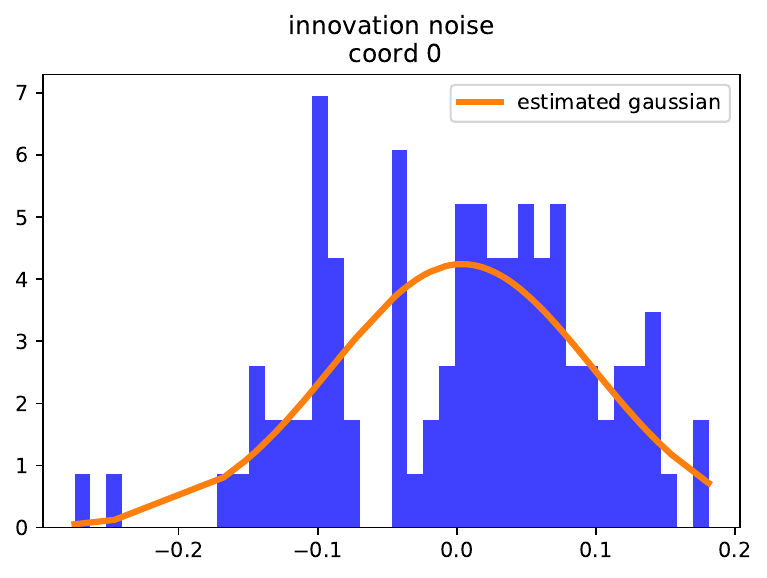}
         \caption{10 modes}
     \end{subfigure}
     \hfill
     \begin{subfigure}[b]{0.45\textwidth}
         \centering
         \includegraphics[width=\textwidth]{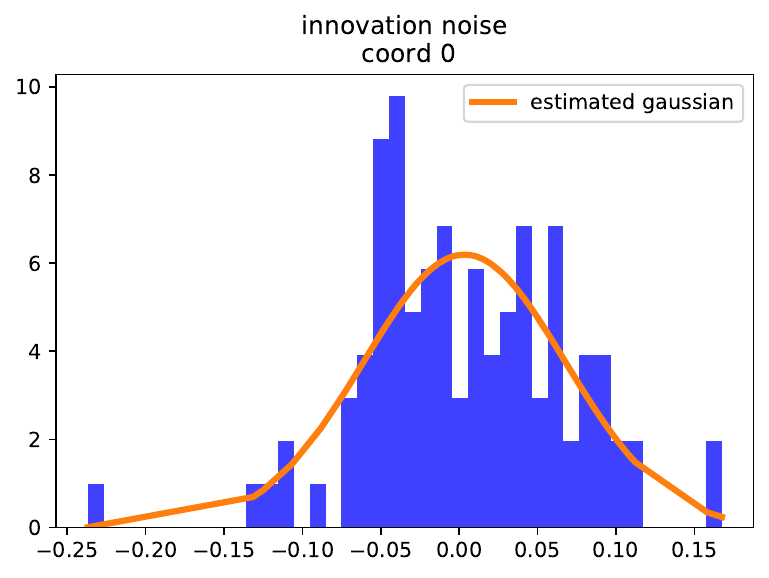}
         \caption{20 modes}
     \end{subfigure} 
     \hfill
     \begin{subfigure}[b]{0.45\textwidth}
         \centering
         \includegraphics[width=\textwidth]{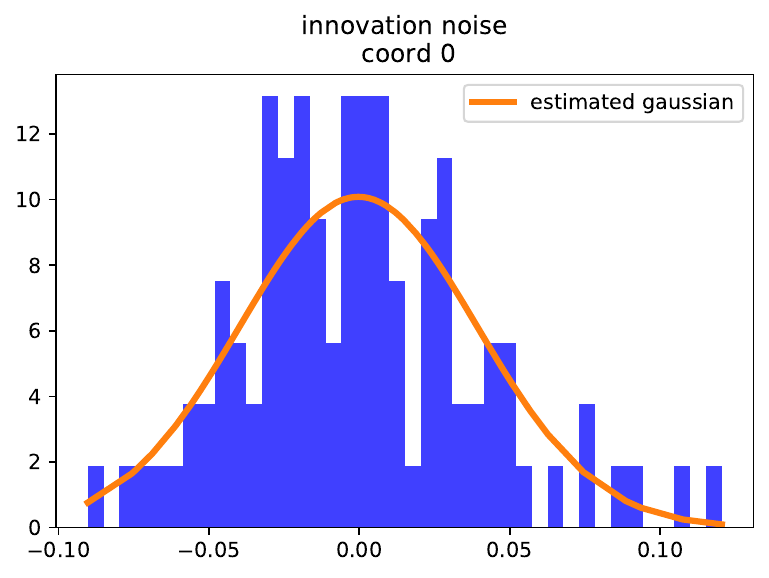}
         \caption{50 modes}
     \end{subfigure}
     \hfill
     \begin{subfigure}[b]{0.45\textwidth}
         \centering
         \includegraphics[width=\textwidth]{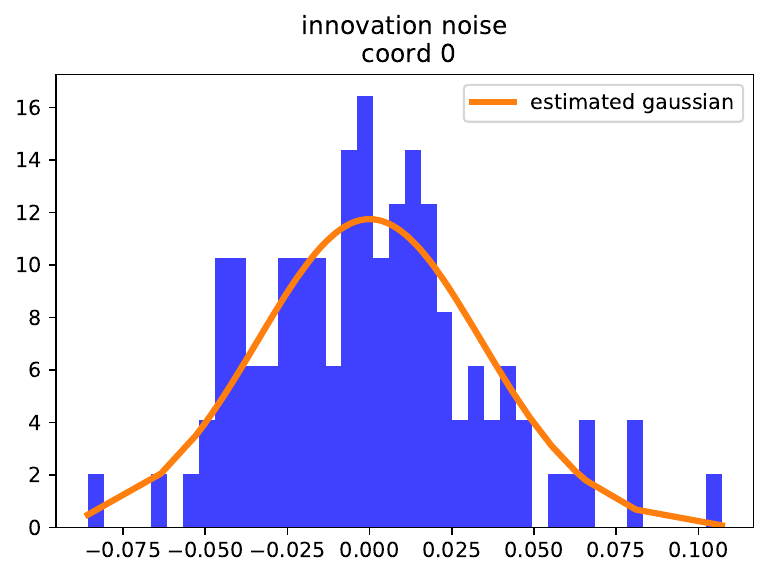}
         \caption{60 modes}
     \end{subfigure}
     \hfill
     \begin{subfigure}[b]{0.45\textwidth}
         \centering
         \includegraphics[width=\textwidth]{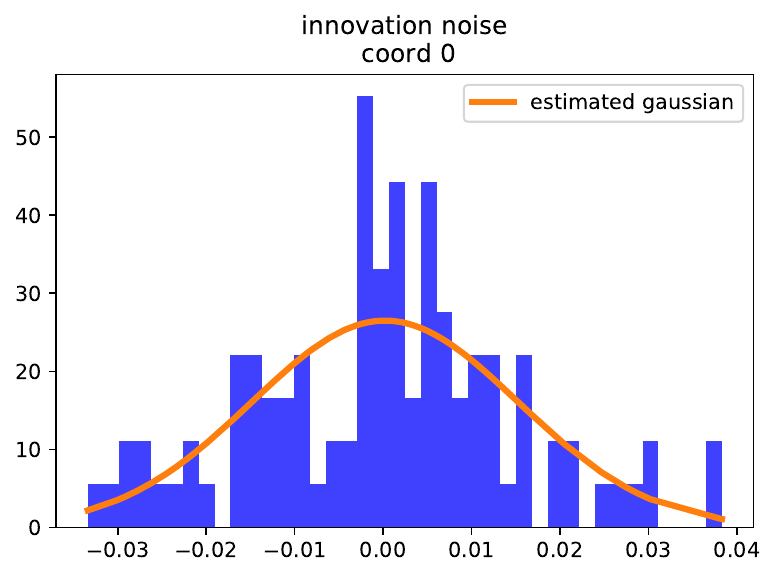}
         \caption{90 modes}
     \end{subfigure}
     \hfill
     \begin{subfigure}[b]{0.45\textwidth}
         \centering
         \includegraphics[width=\textwidth]{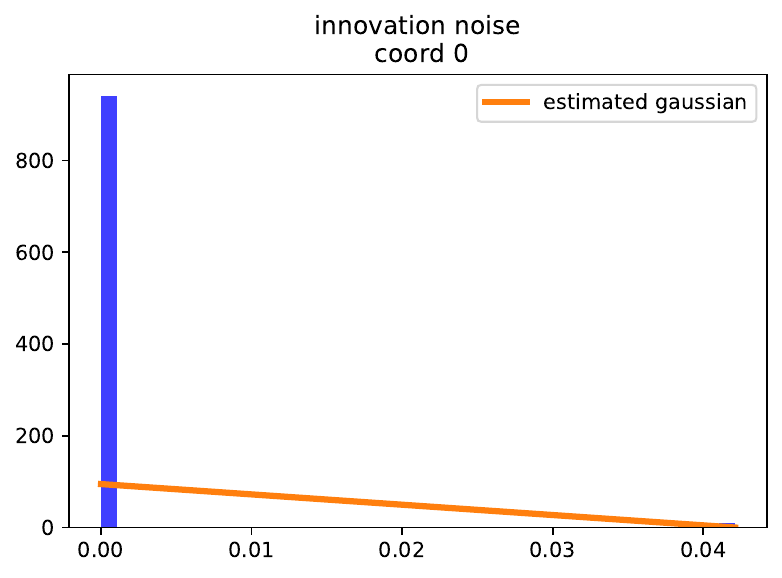}
         \caption{100 modes}
     \end{subfigure}
     \hfill
     \caption{\textbf{Anharmonic Oscillators}. The out-of-plane innovation noise for the for oscillator 0. These results are typical for the other oscillators as well.}
\label{fig:anharmonic-innovation-noise}
\end{figure}

\begin{figure}[ht]
	\centering
	\begin{subfigure}[b]{0.45\textwidth}
         \centering
         \includegraphics[width=\textwidth]{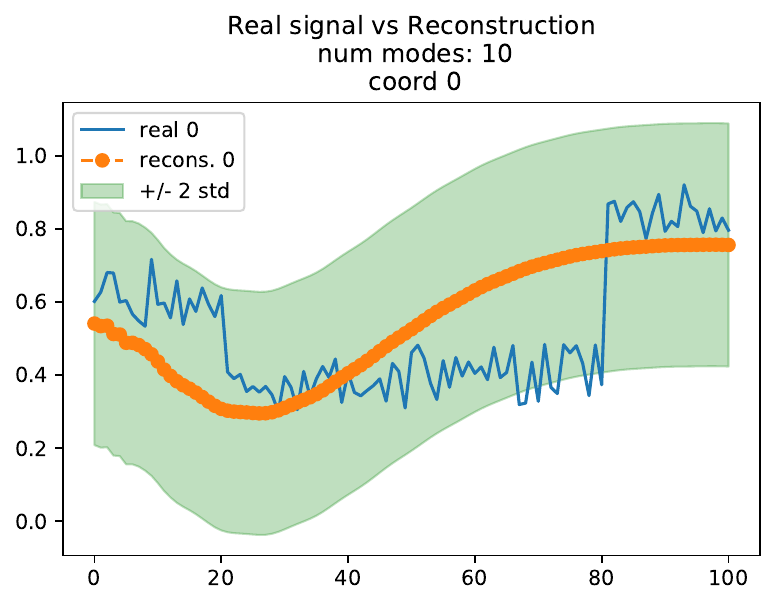}
         \caption{10 modes}
     \end{subfigure}
     \hfill
     \begin{subfigure}[b]{0.45\textwidth}
         \centering
         \includegraphics[width=\textwidth]{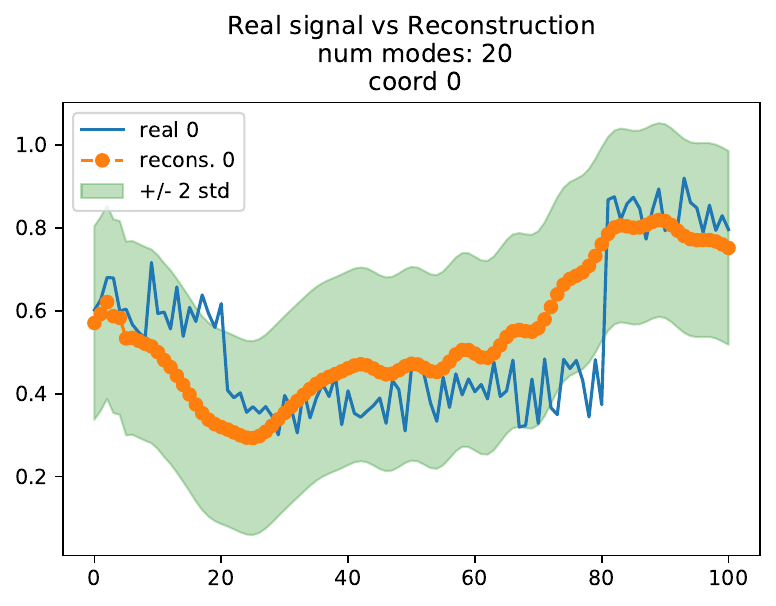}
         \caption{20 modes}
     \end{subfigure} 
     \hfill
     \begin{subfigure}[b]{0.45\textwidth}
         \centering
         \includegraphics[width=\textwidth]{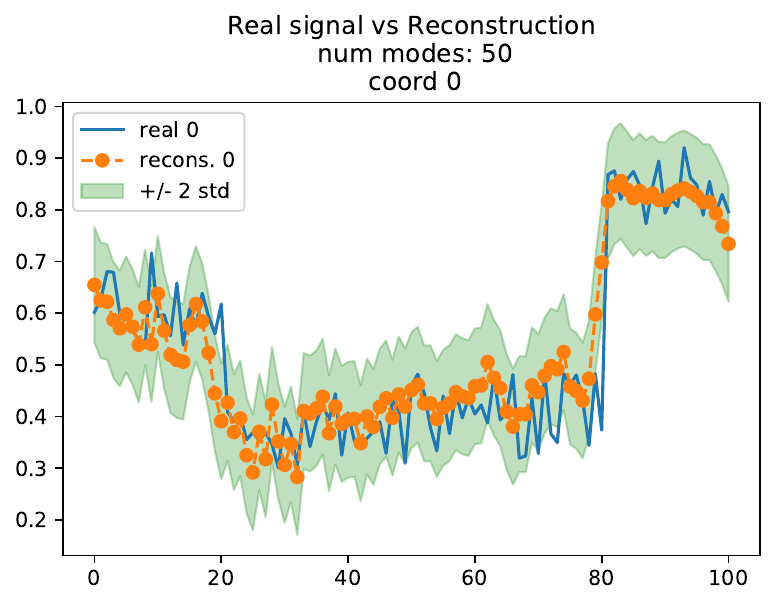}
         \caption{50 modes}
     \end{subfigure}
     \hfill
     \begin{subfigure}[b]{0.45\textwidth}
         \centering
         \includegraphics[width=\textwidth]{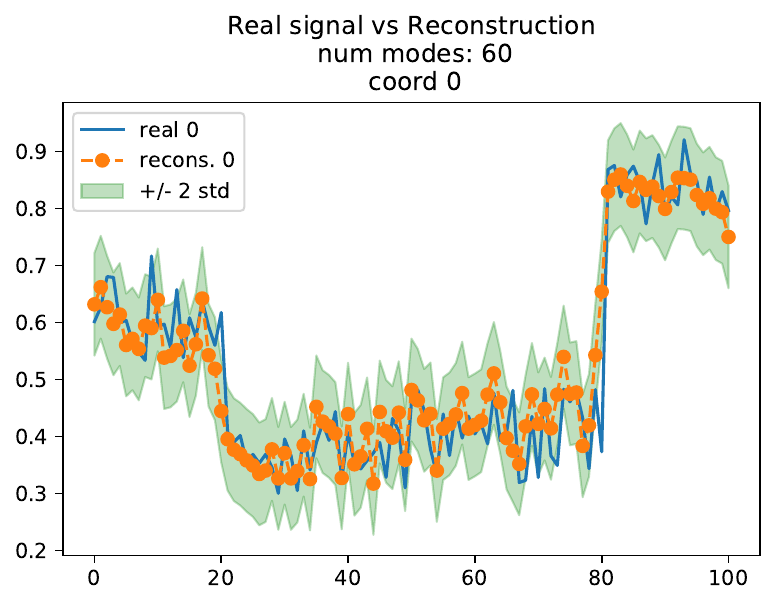}
         \caption{60 modes}
     \end{subfigure}
     \hfill
     \begin{subfigure}[b]{0.45\textwidth}
         \centering
         \includegraphics[width=\textwidth]{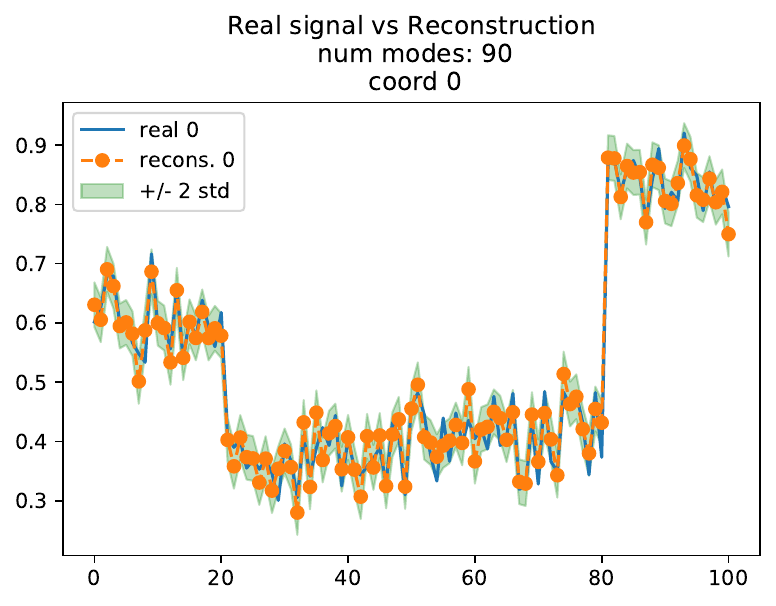}
         \caption{90 modes}
     \end{subfigure}
     \hfill
     \begin{subfigure}[b]{0.45\textwidth}
         \centering
         \includegraphics[width=\textwidth]{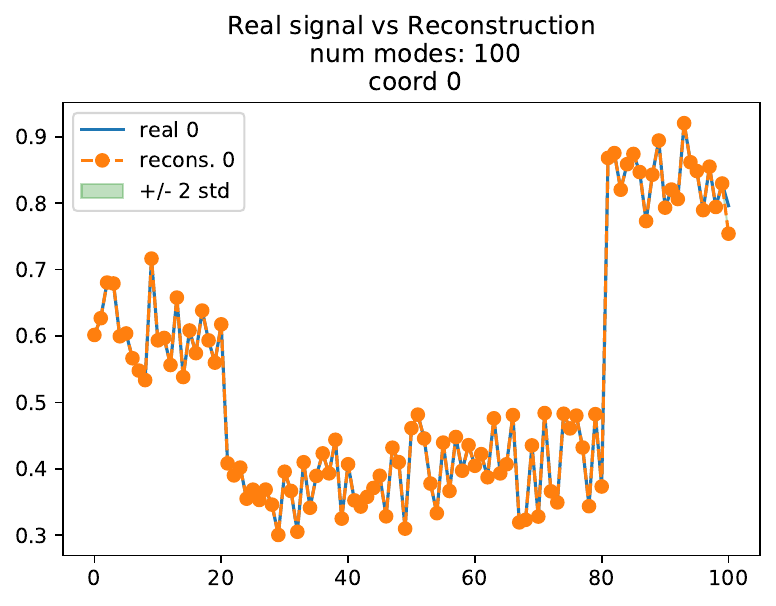}
         \caption{100 modes}
     \end{subfigure}
     \hfill
     \caption{\textbf{Anharmonic Oscillators}. Comparison between the reduced order model (ROM) signal (orange) with the real signal (blue) for the action variable of oscillator 0. The 95\% confidence interval is given by the green band around the orange signal. As the number of reconstruction modes increases, the ROM signal follows the true signal better and the 95\% confidence becomes tighter. These results are typical for the other oscillators as well}
\label{fig:anharmonic-signal-reconstruction-action}
\end{figure}

\begin{figure}[ht]
	\centering
	\begin{subfigure}[b]{0.45\textwidth}
         \centering
         \includegraphics[width=\textwidth]{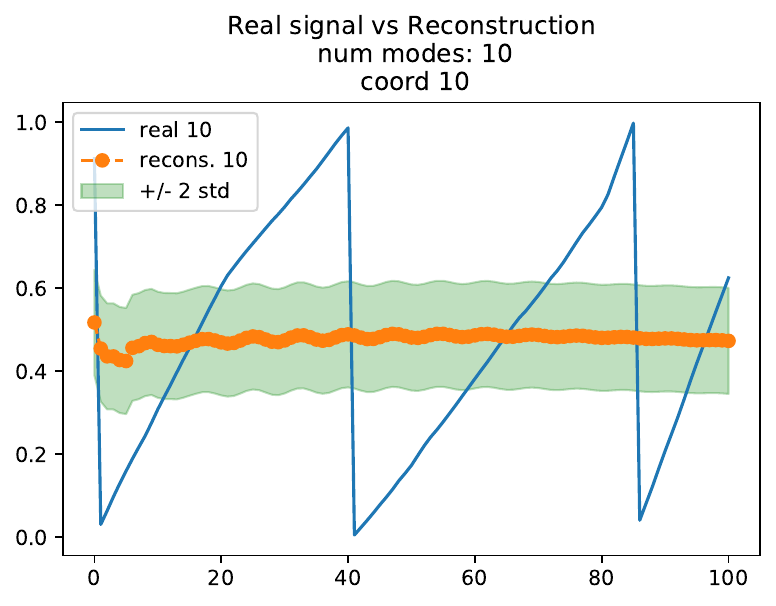}
         \caption{10 modes}
     \end{subfigure}
     \hfill
     \begin{subfigure}[b]{0.45\textwidth}
         \centering
         \includegraphics[width=\textwidth]{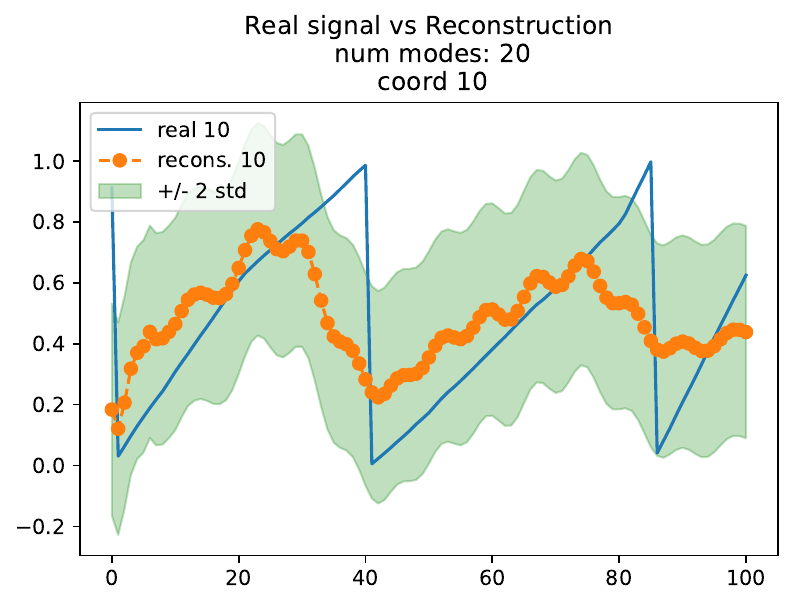}
         \caption{20 modes}
     \end{subfigure} 
     \hfill
     \begin{subfigure}[b]{0.45\textwidth}
         \centering
         \includegraphics[width=\textwidth]{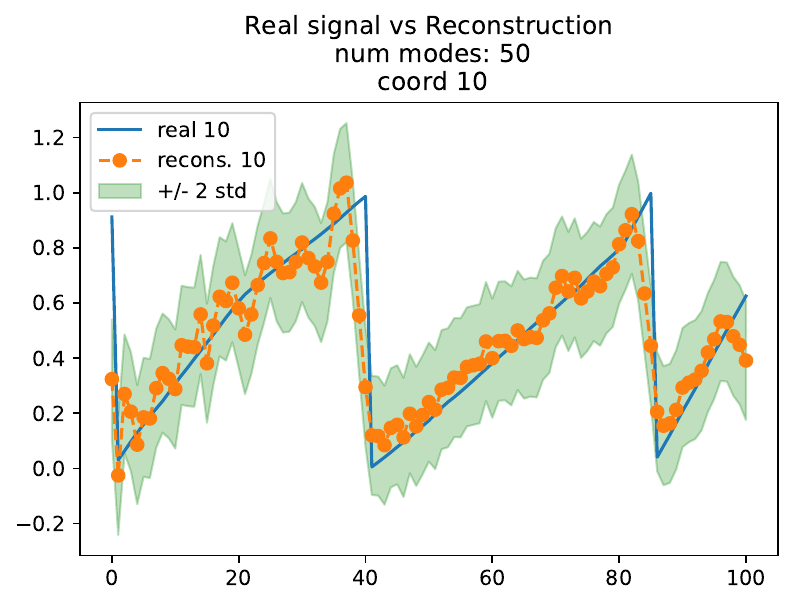}
         \caption{50 modes}
     \end{subfigure}
     \hfill
     \begin{subfigure}[b]{0.45\textwidth}
         \centering
         \includegraphics[width=\textwidth]{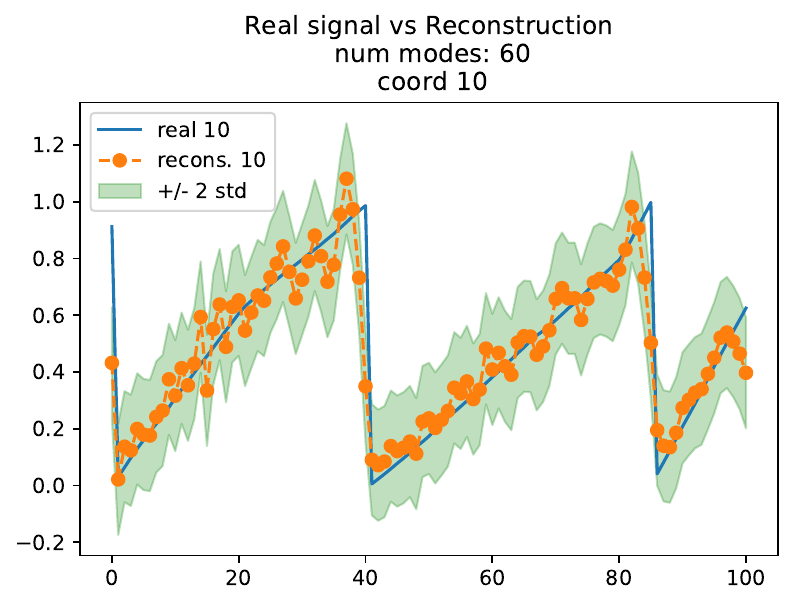}
         \caption{60 modes}
     \end{subfigure}
     \hfill
     \begin{subfigure}[b]{0.45\textwidth}
         \centering
         \includegraphics[width=\textwidth]{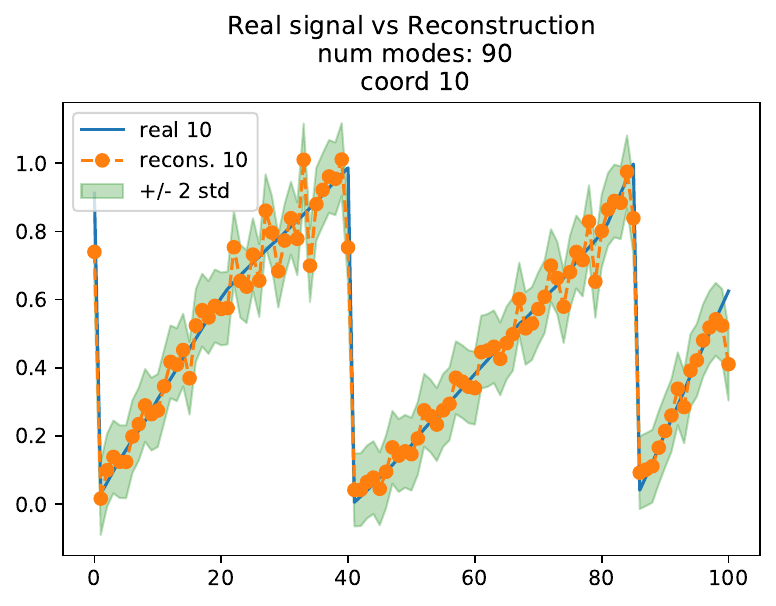}
         \caption{90 modes}
     \end{subfigure}
     \hfill
     \begin{subfigure}[b]{0.45\textwidth}
         \centering
         \includegraphics[width=\textwidth]{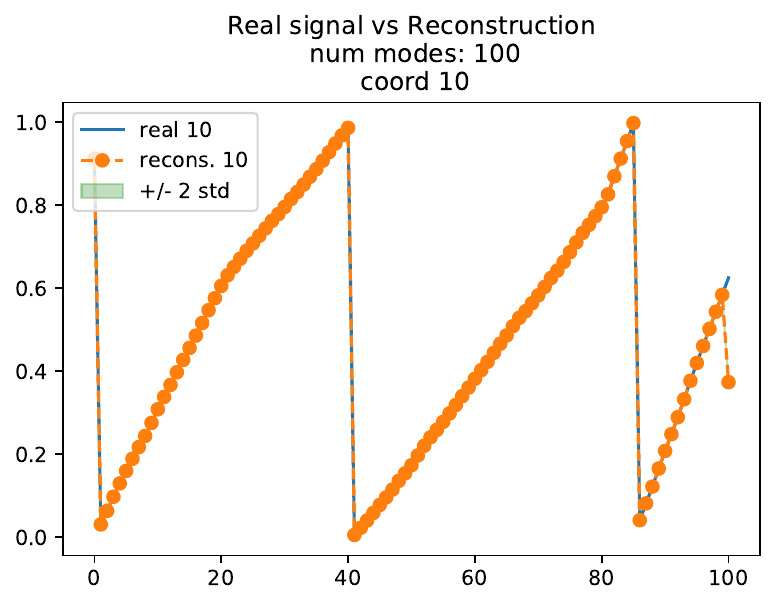}
         \caption{100 modes}
     \end{subfigure}
     \hfill
     \caption{\textbf{Anharmonic Oscillators}. Comparison between the reduced order model (ROM) signal (orange) with the real signal (blue) for the \textbf{angle variable} of oscillator 0. The 95\% confidence interval is given by the green band around the orange signal. As the number of reconstruction modes increases, the ROM signal follows the true signal better and the 95\% confidence becomes tighter. These results are typical for the other oscillators as well}
\label{fig:anharmonic-signal-reconstruction-angle}
\end{figure}

\begin{figure}[ht]
	\centering
	\begin{subfigure}[b]{0.45\textwidth}
         \centering
         \includegraphics[width=\textwidth]{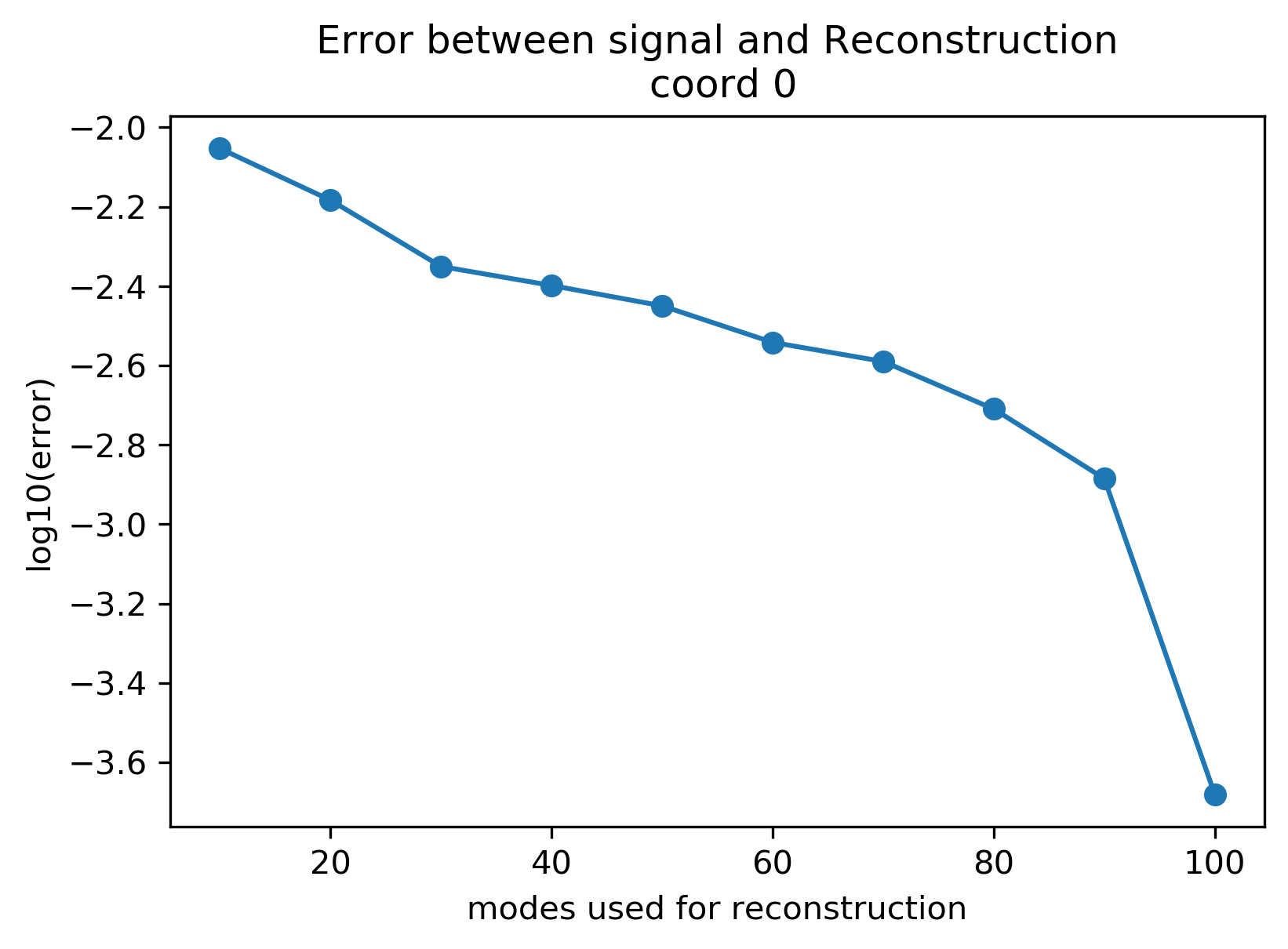}
         \caption{Oscillator 0 action variable.}
     \end{subfigure}
     \hfill
     \begin{subfigure}[b]{0.45\textwidth}
         \centering
         \includegraphics[width=\textwidth]{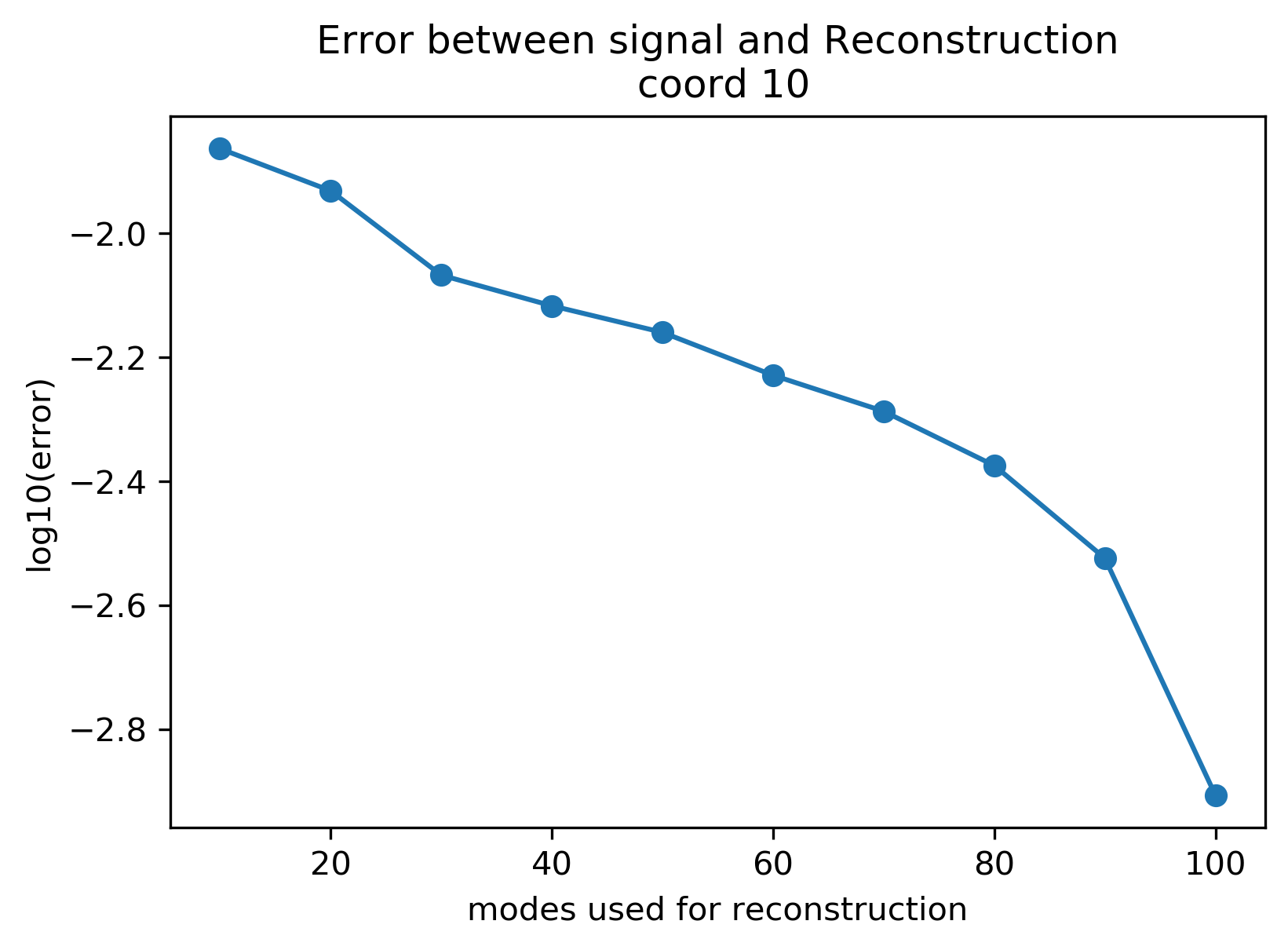}
         \caption{Oscillator 0 angle variable.}
     \end{subfigure} 
     \hfill
     \begin{subfigure}[b]{0.45\textwidth}
         \centering
         \includegraphics[width=\textwidth]{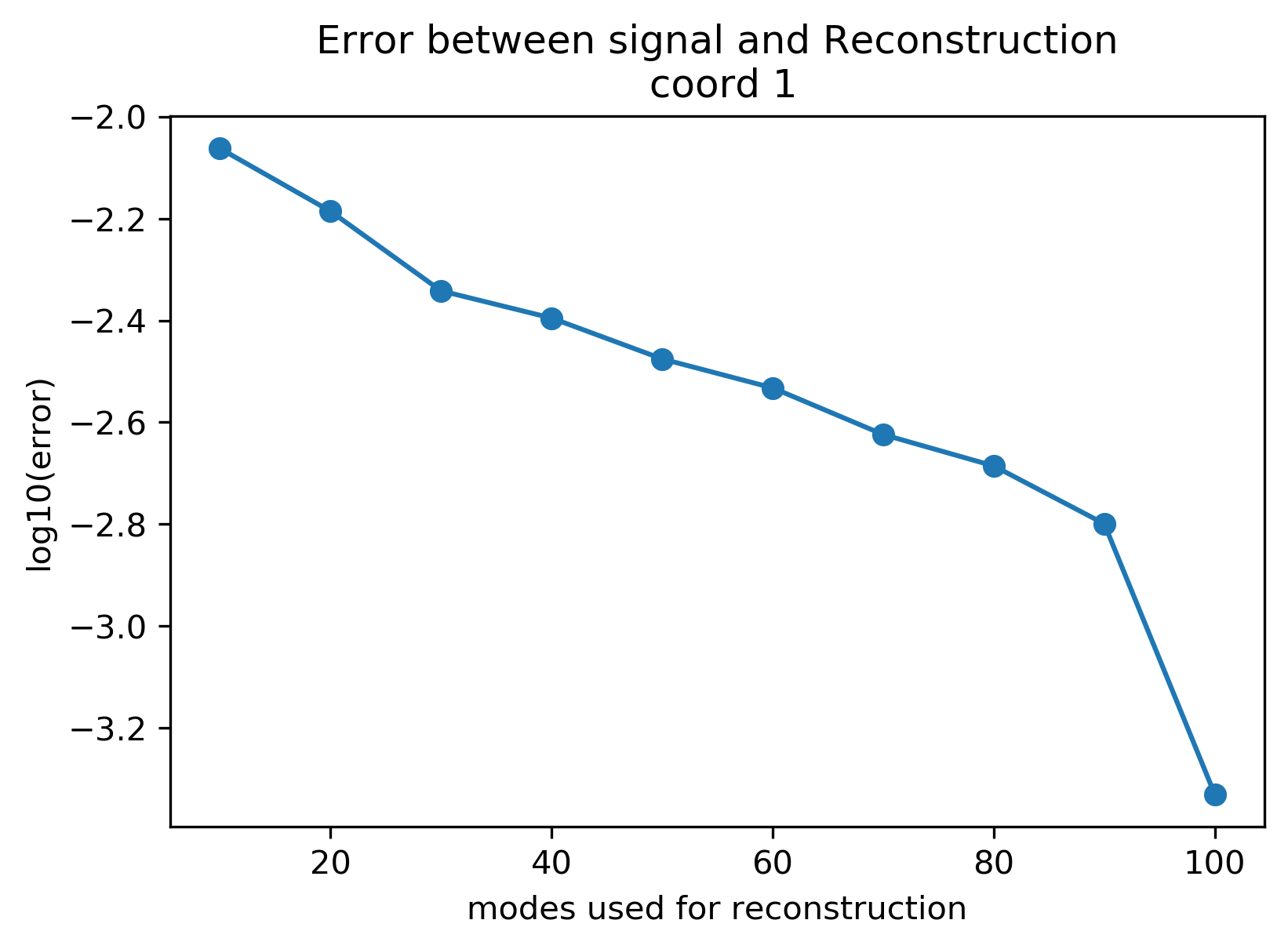}
         \caption{Oscillator 1 action variable.}
     \end{subfigure}
     \hfill
     \begin{subfigure}[b]{0.45\textwidth}
         \centering
         \includegraphics[width=\textwidth]{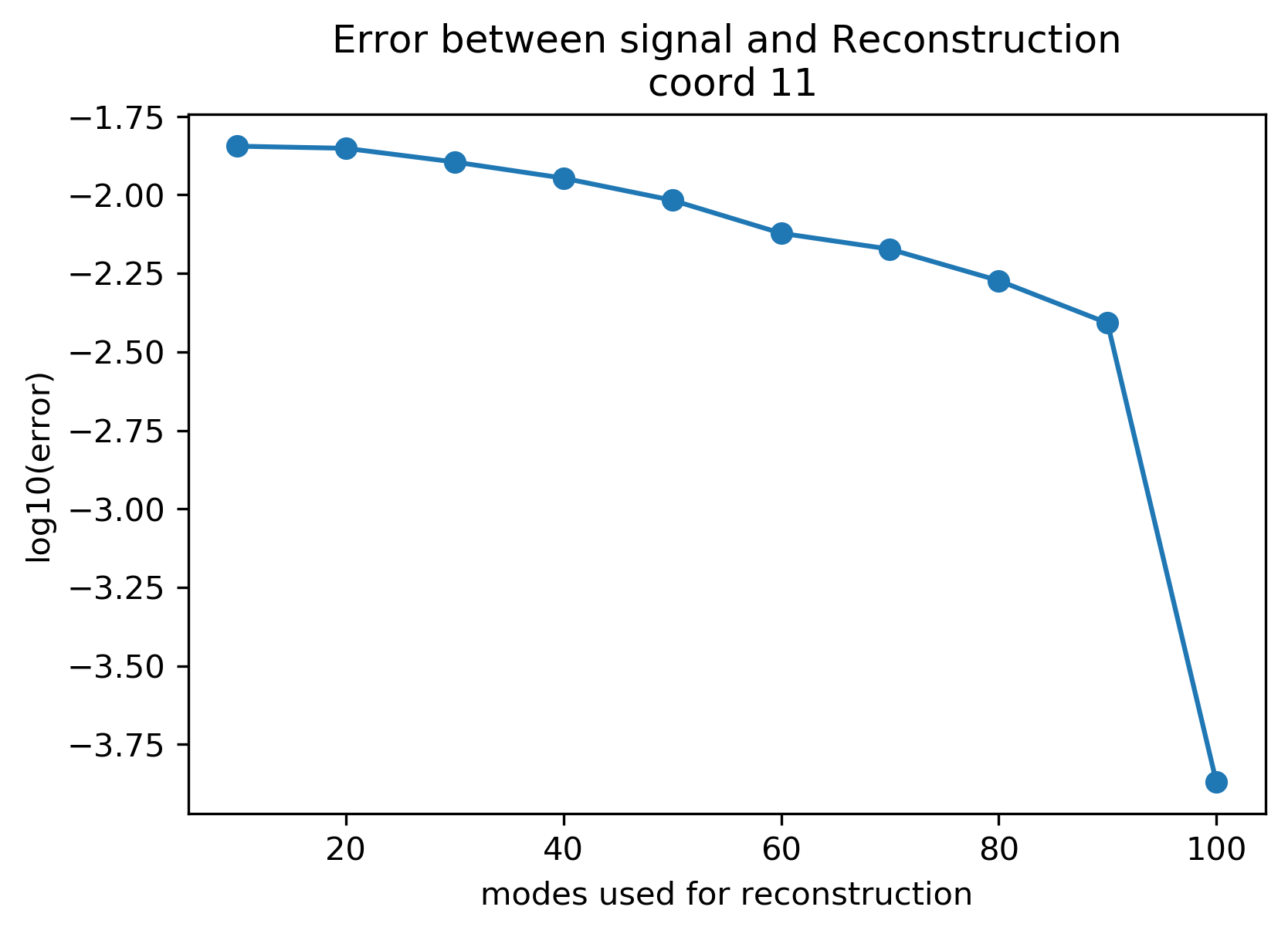}
         \caption{Oscillator 1 angle variable.}
     \end{subfigure}
     \hfill
     \begin{subfigure}[b]{0.45\textwidth}
         \centering
         \includegraphics[width=\textwidth]{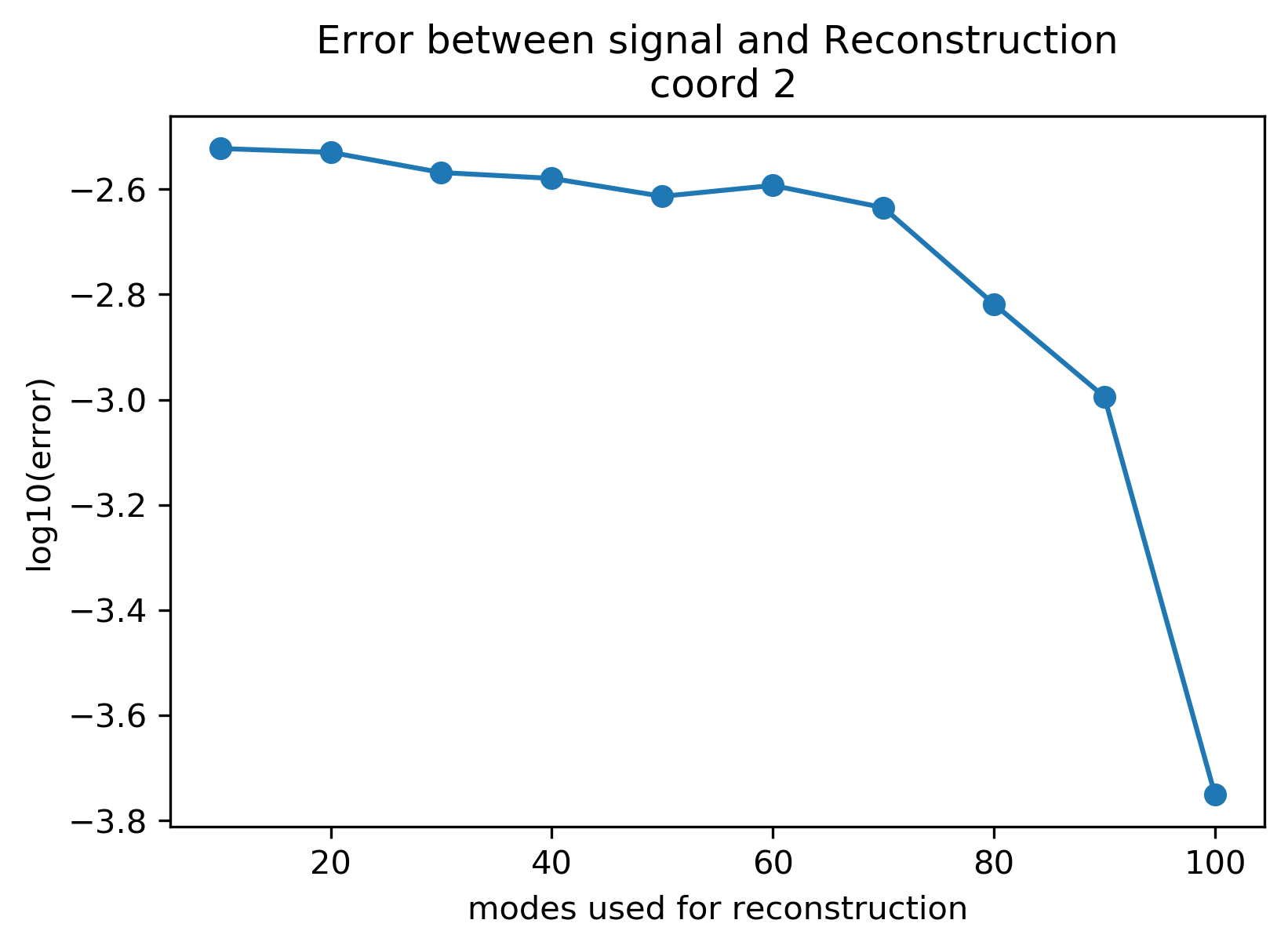}
         \caption{Oscillator 2 action variable.}
     \end{subfigure}
     \hfill
     \begin{subfigure}[b]{0.45\textwidth}
         \centering
         \includegraphics[width=\textwidth]{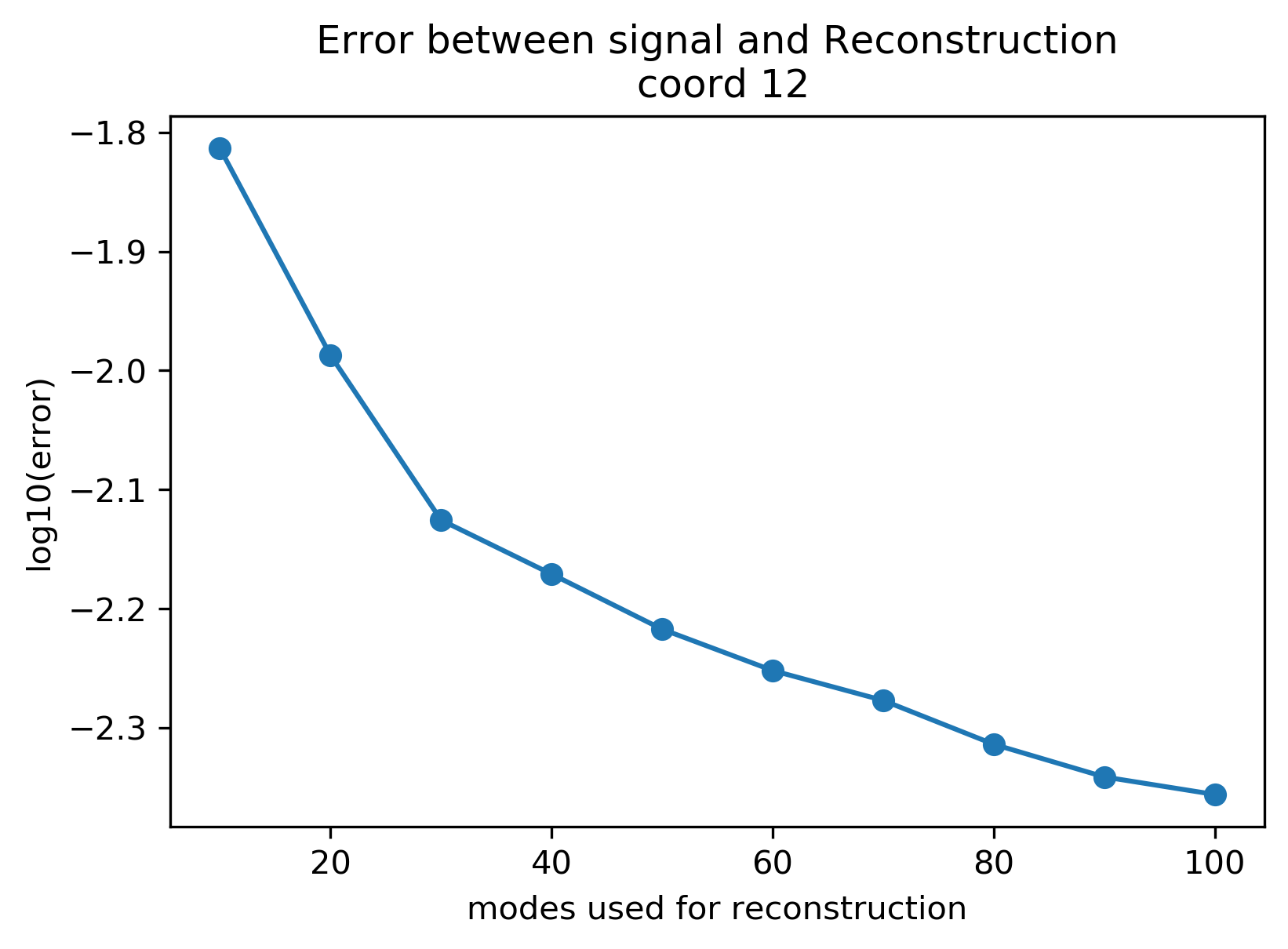}
         \caption{Oscillator 2 action variable.}
     \end{subfigure}
     \hfill
     \caption{\textbf{Anharmonic Oscillators}. Reconstruction errors of action and angle variables vs. number of modes used for the reconstruction for 3 of the oscillators.}
\label{fig:anharmonic-error}
\end{figure}

\begin{figure}[ht]
	\centering
	\begin{subfigure}[b]{0.45\textwidth}
         \centering
         \includegraphics[width=\textwidth]{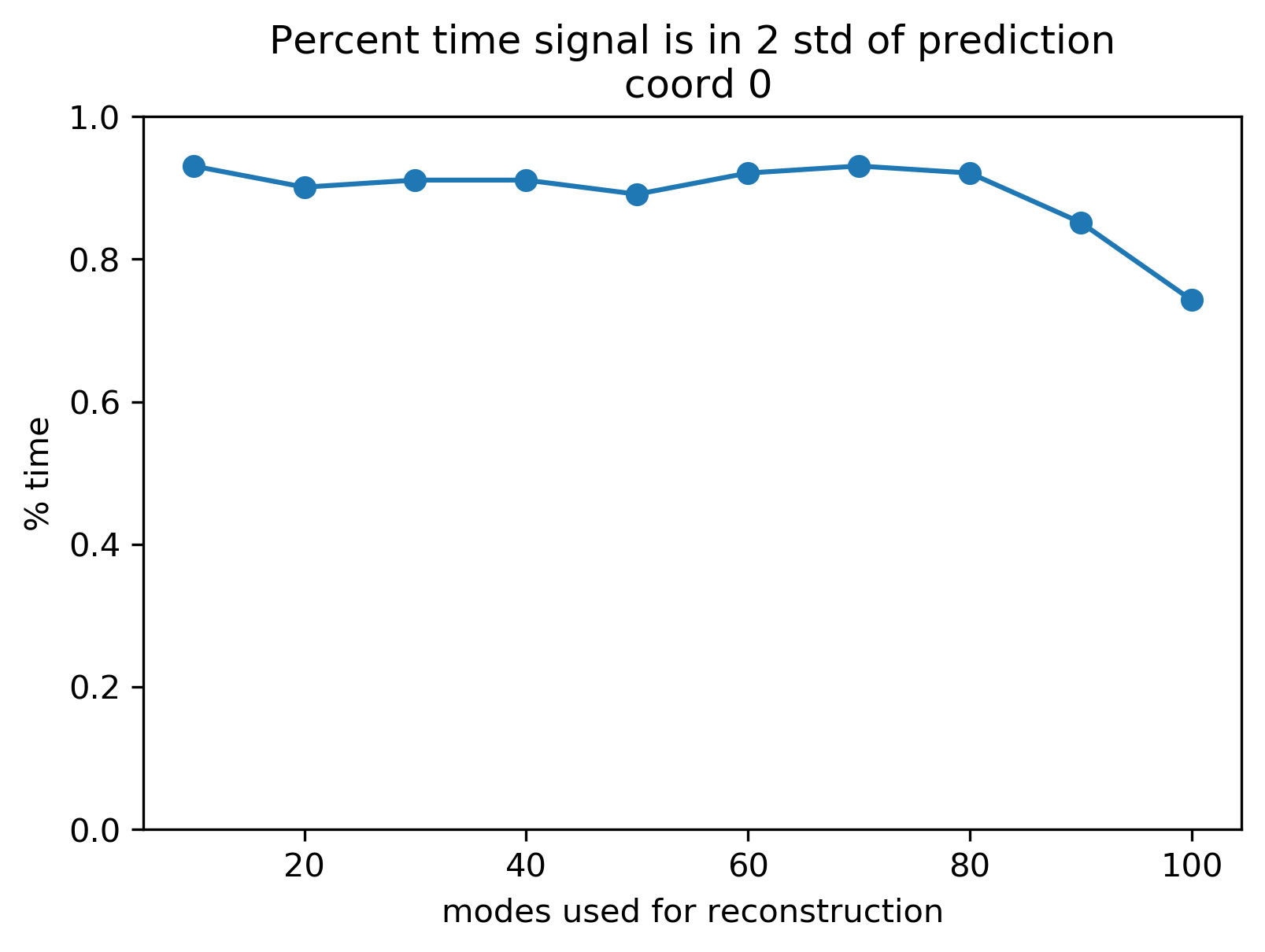}
         \caption{Oscillator 0 action variable.}
     \end{subfigure}
     \hfill
     \begin{subfigure}[b]{0.45\textwidth}
         \centering
         \includegraphics[width=\textwidth]{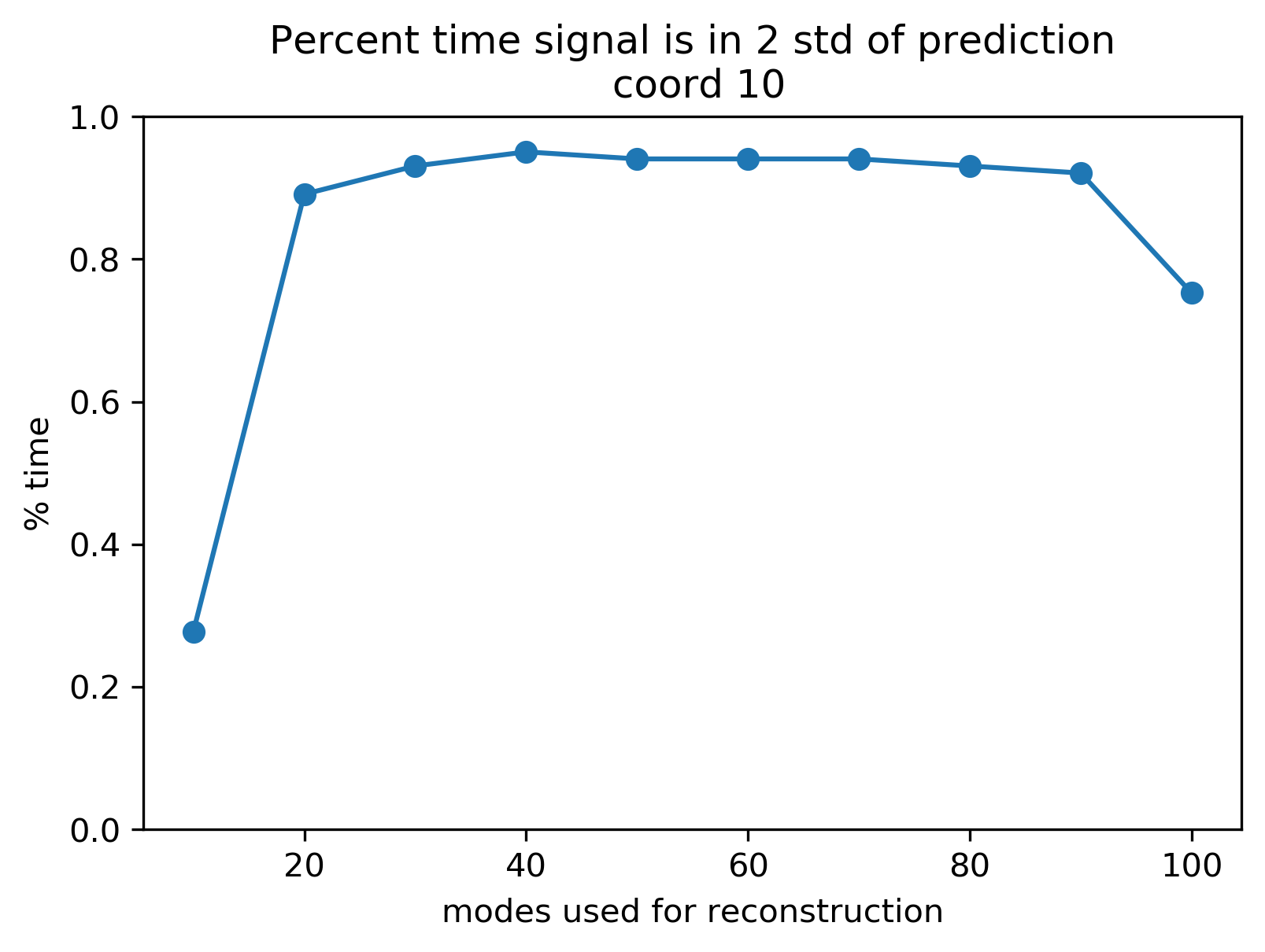}
         \caption{Oscillator 0 angle variable.}
     \end{subfigure} 
     \hfill
     \begin{subfigure}[b]{0.45\textwidth}
         \centering
         \includegraphics[width=\textwidth]{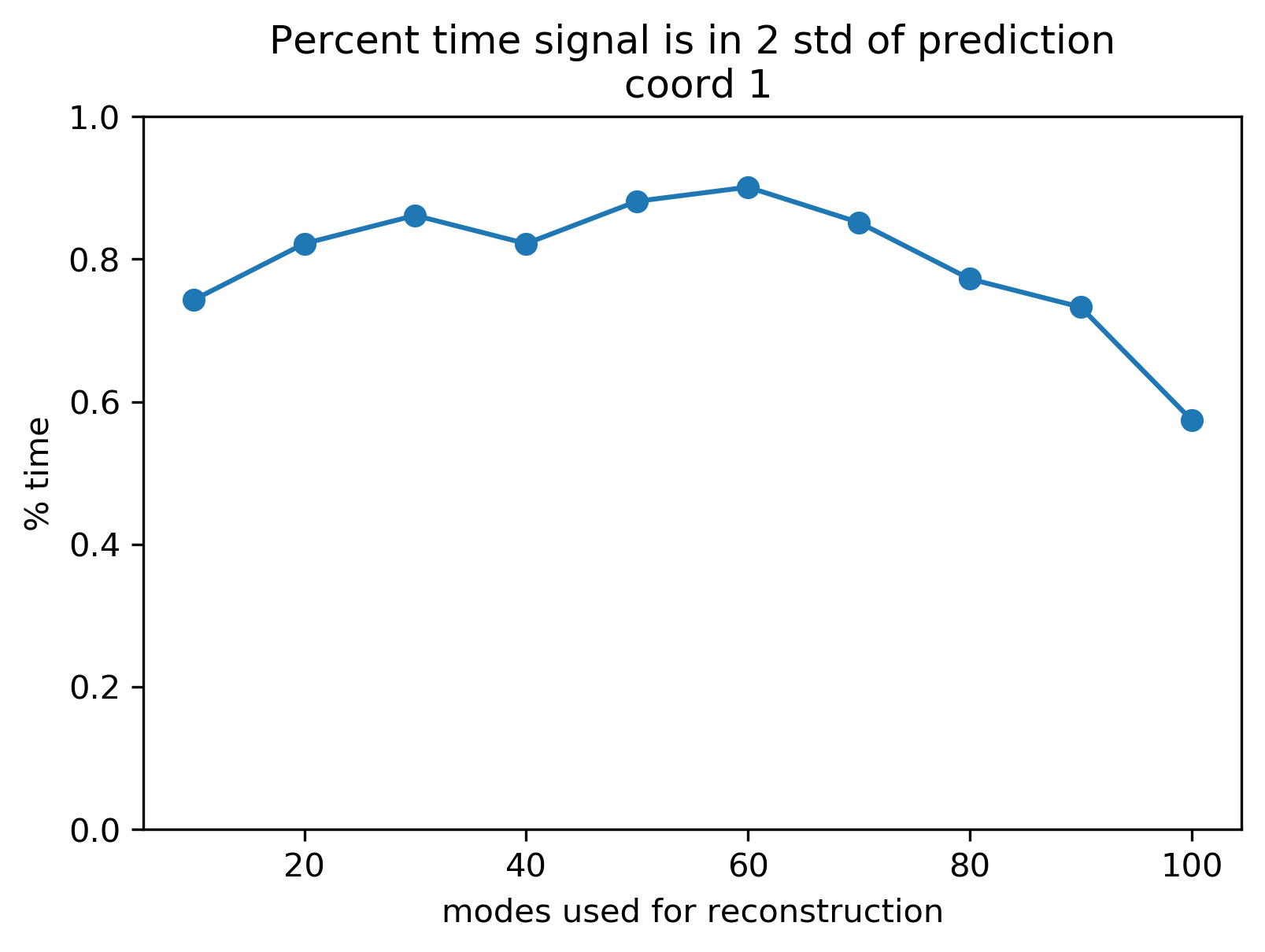}
         \caption{Oscillator 1 action variable.}
     \end{subfigure}
     \hfill
     \begin{subfigure}[b]{0.45\textwidth}
         \centering
         \includegraphics[width=\textwidth]{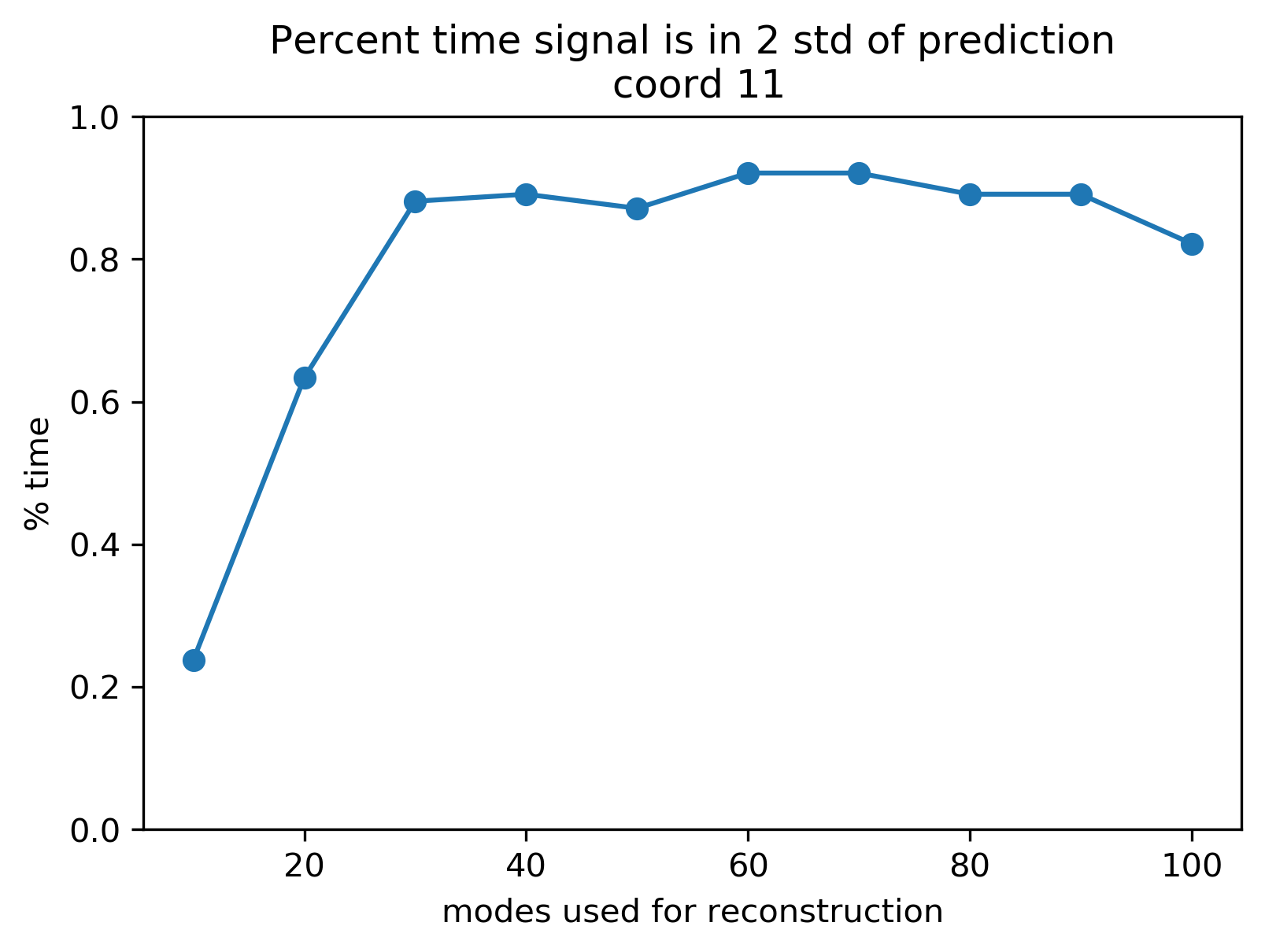}
         \caption{Oscillator 1 angle variable.}
     \end{subfigure}
     \hfill
     \begin{subfigure}[b]{0.45\textwidth}
         \centering
         \includegraphics[width=\textwidth]{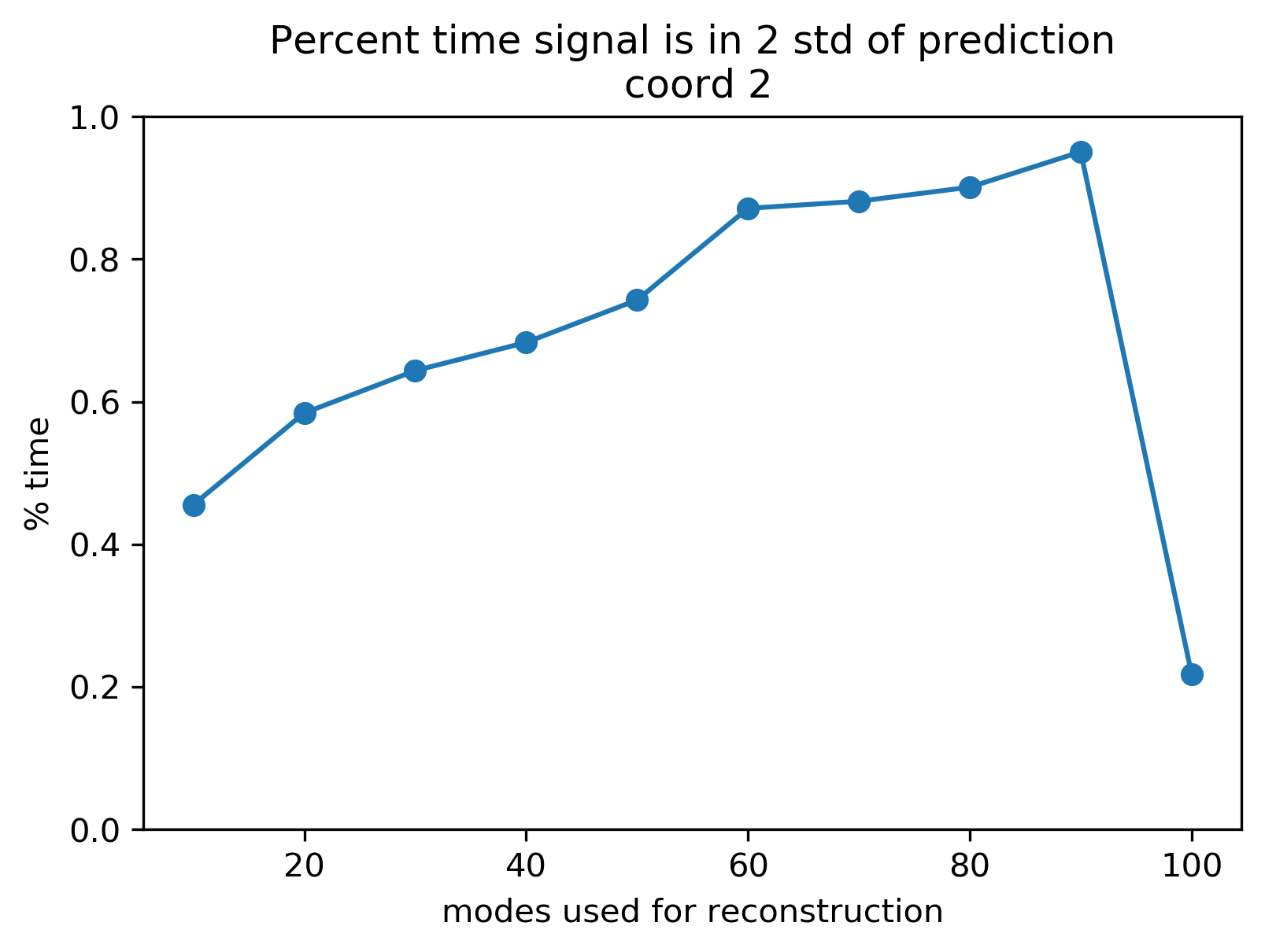}
         \caption{Oscillator 2 action variable.}
     \end{subfigure}
     \hfill
     \begin{subfigure}[b]{0.45\textwidth}
         \centering
         \includegraphics[width=\textwidth]{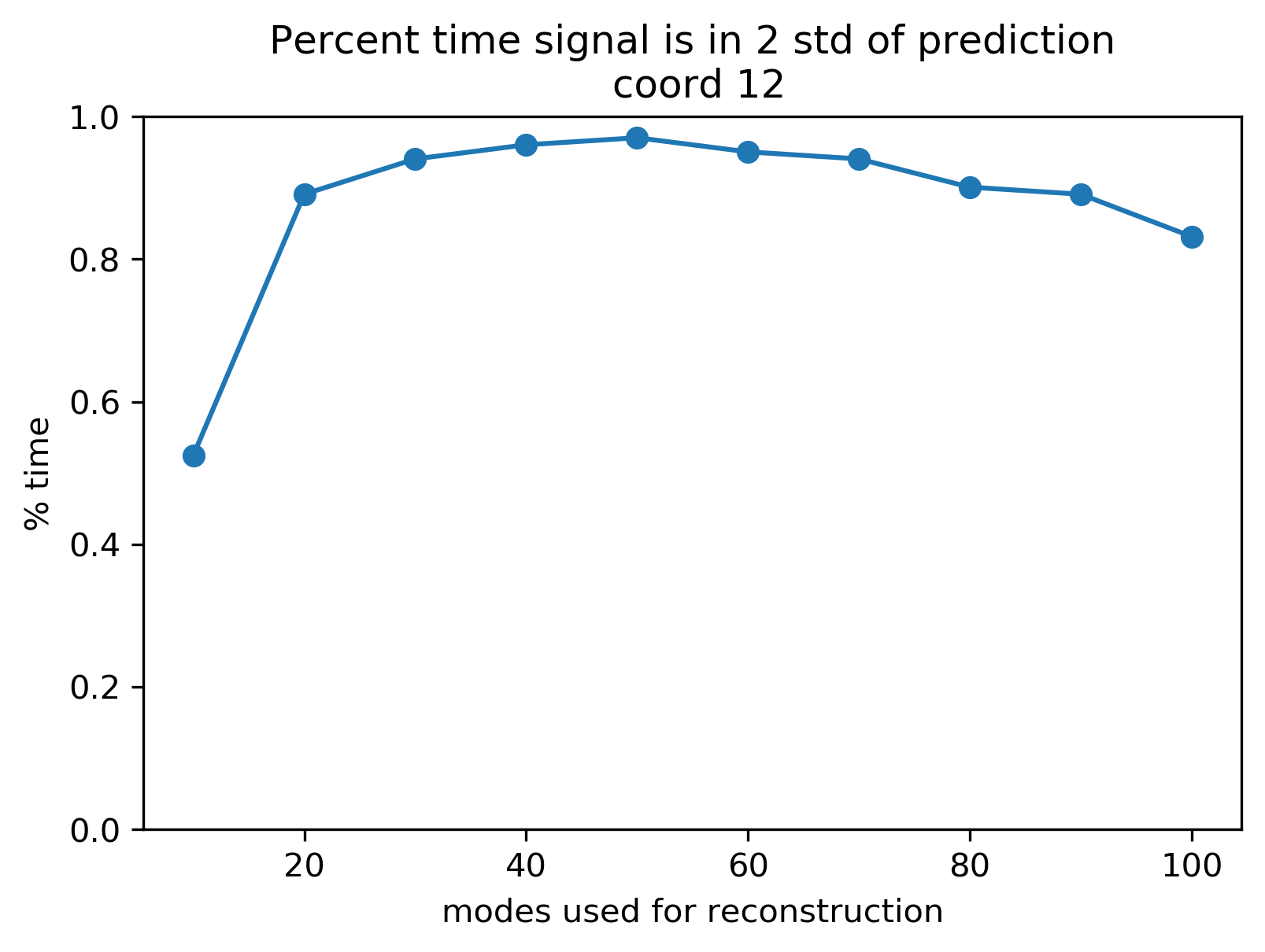}
         \caption{Oscillator 2 angle variable.}
     \end{subfigure}
     \hfill
     \caption{\textbf{Anharmonic Oscillators}. Fraction of time that the true signal is within the 95\% confidence interval. The large drop in the residence times when using 100 modes is due to the relative error between the signal and the standard deviation of the modal noise. When using a large number of modes for the ROM, the noise gets subsumed into the deterministic model. When using the full number of modes (100), we see in figure \ref{fig:anharmonic-modal-noise} that the computed modal noise's standard deviation is on the order of $10^{-14}$, meaning that the noise has been fully subsumed into the deterministic dynamics. The error between the ROM and the signal is on the order of $10^{-2}$ to $10^{-4}$ and thus adding a $10^{-14}$ confidence bound rarely contains the true signal.}
\label{fig:anharmonic-residence-time}
\end{figure}

\blue{
Figures \ref{fig:anharmonic-eigenvalues} - \ref{fig:anharmonic-residence-time} show an ablative study where a family of ROMs are constructed using different subsets of the modes. It is notable from figures \ref{fig:anharmonic-signal-reconstruction-action} and \ref{fig:anharmonic-signal-reconstruction-angle} that the signal reconstruction gets better when around 50 modes are used. It can also be seen from figure \ref{fig:anharmonic-modal-noise} that with a small number of retained modes the modal noise distribution is not Gaussian, becomes more Gaussian as more modes are retained, until it starts approaching a delta functions as too many modes are retained.

The heuristic in section \ref{sec:heuristic} gives a method of choosing the number of modes retained based on a statistical test that shows how close the modal distribution is to a normal distribution.

The heuristic dictates that 40 modes (in fact, that number is somewhat lower (see figure \ref{fig:heuristic-anharmonic}) is the minimum number of modes that can be used for the model, as this is where the modal noise becomes approximately Gaussian.
}

Figure \ref{fig:anharmonic-error} shows the error of the ROM's prediction vs. the number of modes used in the model. The error of the reconstruction does not, however, reflect the accuracy of the ROM prediction. Namely, the ROM has been computed from a noisy signal, and thus the representation with a full number of modes subsumes noise into coherent modal dynamics. Residence times for the different ROMs are given in figure \ref{fig:anharmonic-residence-time}.

\clearpage
\subsection{Kuramoto Model}

As our final example for applying the method, we turn to the Kuramoto model of coupled oscillators \cite{acebron2005kuramoto}. In the classical model (see eq.\ref{eq:kuramoto}), the strength of the connections between the oscillators is fixed and uniform for all oscillators. As the strength of the connection increases, the oscillators can display behaviors such as synchronization \cite{acebron2005kuramoto}.
	\begin{equation}\label{eq:kuramoto}
	\frac{d\theta_i}{dt} = \omega_i + \xi_i + \frac{K}{N} \sum_{j=1}^{N} \sin(\theta_j - \theta_i).
	\end{equation}

Here, we introduce a modification where instead of the fixed coupling strength $K/N$ we introduce randomness. Our enhanced model looks like
	\begin{equation}\label{eq:kuramoto-random-strength}
	\frac{d\theta_i}{dt} = \omega_i + \xi_i +  \sum_{j=1}^{N} \frac{K + \zeta_{i,j}}{N} \sin(\theta_j - \theta_i).
	\end{equation}
In our examples, each $\zeta_{i,j}$ is distributed according to $N(0, \sigma^2)$, with $\sigma = 1$, $K=5$, and  $N=10$. The $N\times N$ coupling matrix $Z = [\zeta_{i,j}]_{1\leq i,j \leq N}$ is restricted to be symmetric with $\zeta_{i,j} = \zeta_{j,i}$. In our experiments, we have chosen to set $\xi_i \sim N(0, 0.25)$ for all $i$. The natural frequencies were drawn uniformly at random from the interval $[0.25, 0.75]$, whereas the initial condition for each coordinate of $\theta$ was chosen uniformly at random from the interval $[0, 2\pi]$. When constructing the data matrix we use a complex representation of the angles ($\theta \mapsto z = \exp(i \theta)$).

\begin{table}[h]
\caption{Fixed parameters for Kuramoto model simulation }
\begin{center}
\begin{tabular}{|ll|}
\hline
number of oscillators, $N$ & 10 \\ 
\hline
numerical time step, $dt$ & 0.05 sec \\
\hline
simulation time, $t$ & 20 sec \\
\hline
Nominal coupling strength, $K$ & 5 \\
\hline 
Effective coupling strength, $K/N$ & 0.5 \\
\hline
Randomly coupling strength, $\zeta_{i,j} = \zeta_{j,i}$ & $\sim N(0, 1)$ \\
\hline
Effective Randomly coupling strength, $\zeta_{i,j} /N = \zeta_{j,i} / N$ & $\sim N(0, 1) / 10 $ \\
\hline
Additive Gaussian noise, $\xi_{i}$ & $\sim N(0, 0.25)$ \\
\hline
observable, $\vec f$ & $z = \exp(i\theta)$ \\ 
\hline
Hankel delay embedding & 300 \\
\hline
Natural frequencies, $\omega_i$ & $\sim U([0.25, 0.75])$ \\
\hline
Initial conditions, $\theta_i(0)$ & $\sim U([0, 2\pi])$ \\
\hline
\end{tabular}
\end{center}
\label{tab:kuramoto}
\end{table}%

\begin{figure}[htbp]
\begin{center}
\includegraphics[width=0.45\textwidth]{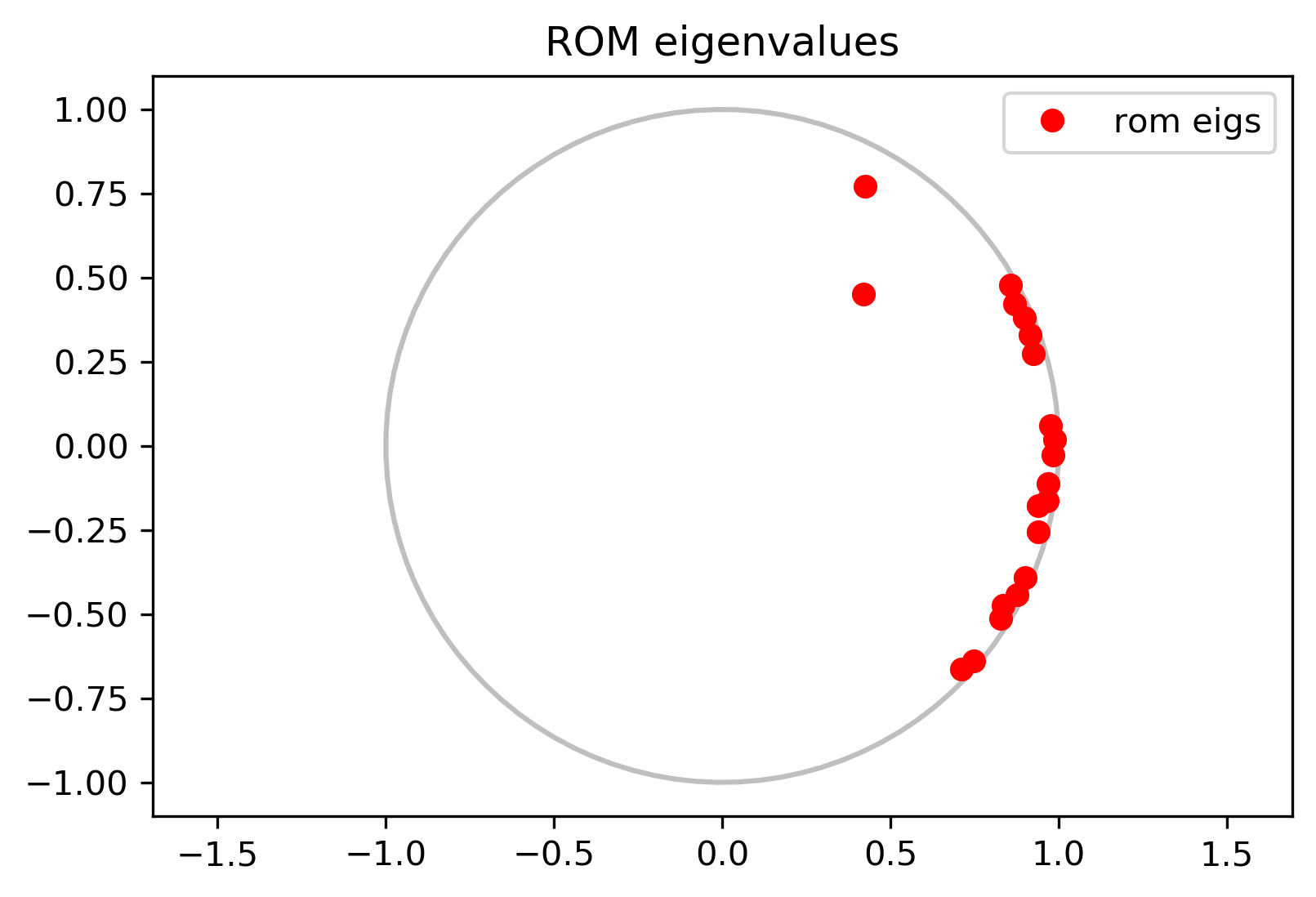}
\includegraphics[width=0.45\textwidth]{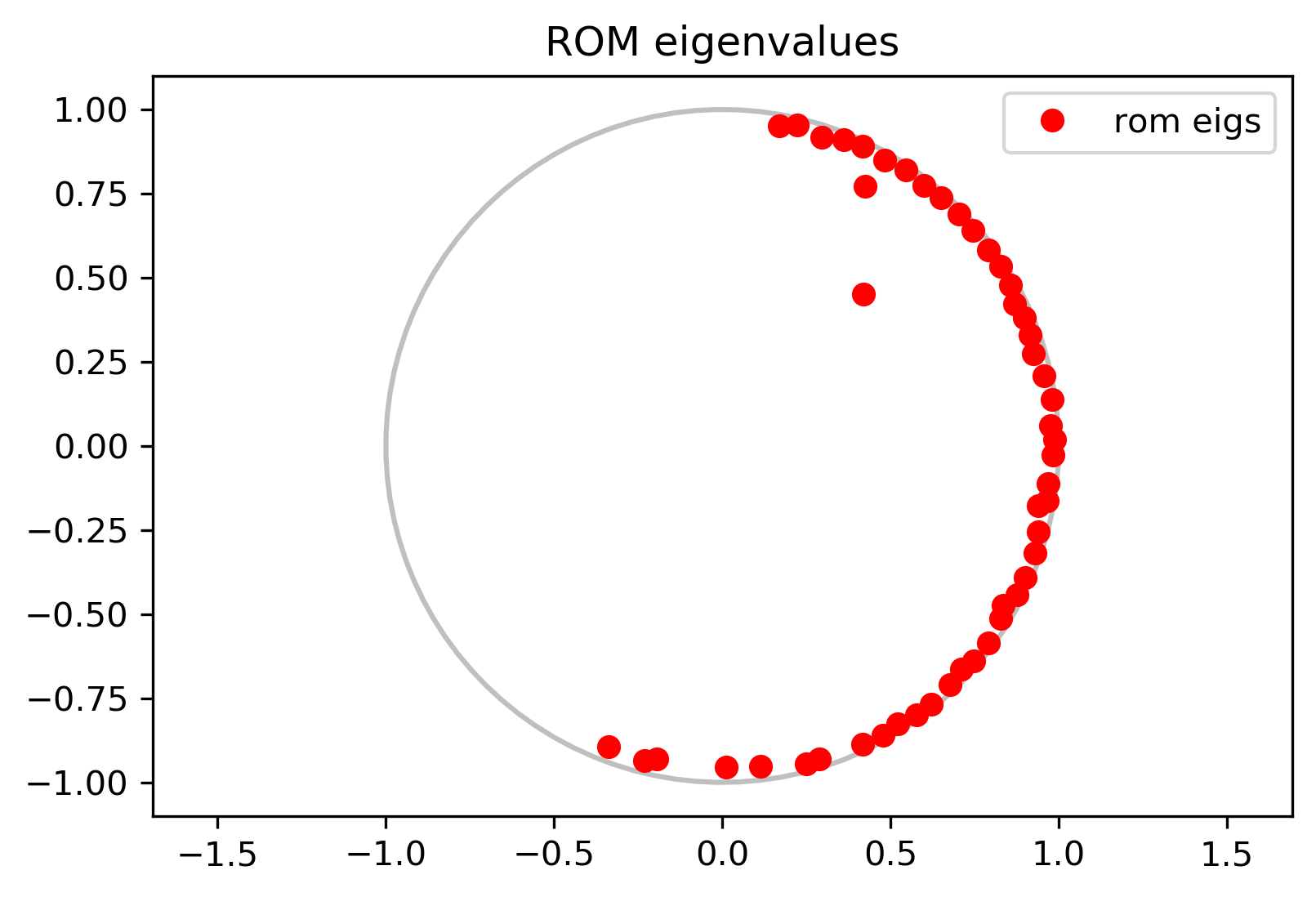}
\includegraphics[width=0.45\textwidth]{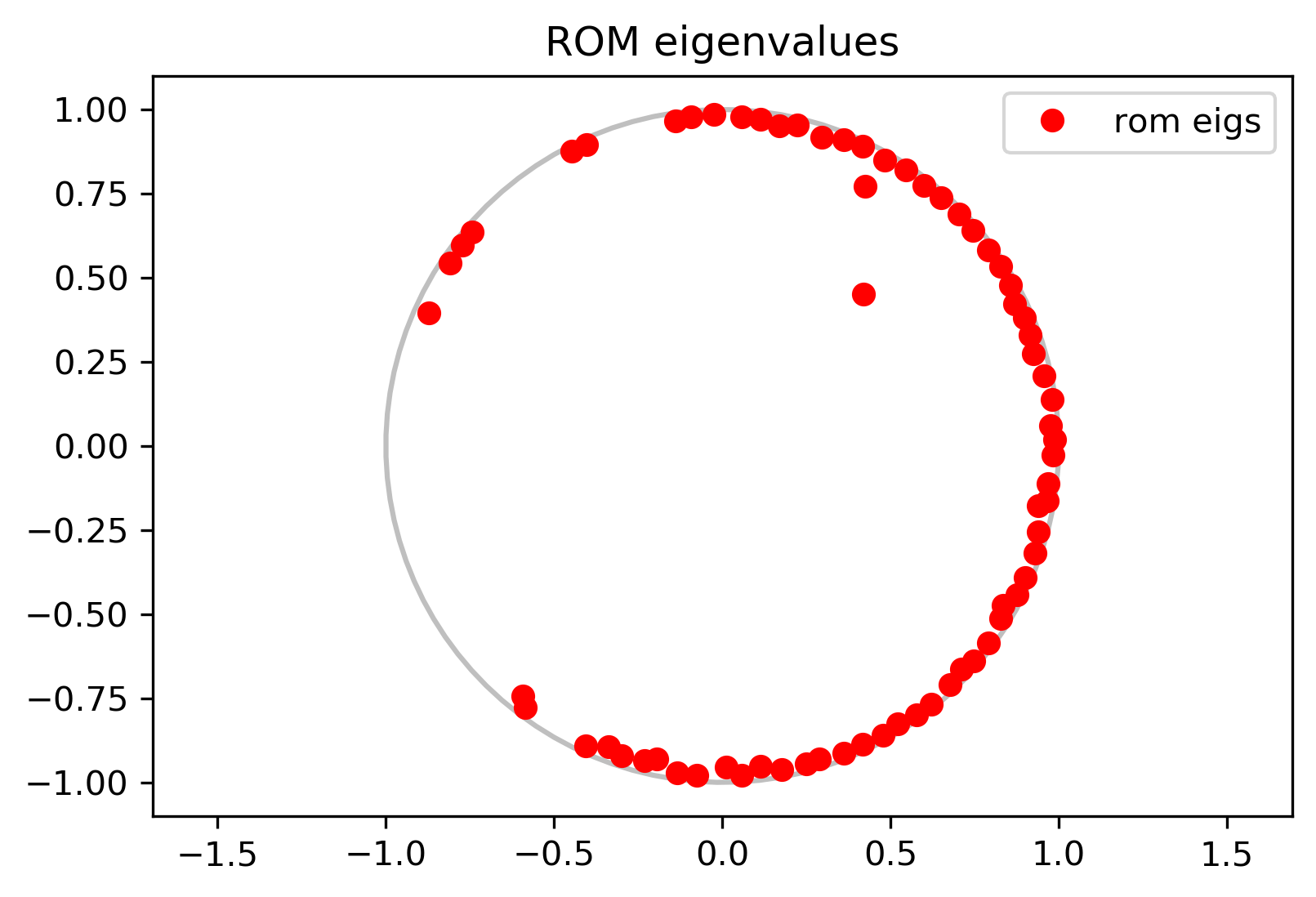}
\includegraphics[width=0.45\textwidth]{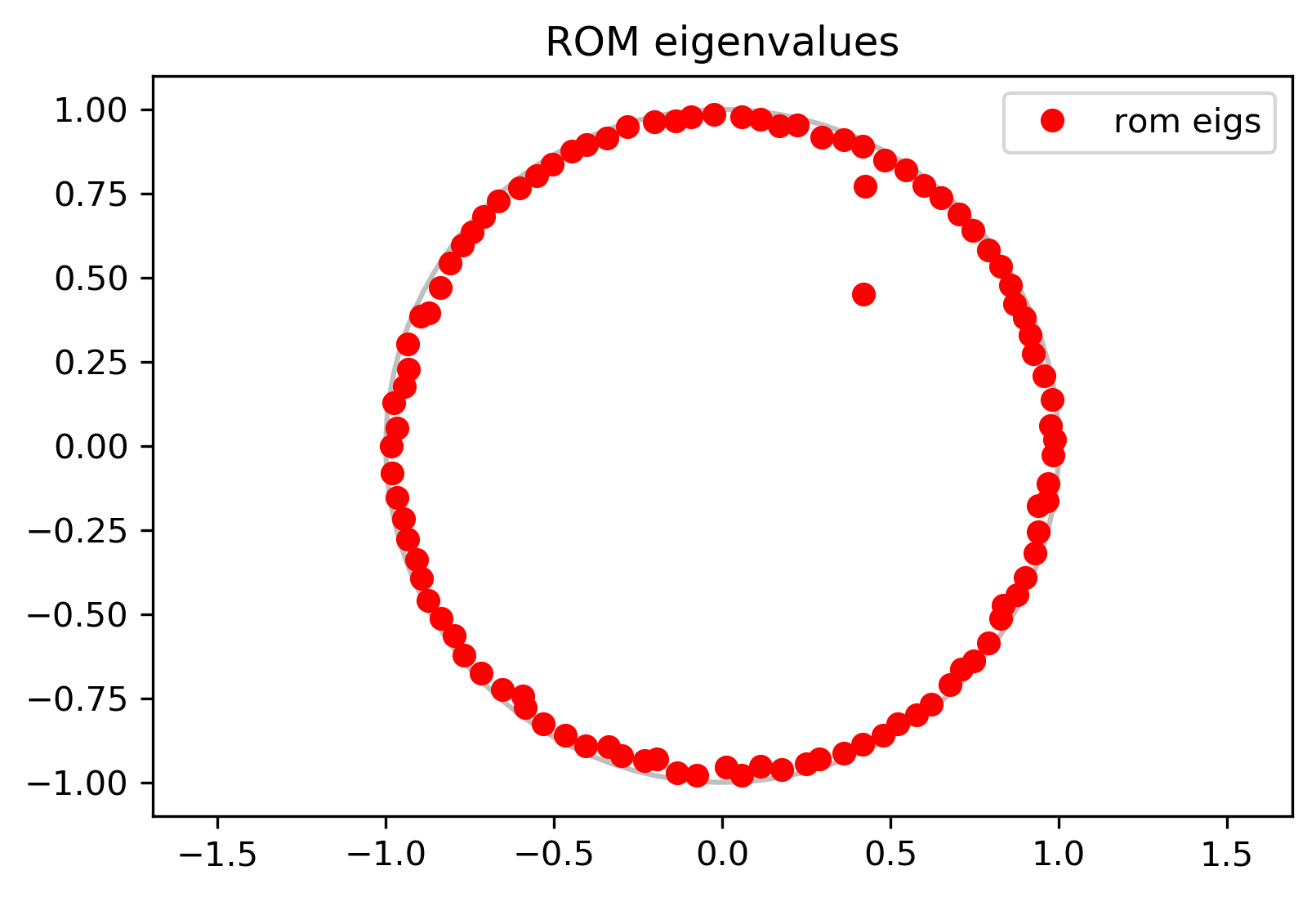}
\caption{\textbf{Kuramoto models, ROM eigenvalues}: Top left: ROM with 20 modes. Top right: ROM with 50 modes. Bottom left: ROM with 70 modes. Bottom right: ROM with 100 modes (full reconstruction). Eigenvalues of lower-order models are strict subsets of higher-order models. We note that we use the complexification of the real angles, $\theta \mapsto \exp(i2\pi\theta)$, in the construction of the data matrix. Since the data matrix is complex rather than real, complex-conjugate pairs of eigenvalues are not enforced, resulting in an asymmetric spectrum.}
\label{fig:kuramoto-eigenvalues-complex-angles}
\end{center}
\end{figure}

\begin{figure}[ht]
\centering
    \begin{subfigure}{0.45\textwidth}
        \centering
        \includegraphics[width = \textwidth, height=0.15\textheight]{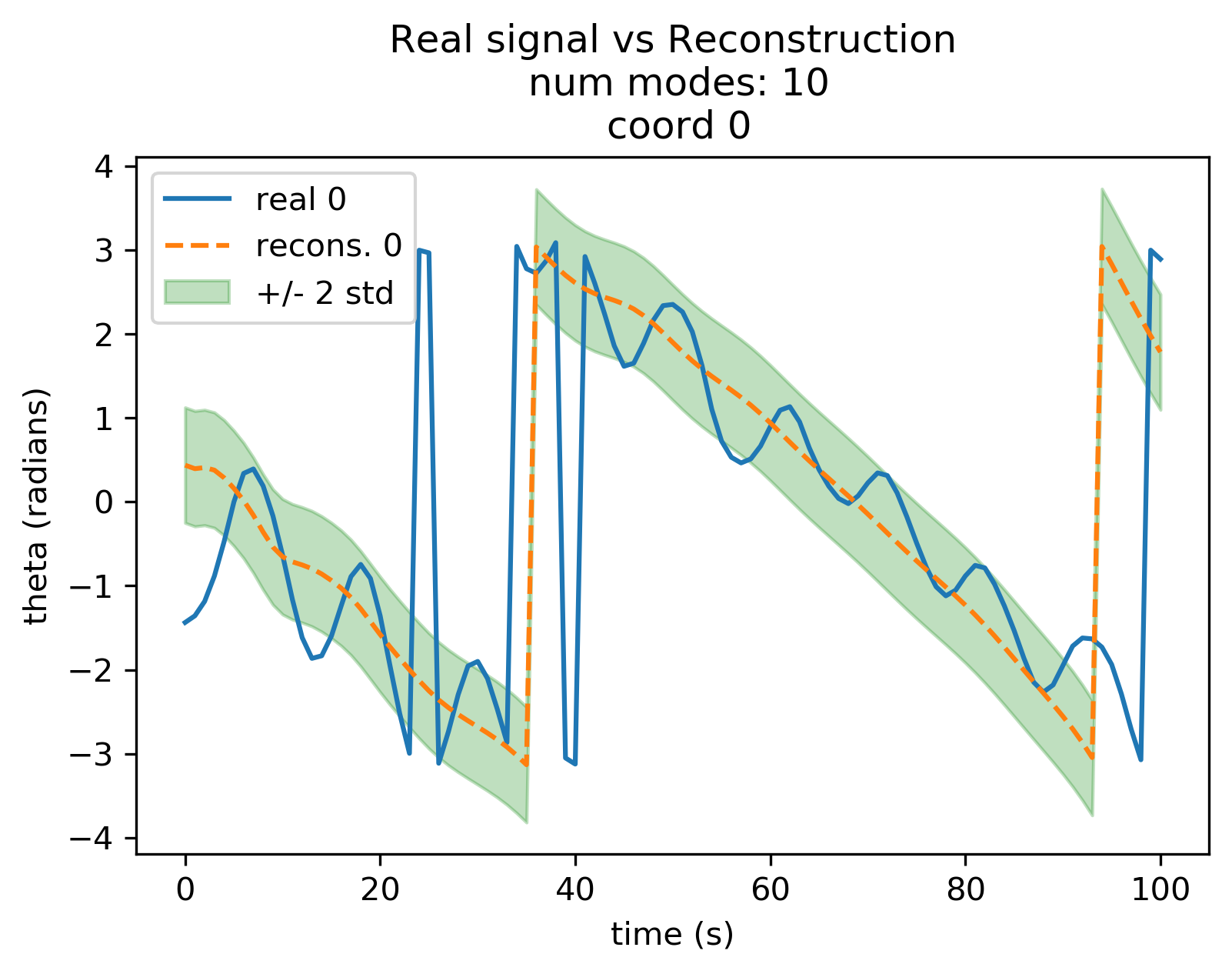}
        \caption{10 modes}
    \end{subfigure}
    \hfill
    \begin{subfigure}{0.45\textwidth}
        \centering
        \includegraphics[width = \textwidth, height=0.15\textheight]{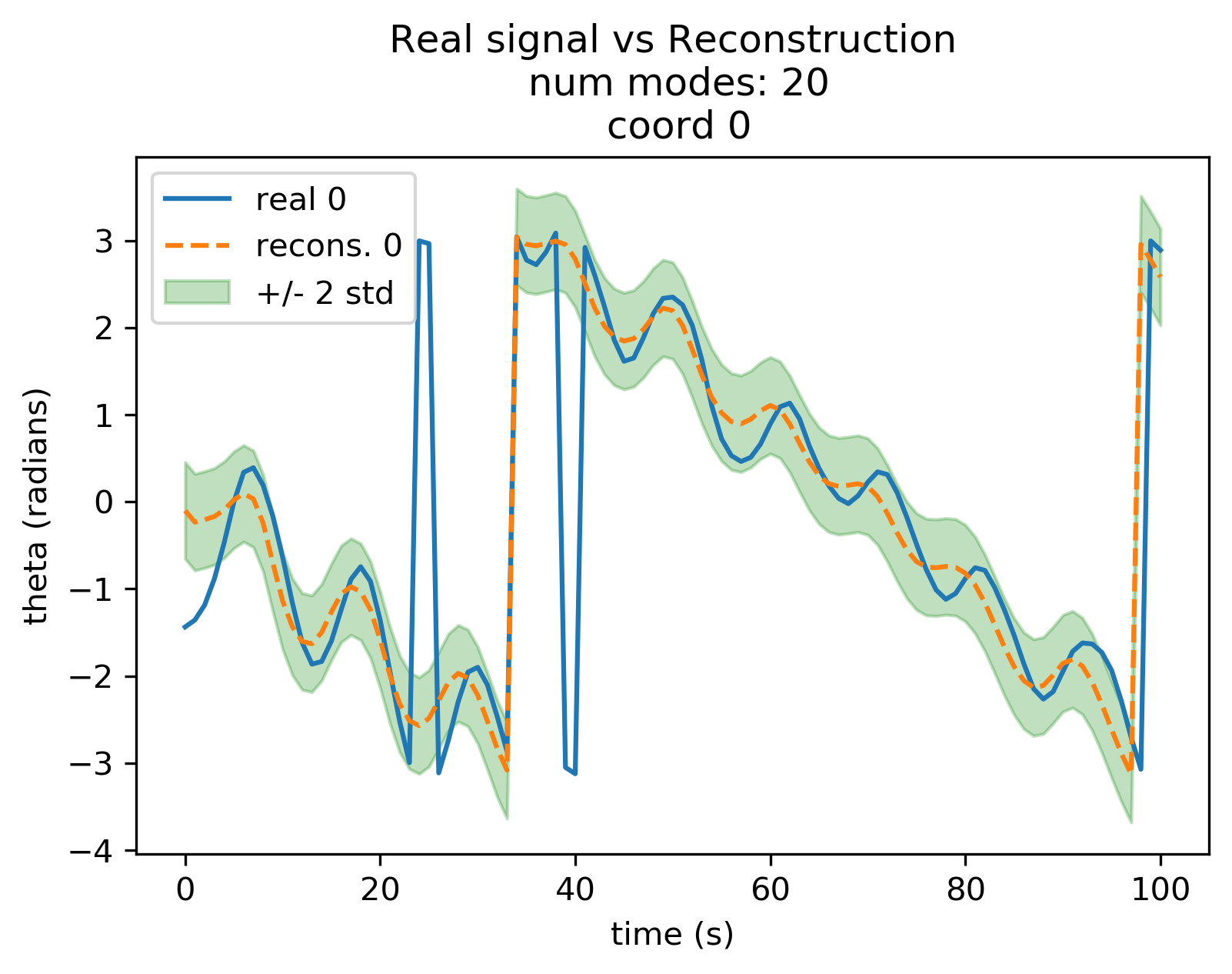}
        \caption{20 modes}
    \end{subfigure}
    \hfill
    \begin{subfigure}{0.45\textwidth}
        \centering
        \includegraphics[width = \textwidth, height=0.15\textheight]{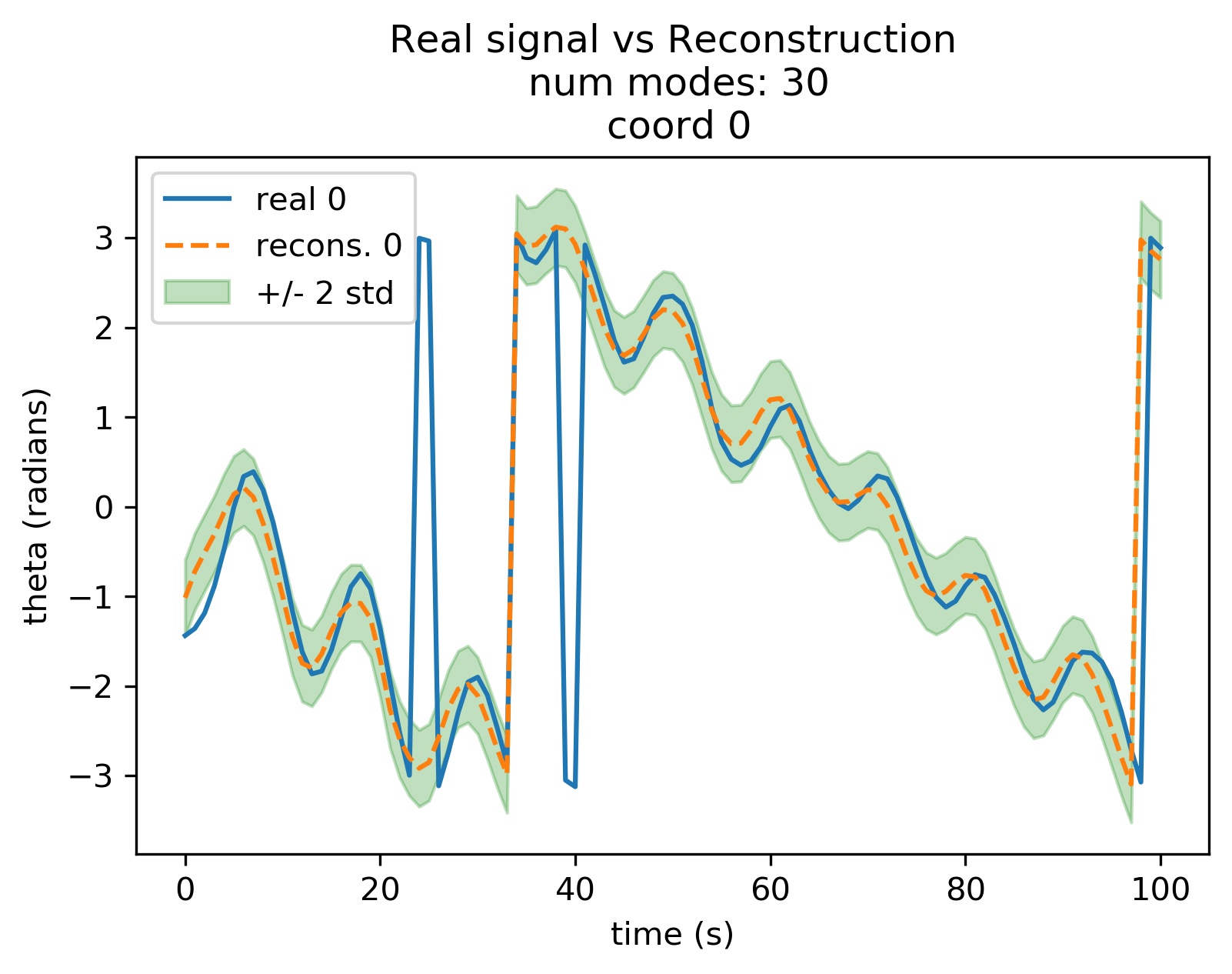}
        \caption{30 modes}
    \end{subfigure}
    \hfill
    \begin{subfigure}{0.45\textwidth}
        \centering
        \includegraphics[width = \textwidth, height=0.15\textheight]{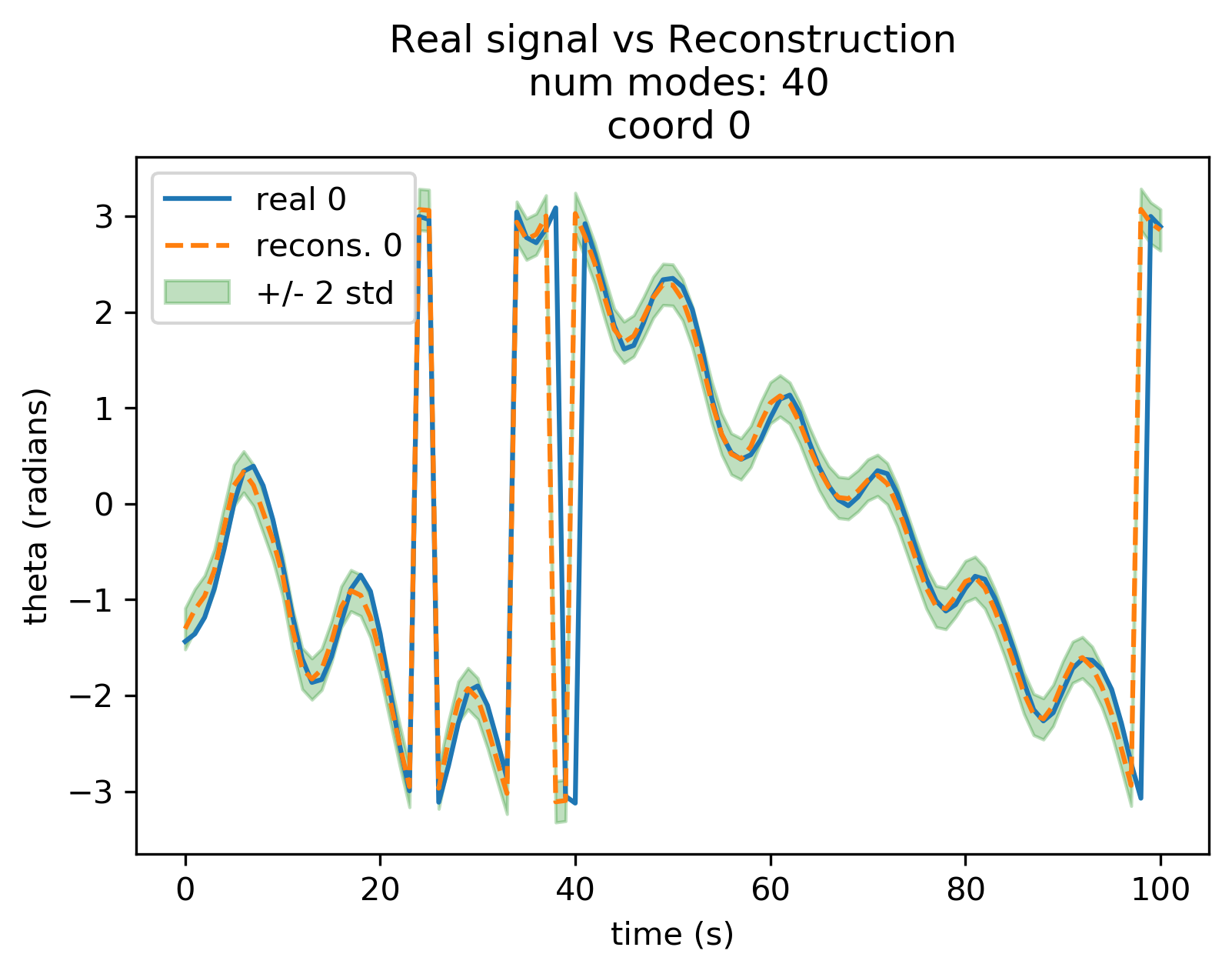}
        \caption{40 modes}
    \end{subfigure}
    \hfill
    \begin{subfigure}{0.45\textwidth}
        \centering
        \includegraphics[width = \textwidth, height=0.15\textheight]{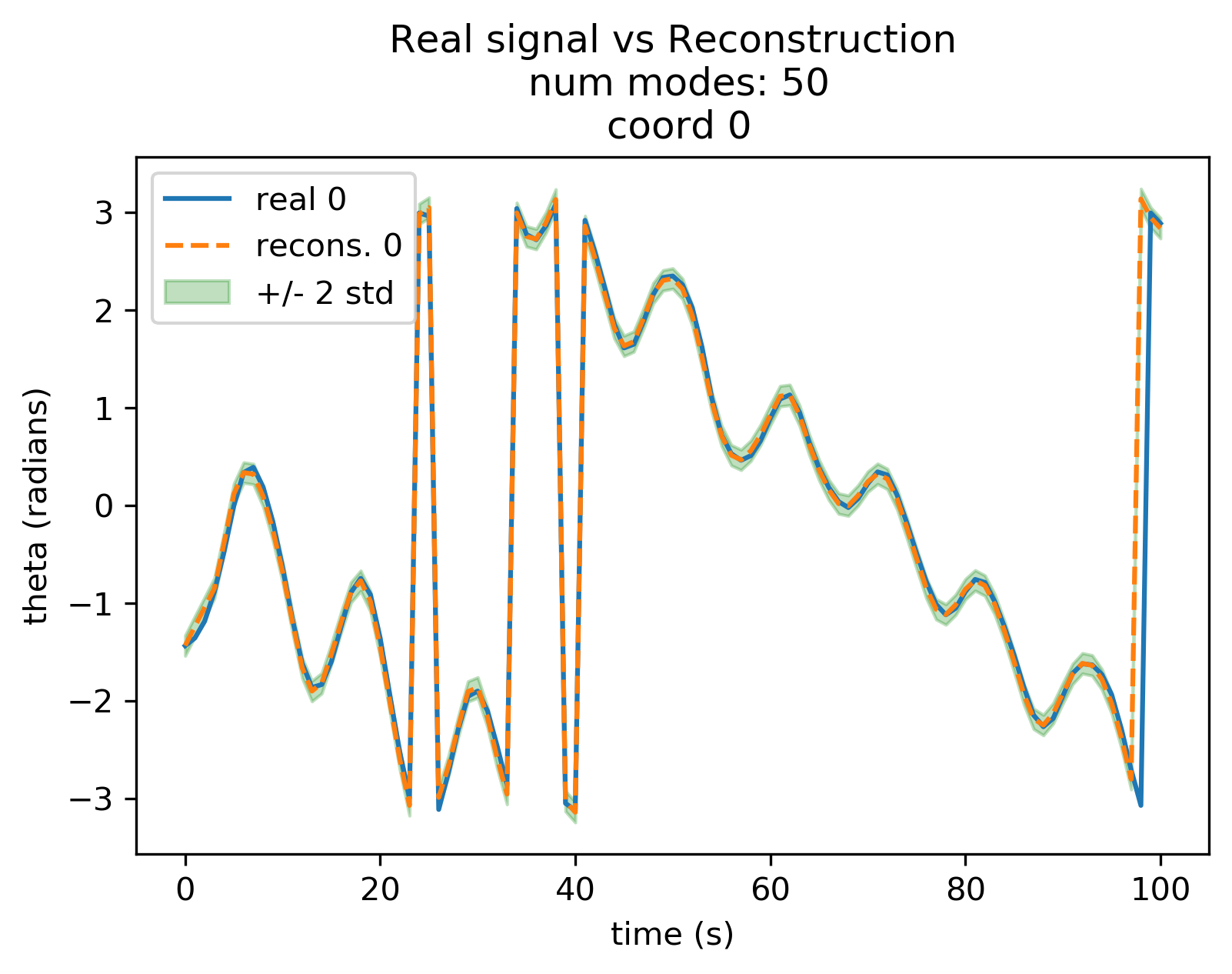}
        \caption{50 modes}
    \end{subfigure}
    \hfill
    \begin{subfigure}{0.45\textwidth}
        \centering
        \includegraphics[width=\textwidth, height=0.15\textheight]{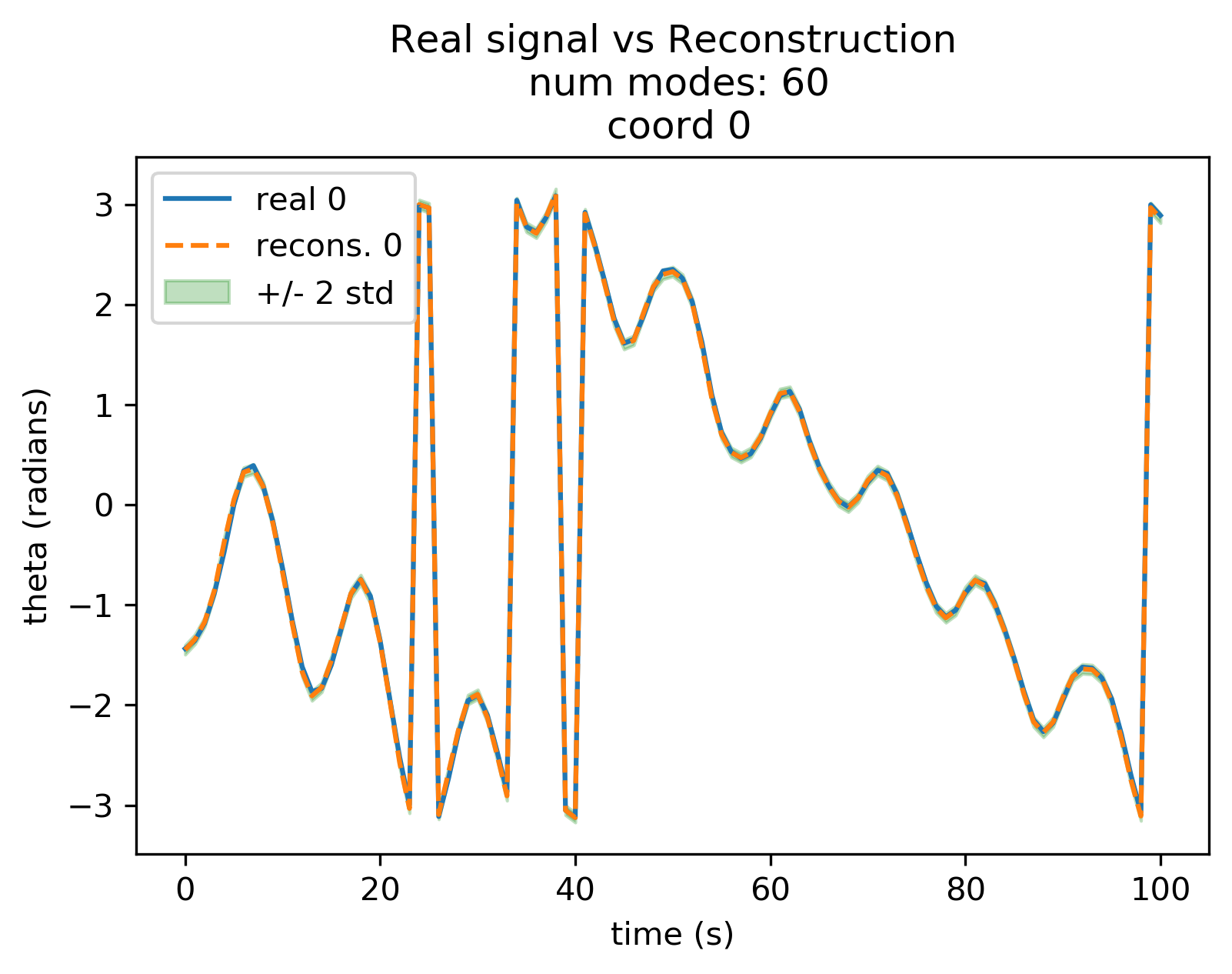}
        \caption{60 modes}
    \end{subfigure}
    \hfill
    \begin{subfigure}{0.45\textwidth}
        \centering
        \includegraphics[width=\textwidth, height=0.15\textheight]{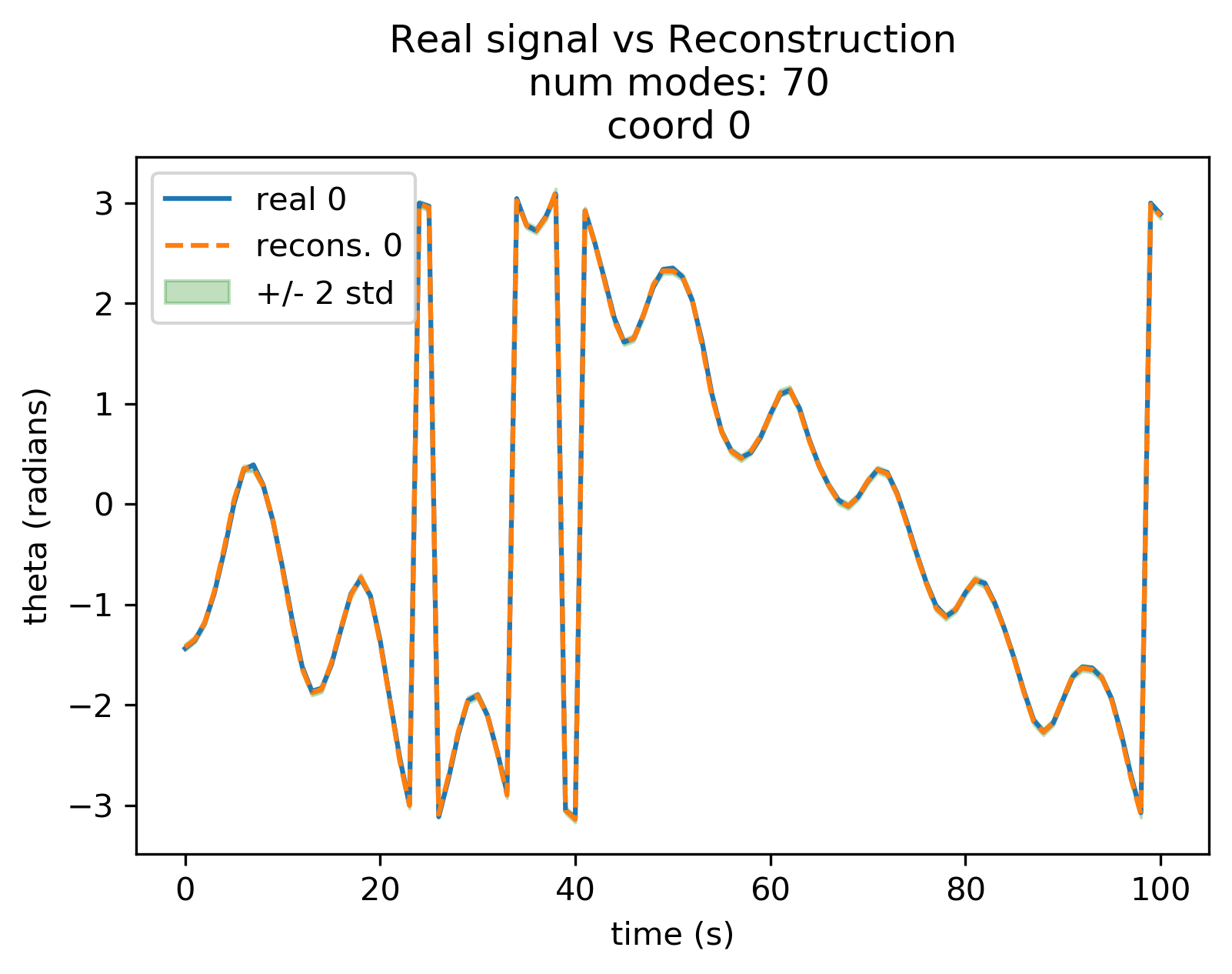}
        \caption{70 modes}
    \end{subfigure}
    \hfill
    \begin{subfigure}{0.45\textwidth}
        \centering
        \includegraphics[width=\textwidth, height=0.15\textheight]{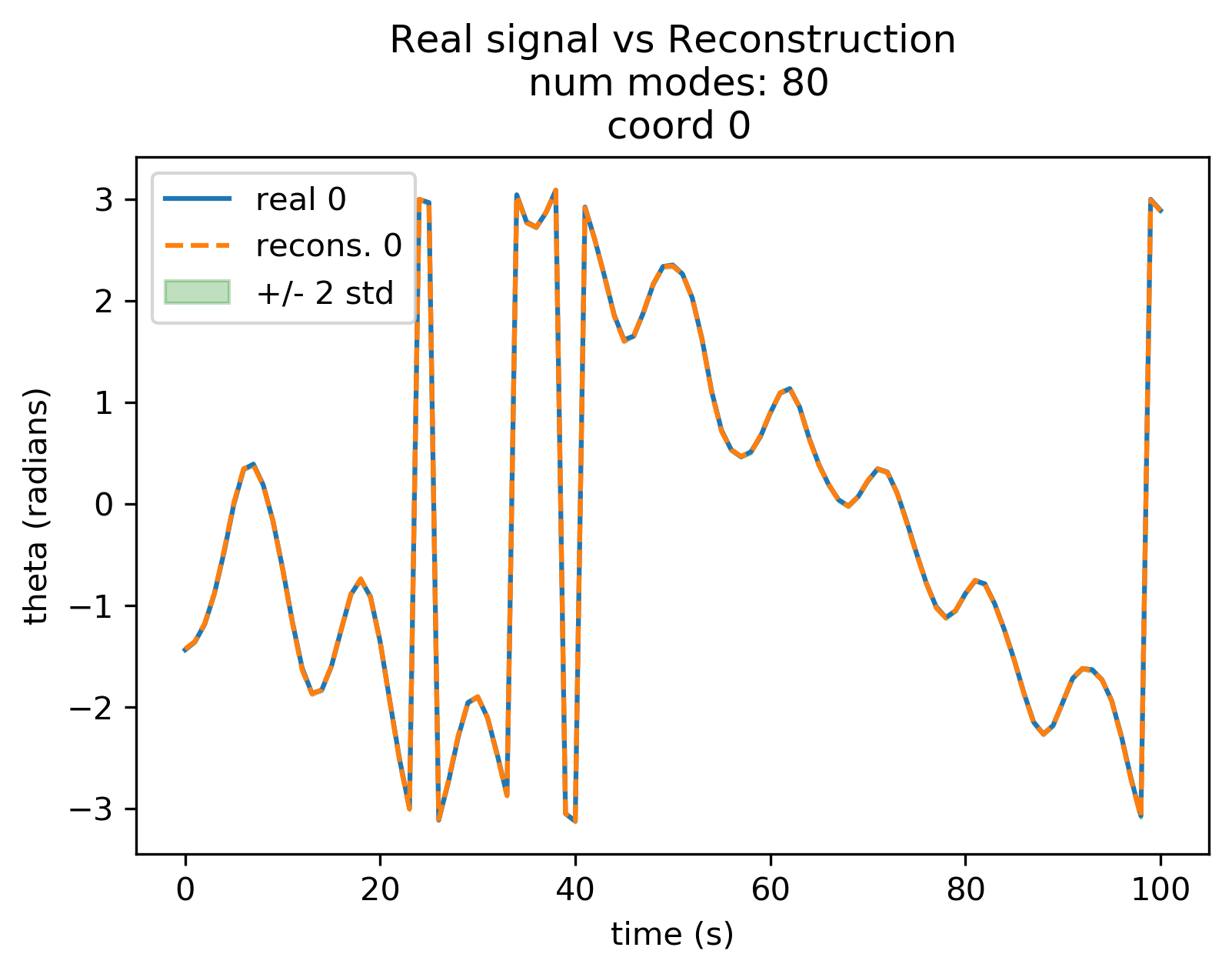}
        \caption{80 modes}
    \end{subfigure}
    \hfill
    \begin{subfigure}{0.45\textwidth}
        \centering
        \includegraphics[width=\textwidth, height=0.15\textheight]{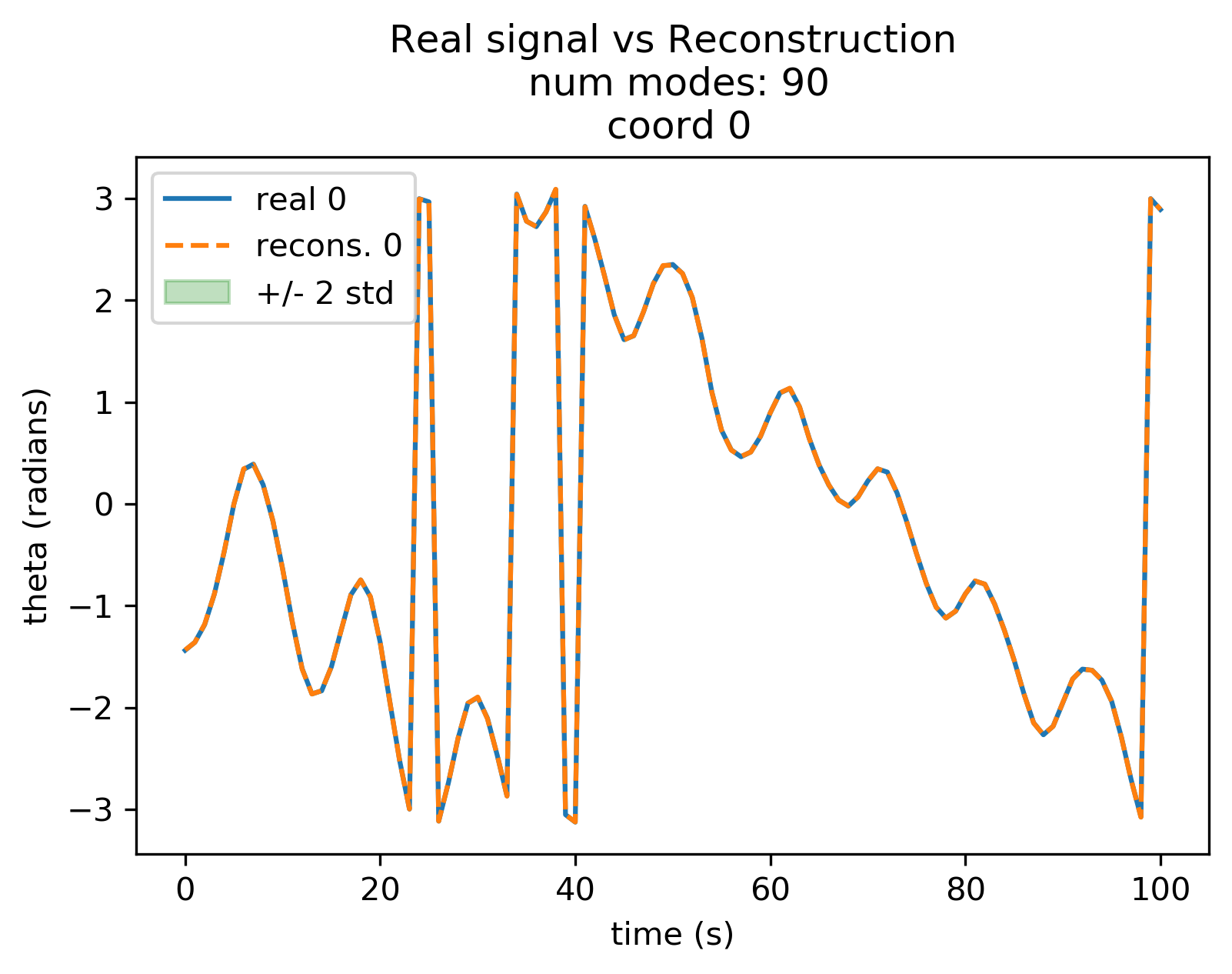}
        \caption{90 modes}
    \end{subfigure}
    \hfill
    \begin{subfigure}{0.45\textwidth}
        \centering
        \includegraphics[width=\textwidth, height=0.15\textheight]{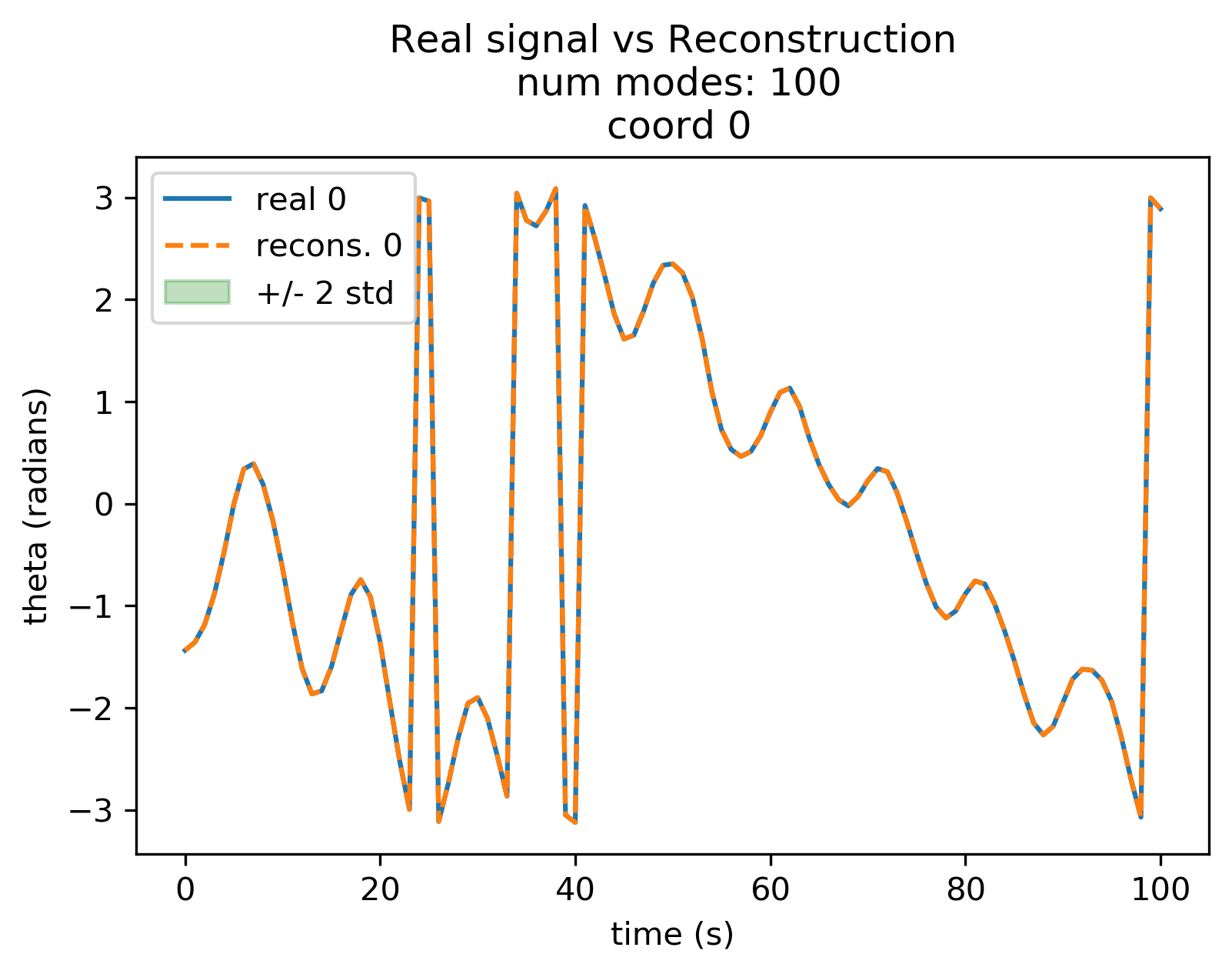}
        \caption{100 modes}
    \end{subfigure}
\caption{\textbf{Kuramoto models, signal vs reconstruction}: Comparison of the true signal with reconstruction signals using $10, 20, 30, \dots, 100$ modes for reconstruction for oscillator 0. The true signal is in blue. The reconstruction signal is orange and the green band is $\pm 2$ standard deviations of the modal noise.}
\label{fig:kuramoto-reconstruction-complex-angles}
\end{figure}

\begin{figure}[htbp]
\centering
\begin{subfigure}[b]{0.45\textwidth}
    \centering
    \includegraphics[width = \textwidth, height=0.15\textheight]{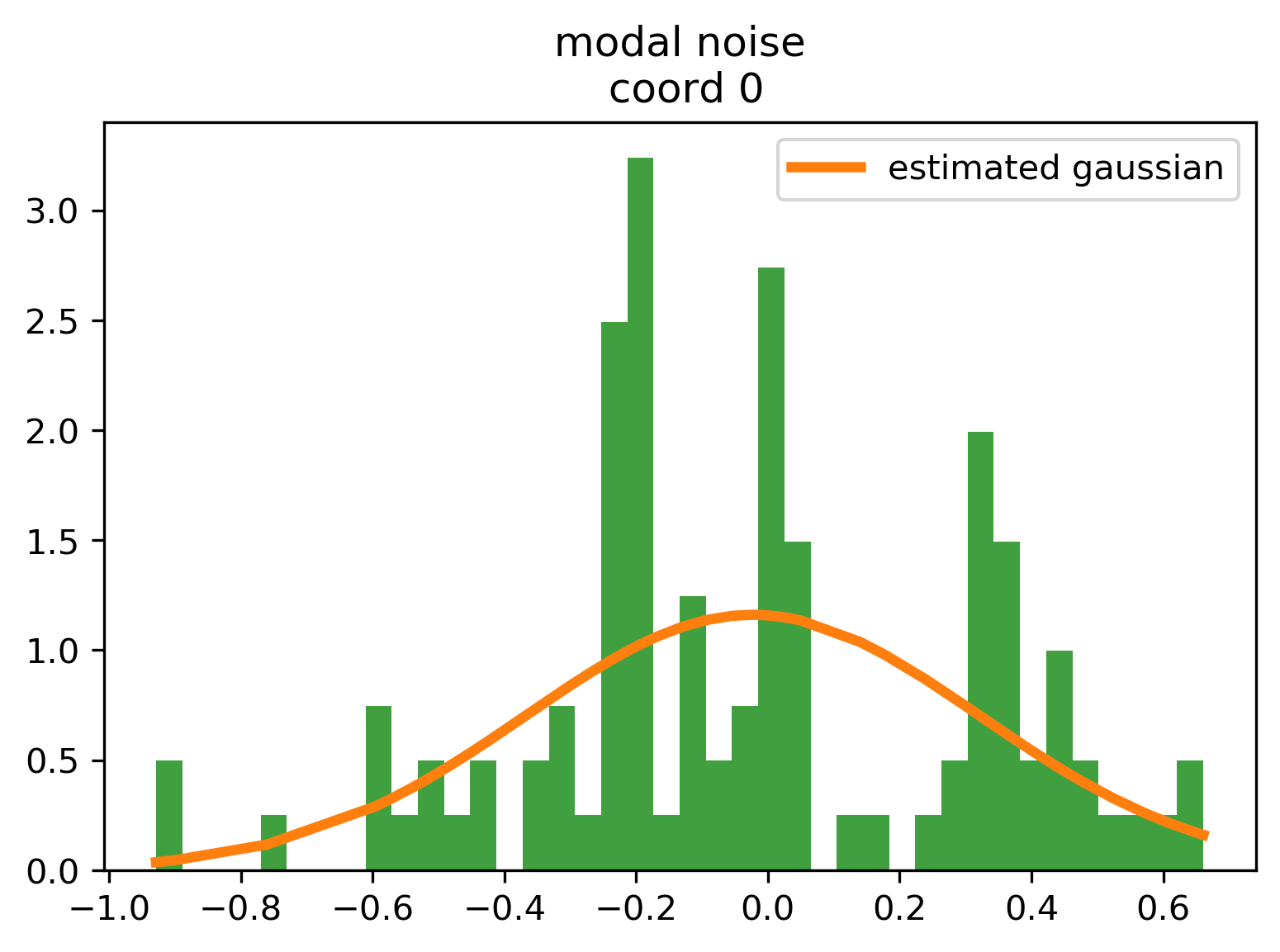}
    \caption{10 modes}
\end{subfigure}
\begin{subfigure}[b]{0.45\textwidth}
    \centering
    \includegraphics[width = \textwidth, height=0.15\textheight]{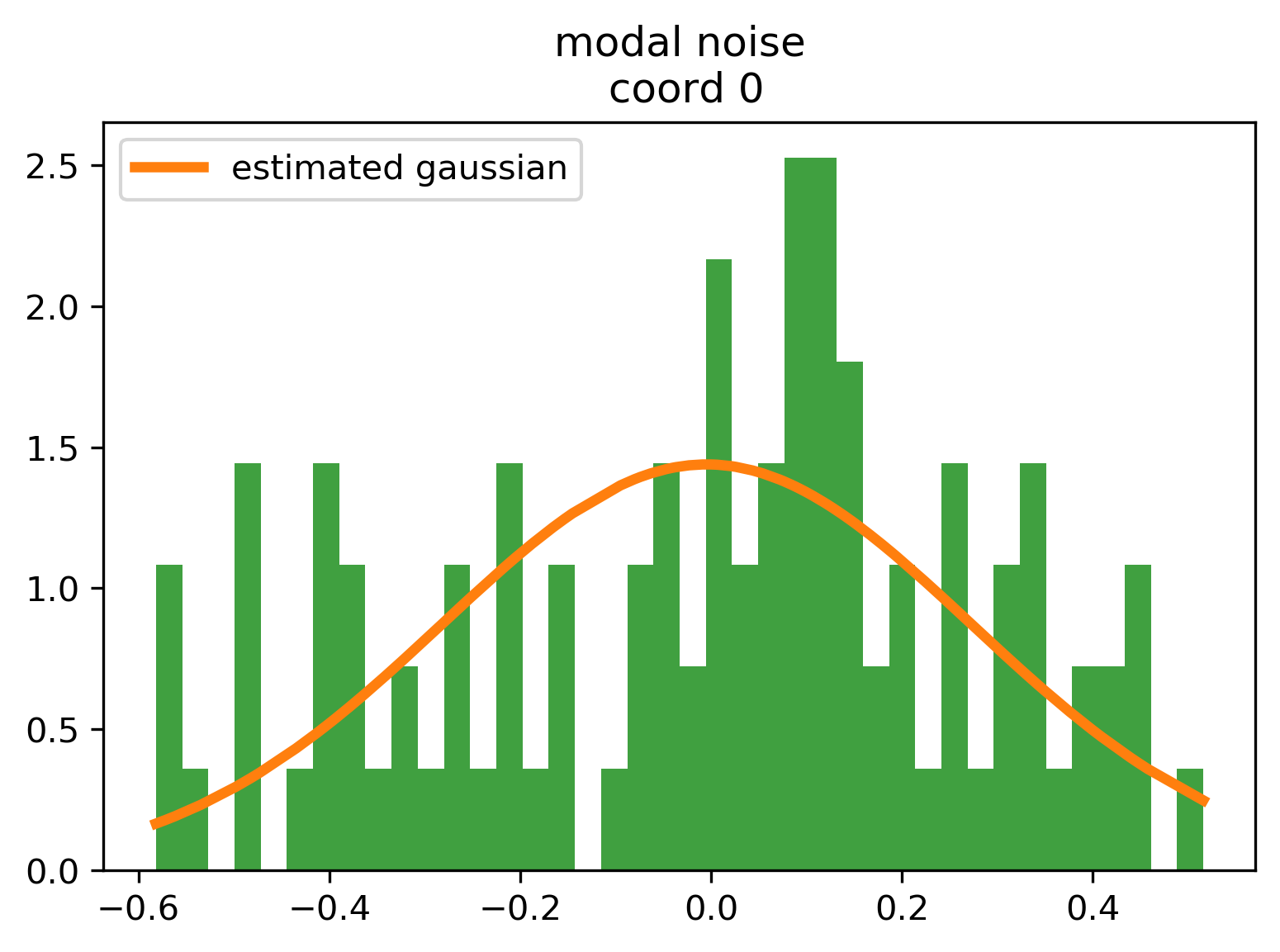}
    \caption{20 modes}
\end{subfigure}
\begin{subfigure}[b]{0.45\textwidth}
    \centering
    \includegraphics[width = \textwidth, height=0.15\textheight]{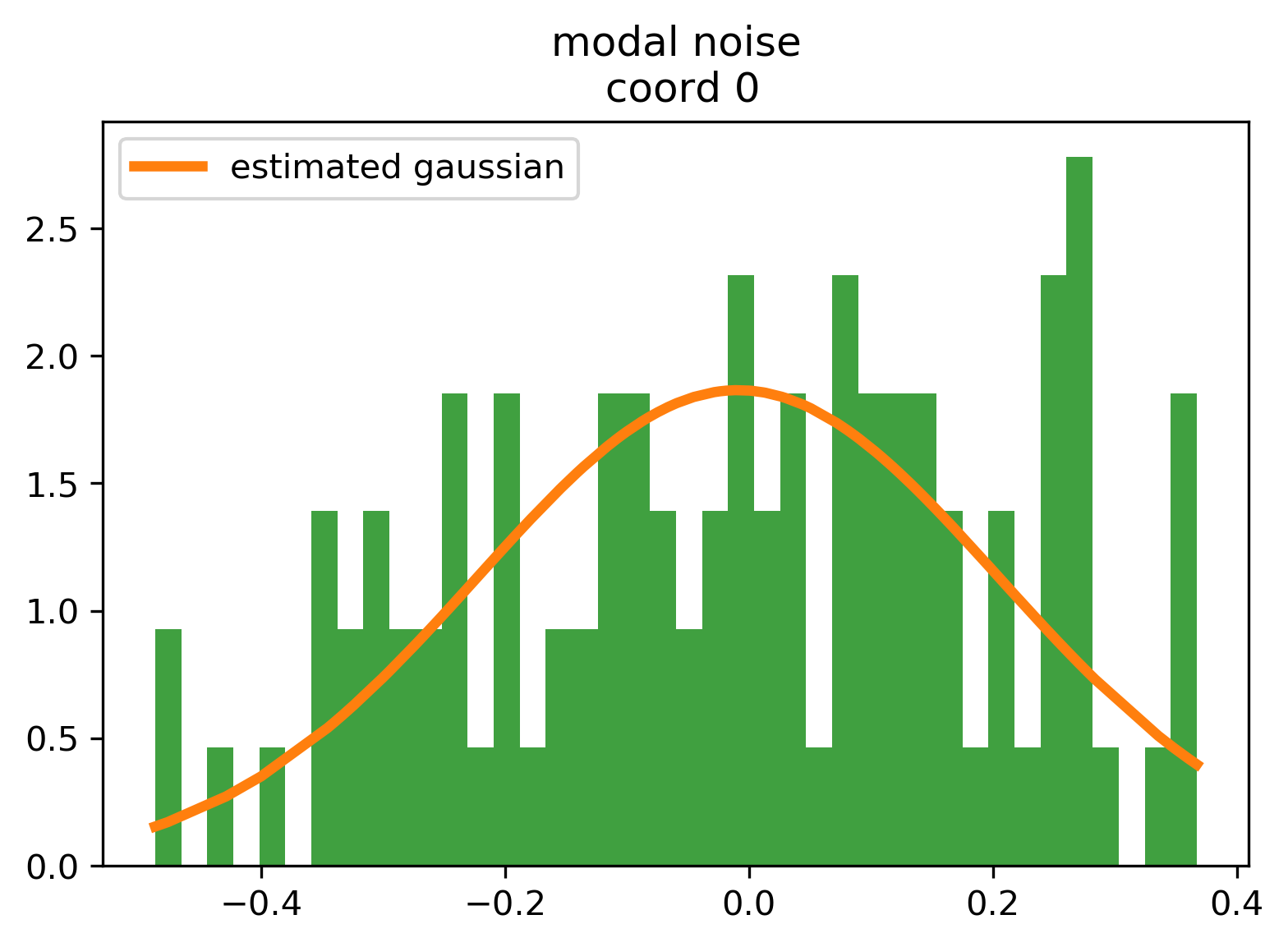}
    \caption{30 modes}
\end{subfigure}
\begin{subfigure}[b]{0.45\textwidth}
    \centering
    \includegraphics[width = \textwidth, height=0.15\textheight]{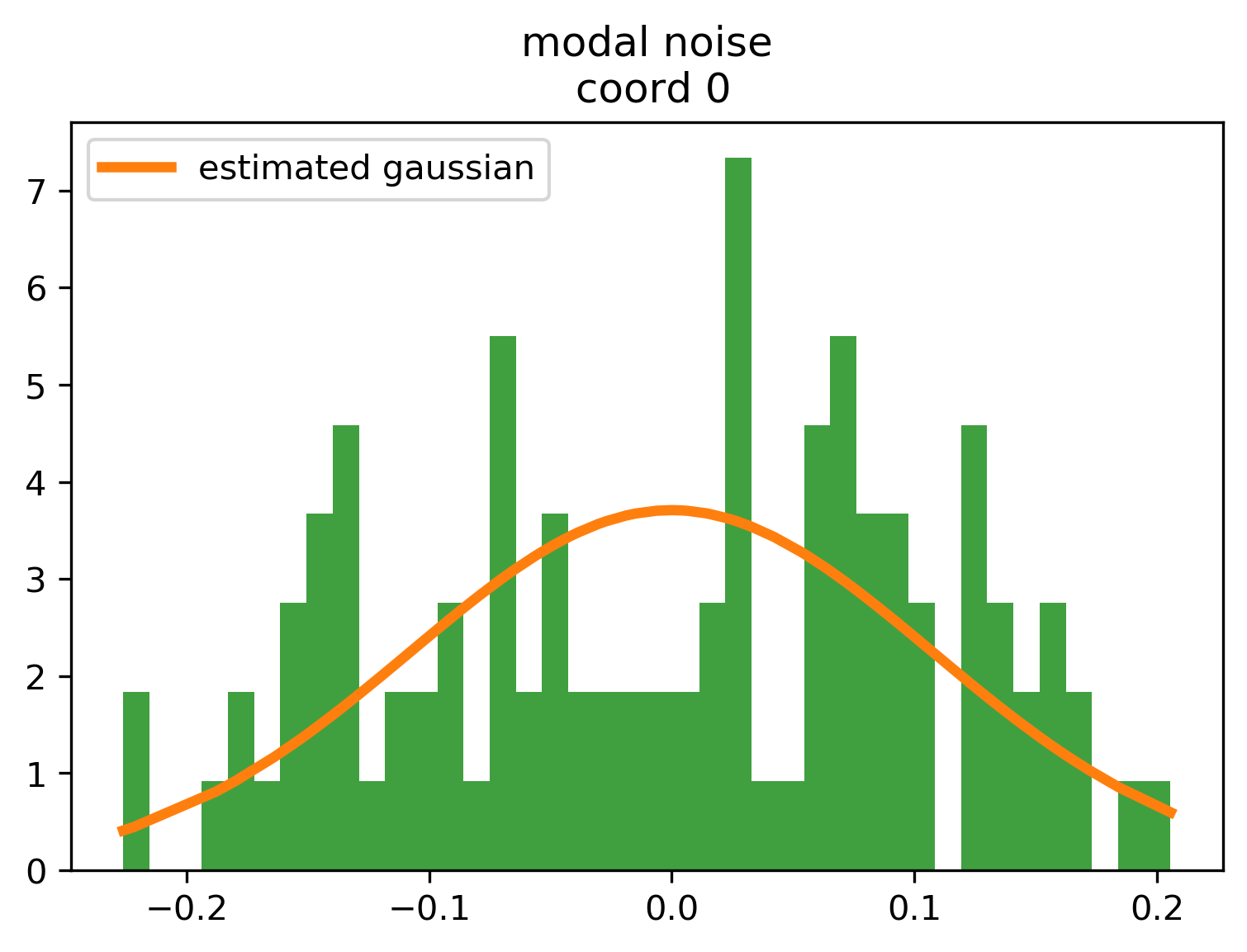}
    \caption{40 modes}
\end{subfigure}
\begin{subfigure}[b]{0.45\textwidth}
    \centering
    \includegraphics[width = \textwidth, height=0.15\textheight]{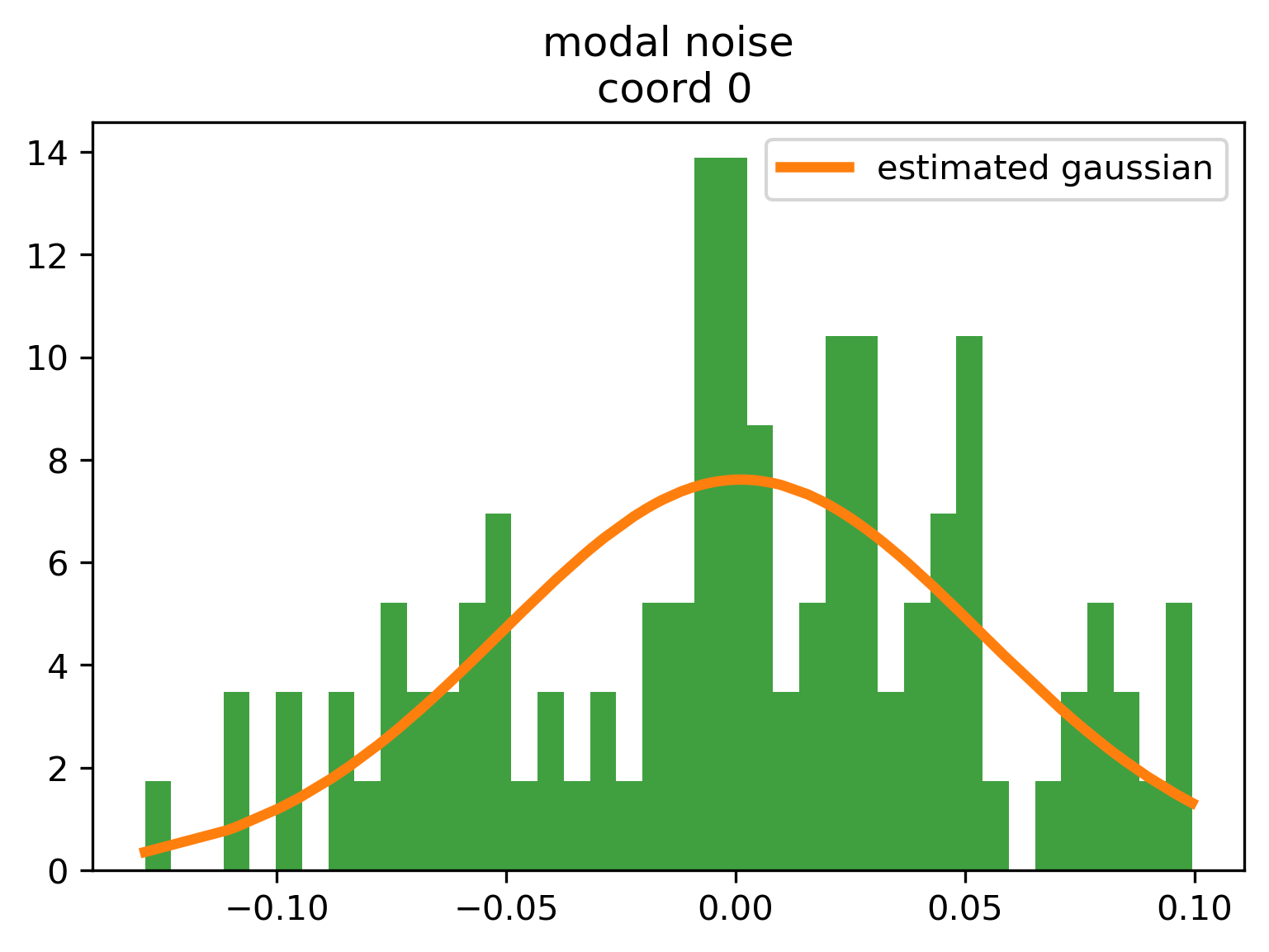}
    \caption{50 modes}
\end{subfigure}
\begin{subfigure}[b]{0.45\textwidth}
    \centering
    \includegraphics[width = \textwidth, height=0.15\textheight]{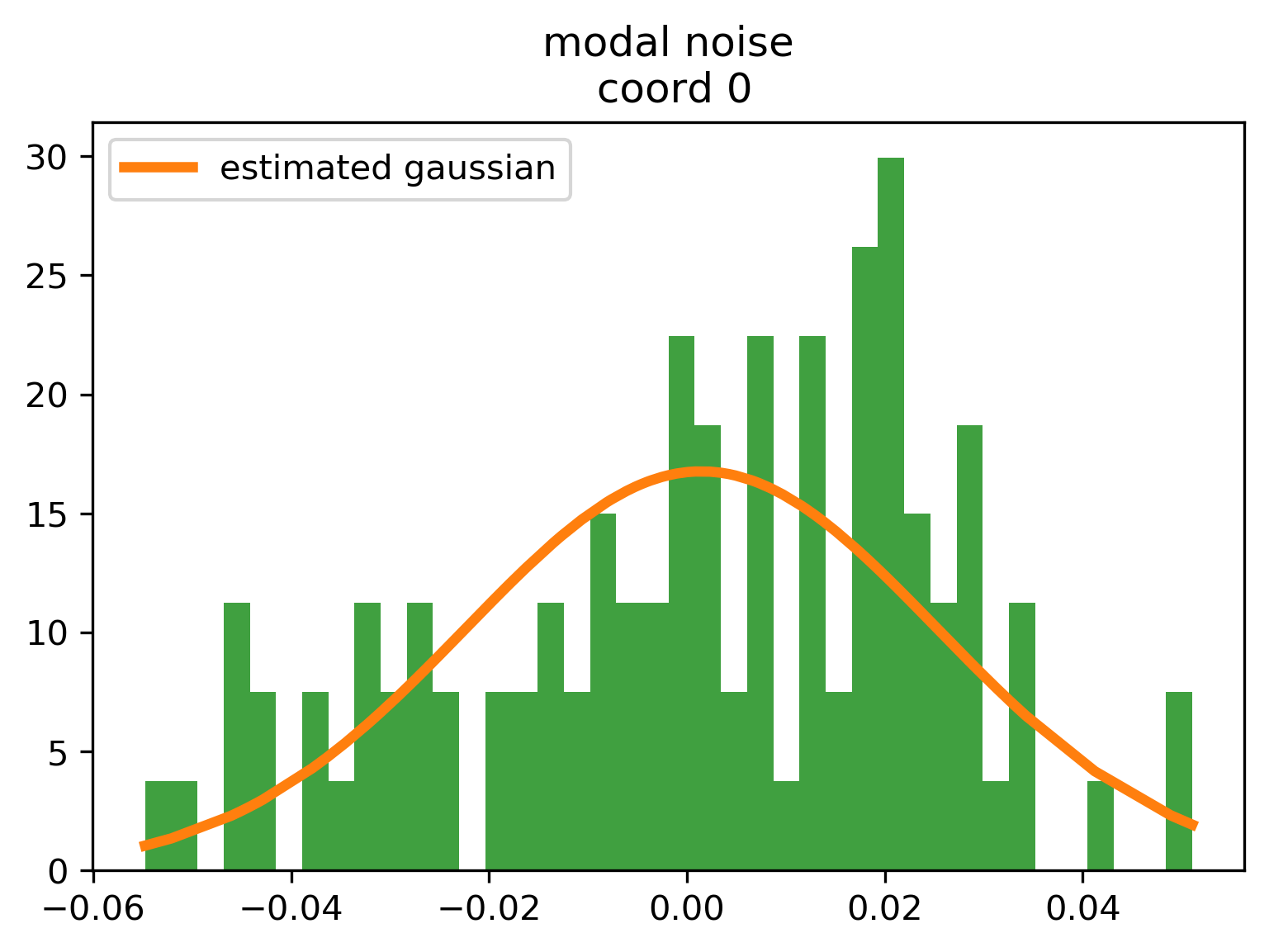}
    \caption{60 modes}
\end{subfigure}
\begin{subfigure}[b]{0.45\textwidth}
    \centering
    \includegraphics[width = \textwidth, height=0.15\textheight]{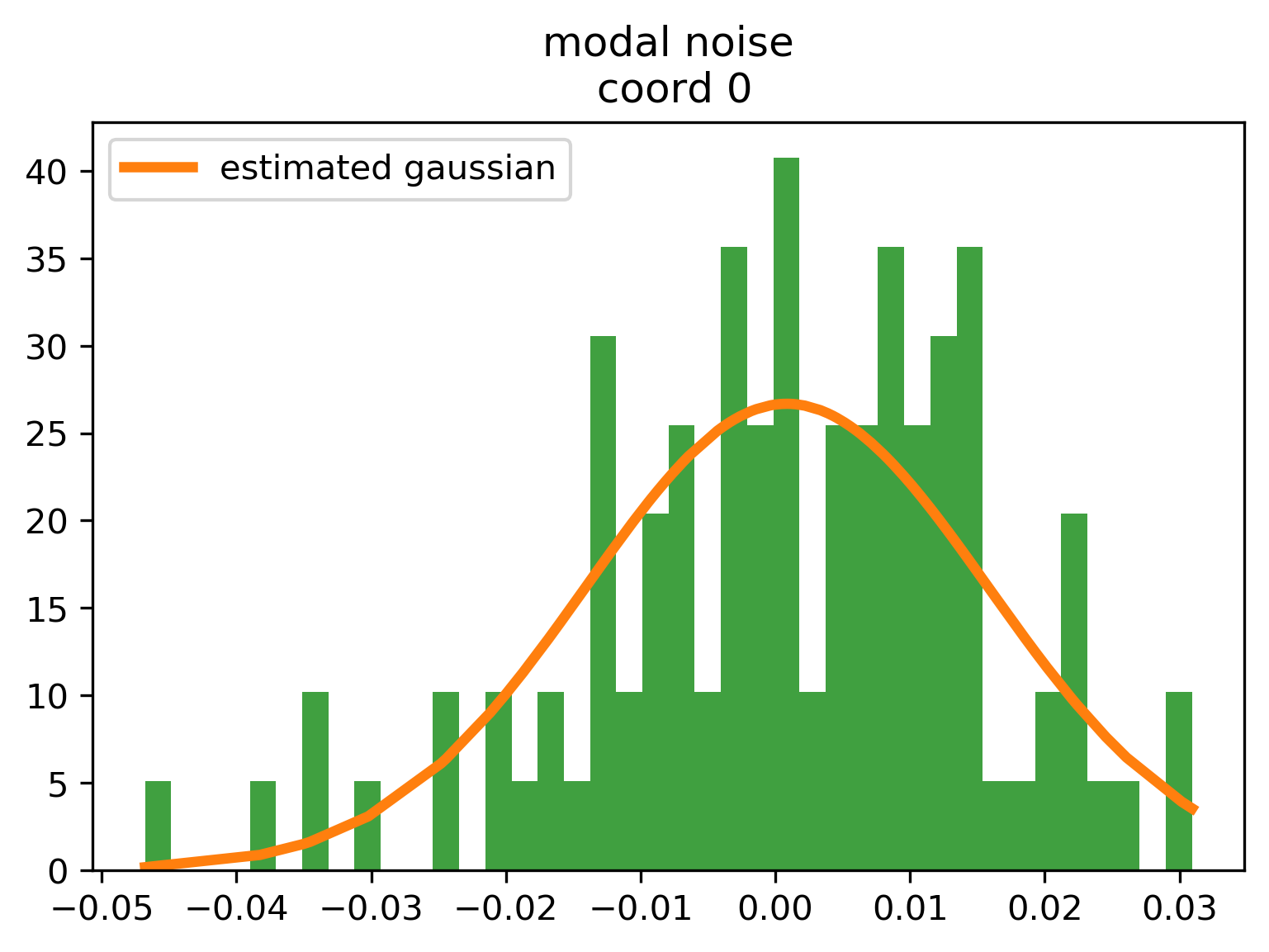}
    \caption{70 modes}
\end{subfigure}
\begin{subfigure}[b]{0.45\textwidth}
    \centering
    \includegraphics[width = \textwidth, height=0.15\textheight]{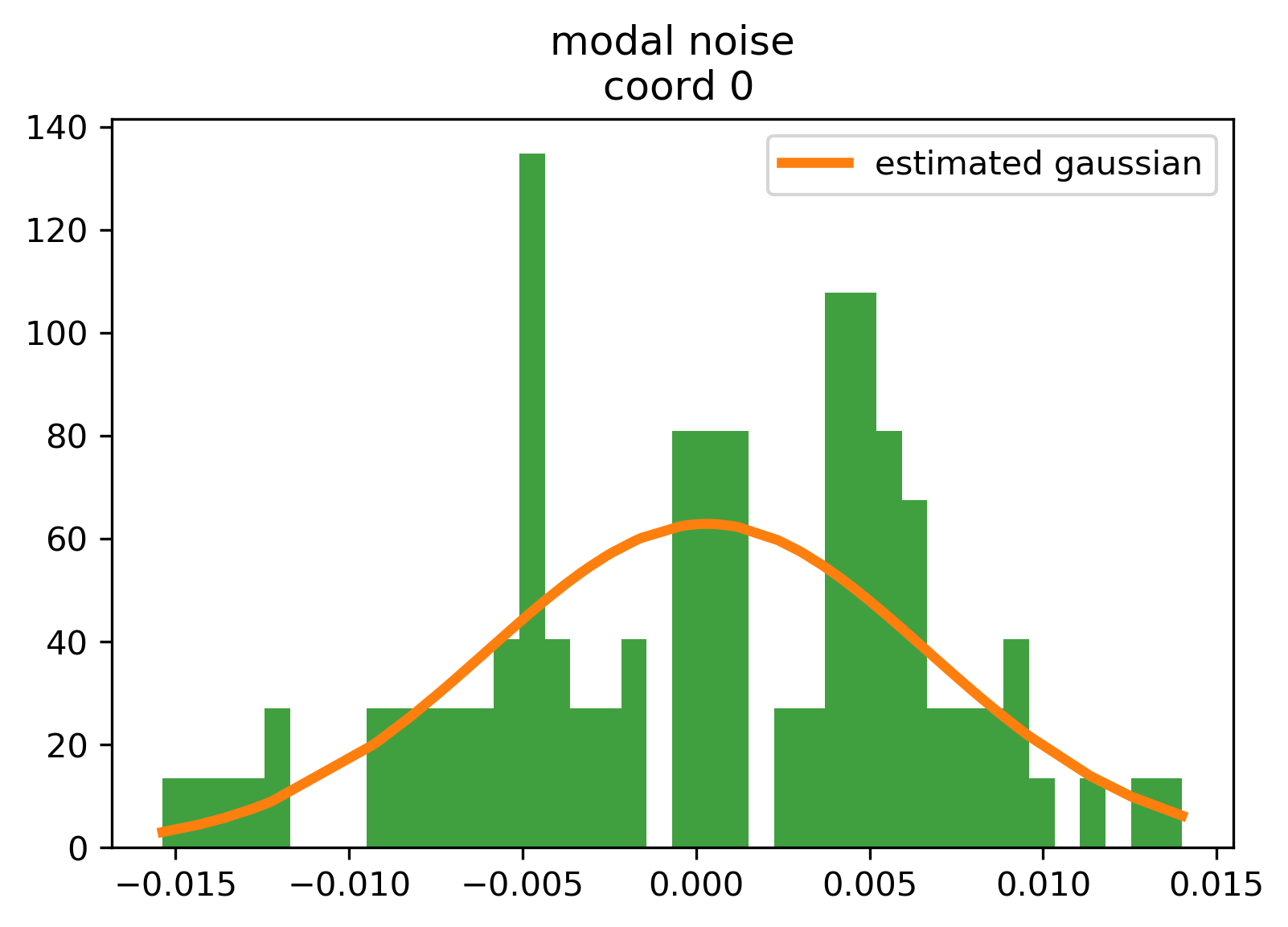}
    \caption{80 modes}
\end{subfigure}
\begin{subfigure}[b]{0.45\textwidth}
    \centering
    \includegraphics[width = \textwidth, height=0.15\textheight]{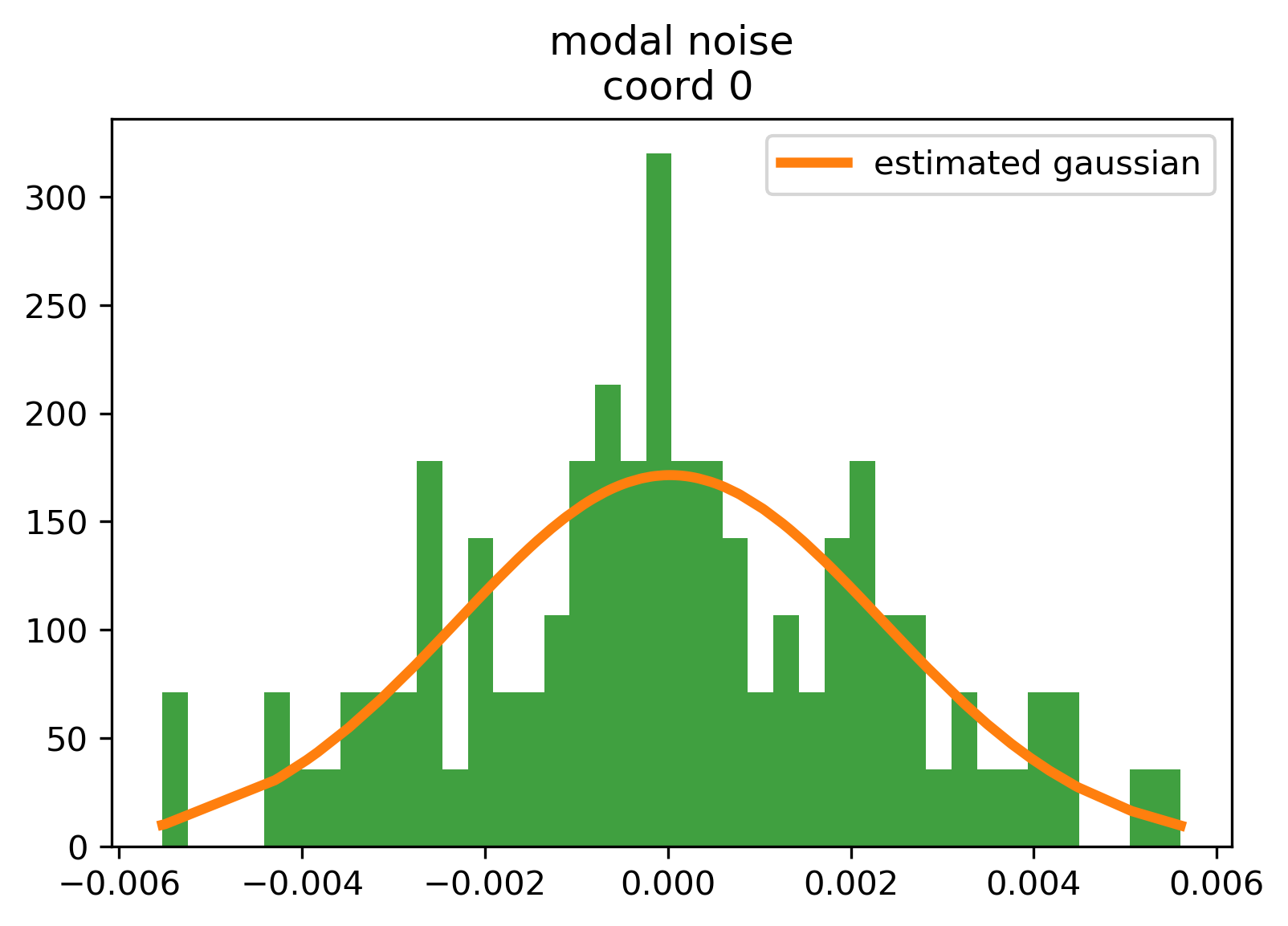}
    \caption{90 modes}
\end{subfigure}
\begin{subfigure}[b]{0.45\textwidth}
    \centering
    \includegraphics[width = \textwidth, height=0.15\textheight]{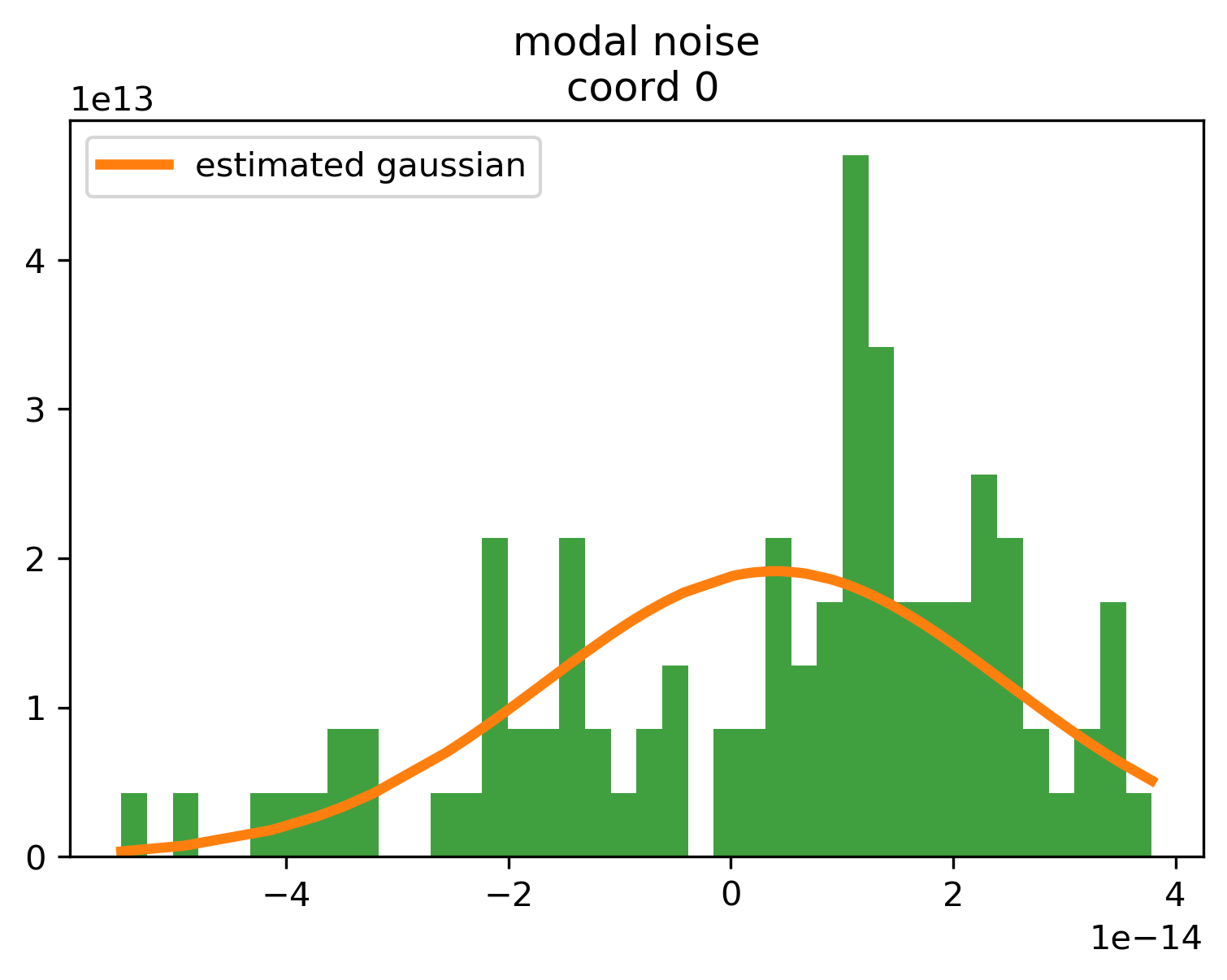}
    \caption{100 modes}
\end{subfigure}
\caption{\textbf{Kuramoto models}: Modal noise distributions vs number of modes used for reconstruction.}
\label{fig:kuramoto-modal-noise-complex-angles}
\end{figure}

\begin{figure}[htbp]
\begin{center}
\includegraphics[width = 0.45\textwidth, height=0.15\textheight]{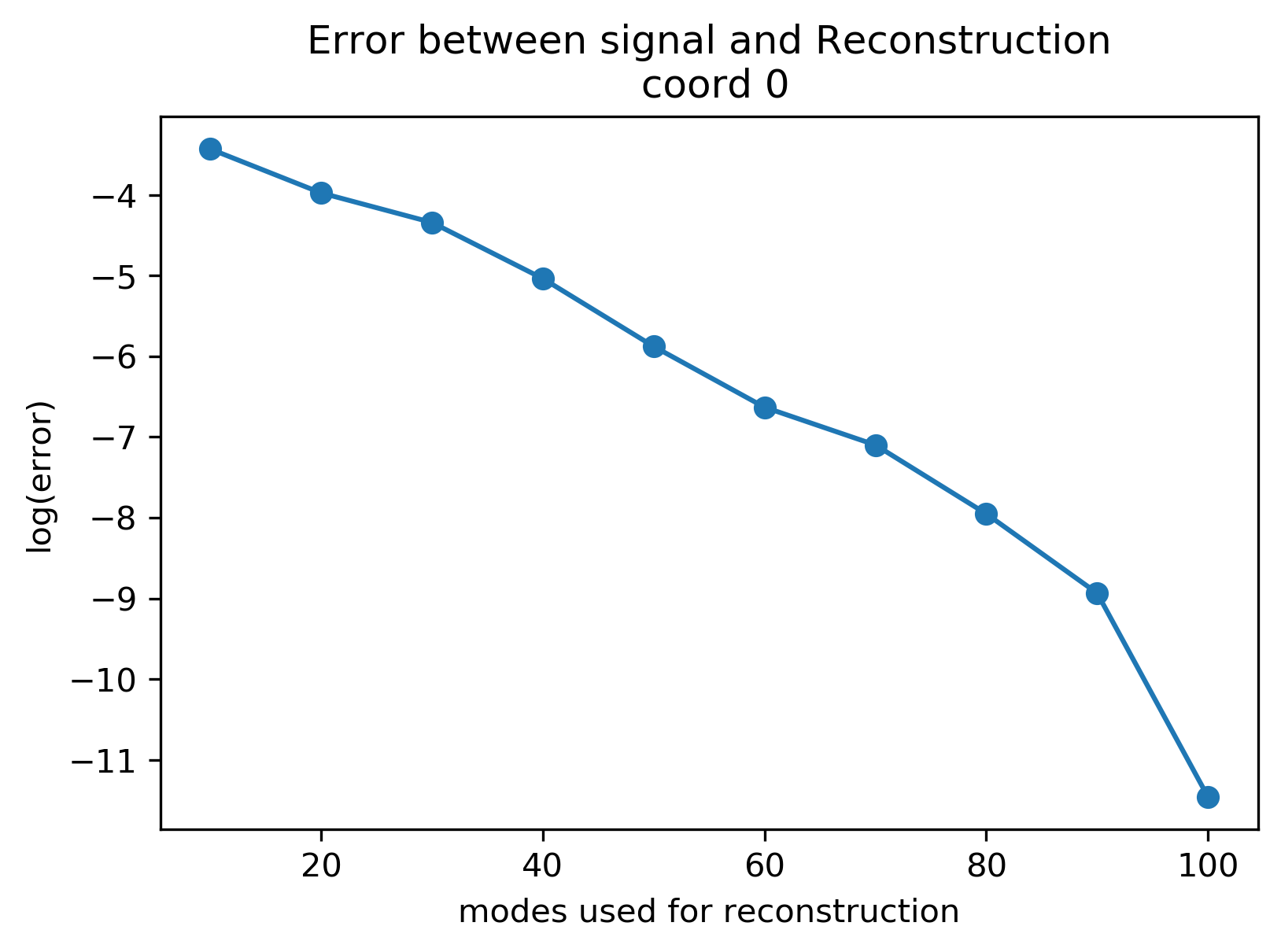}
\includegraphics[width = 0.45\textwidth, height=0.15\textheight]{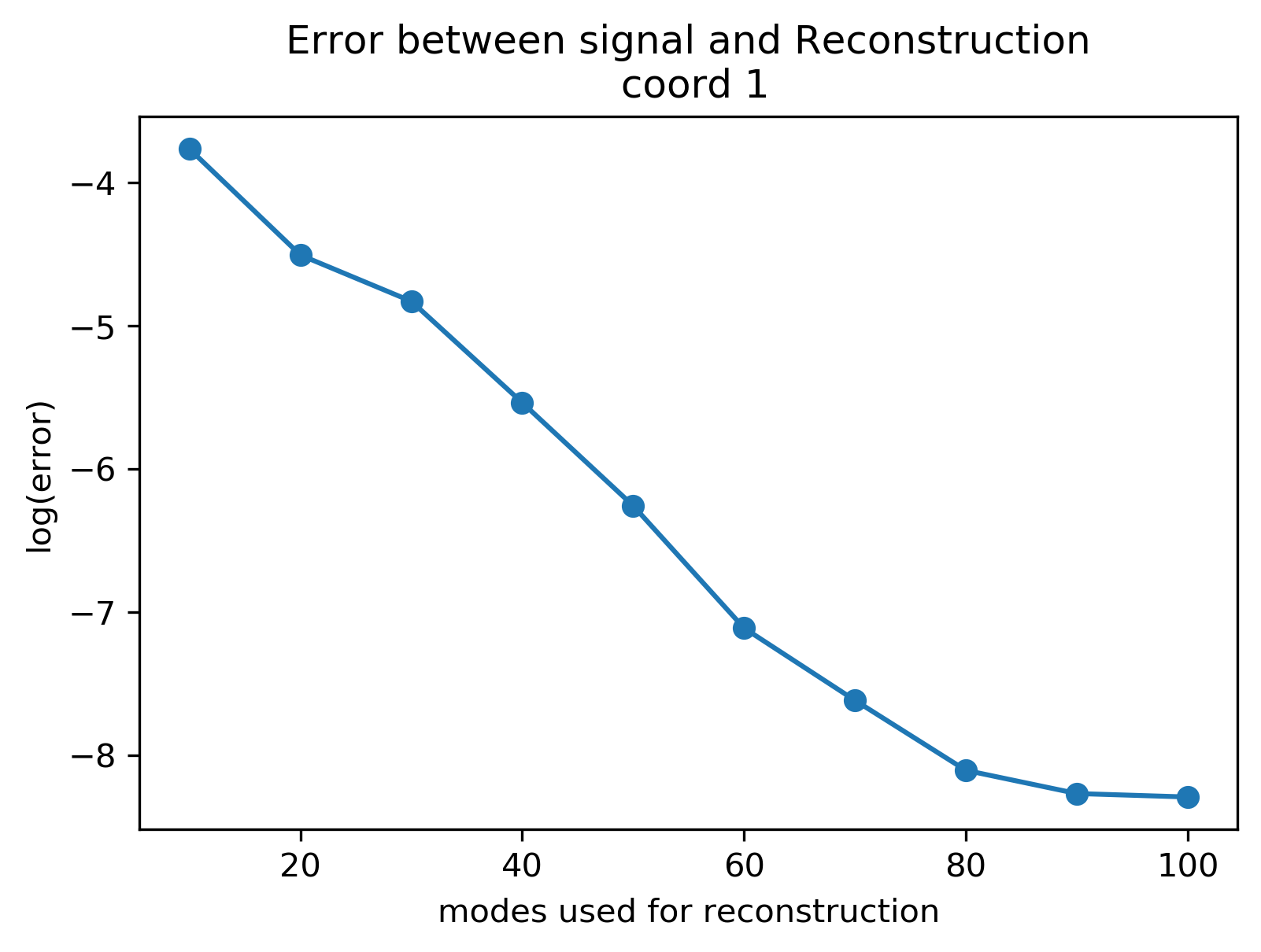}
\includegraphics[width = 0.45\textwidth, height=0.15\textheight]{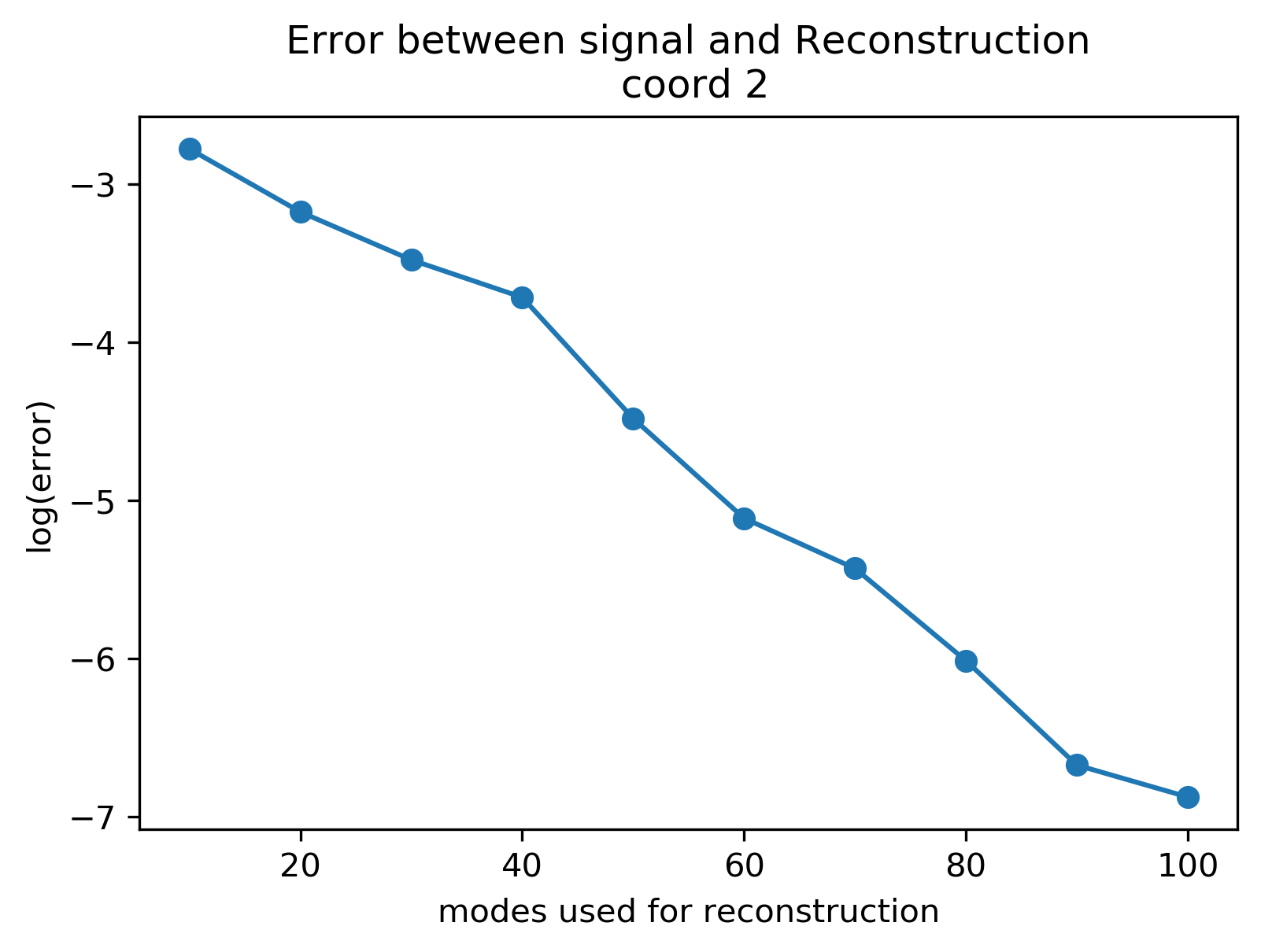}
\includegraphics[width = 0.45\textwidth, height=0.15\textheight]{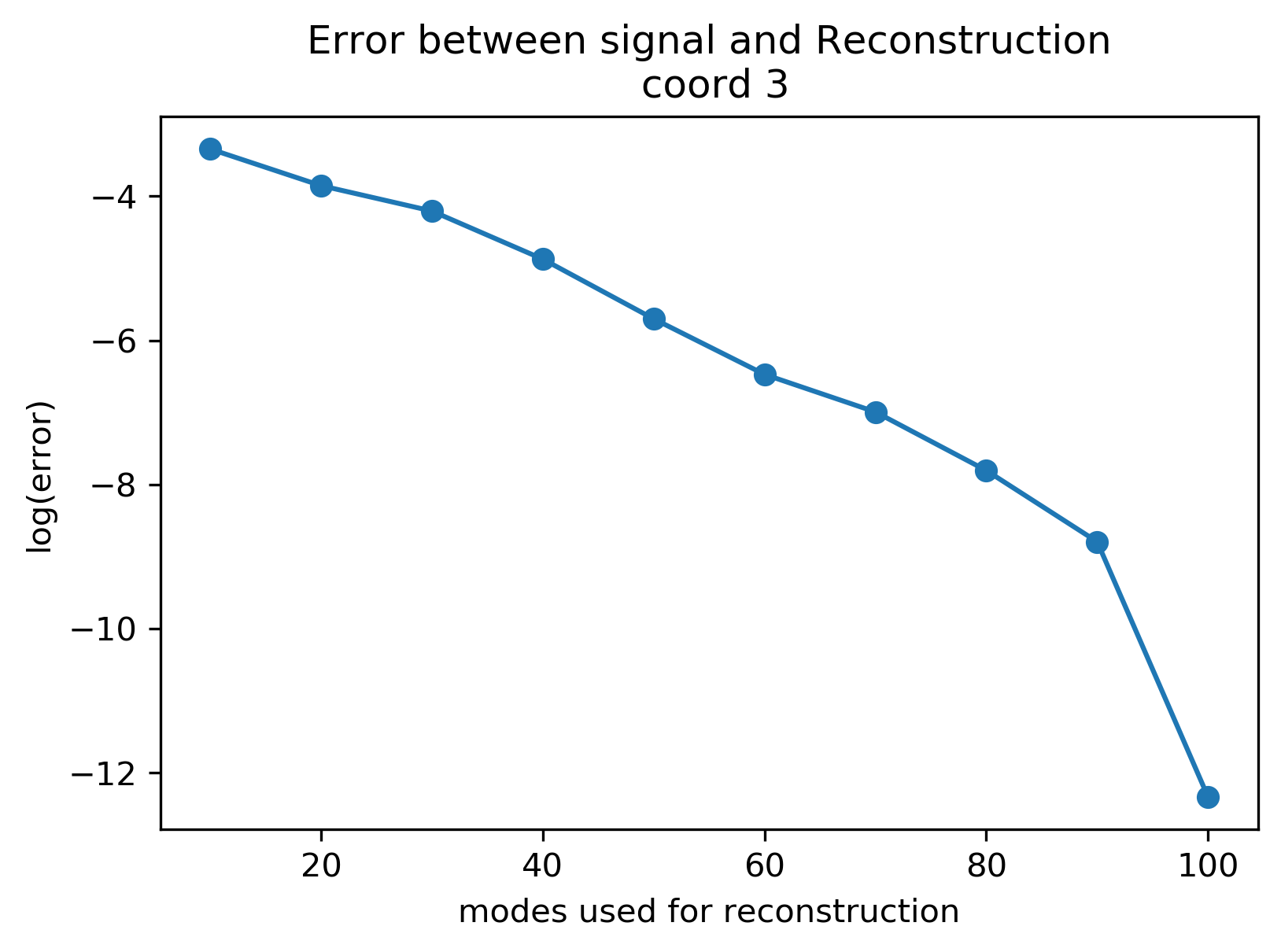}
\includegraphics[width = 0.45\textwidth, height=0.15\textheight]{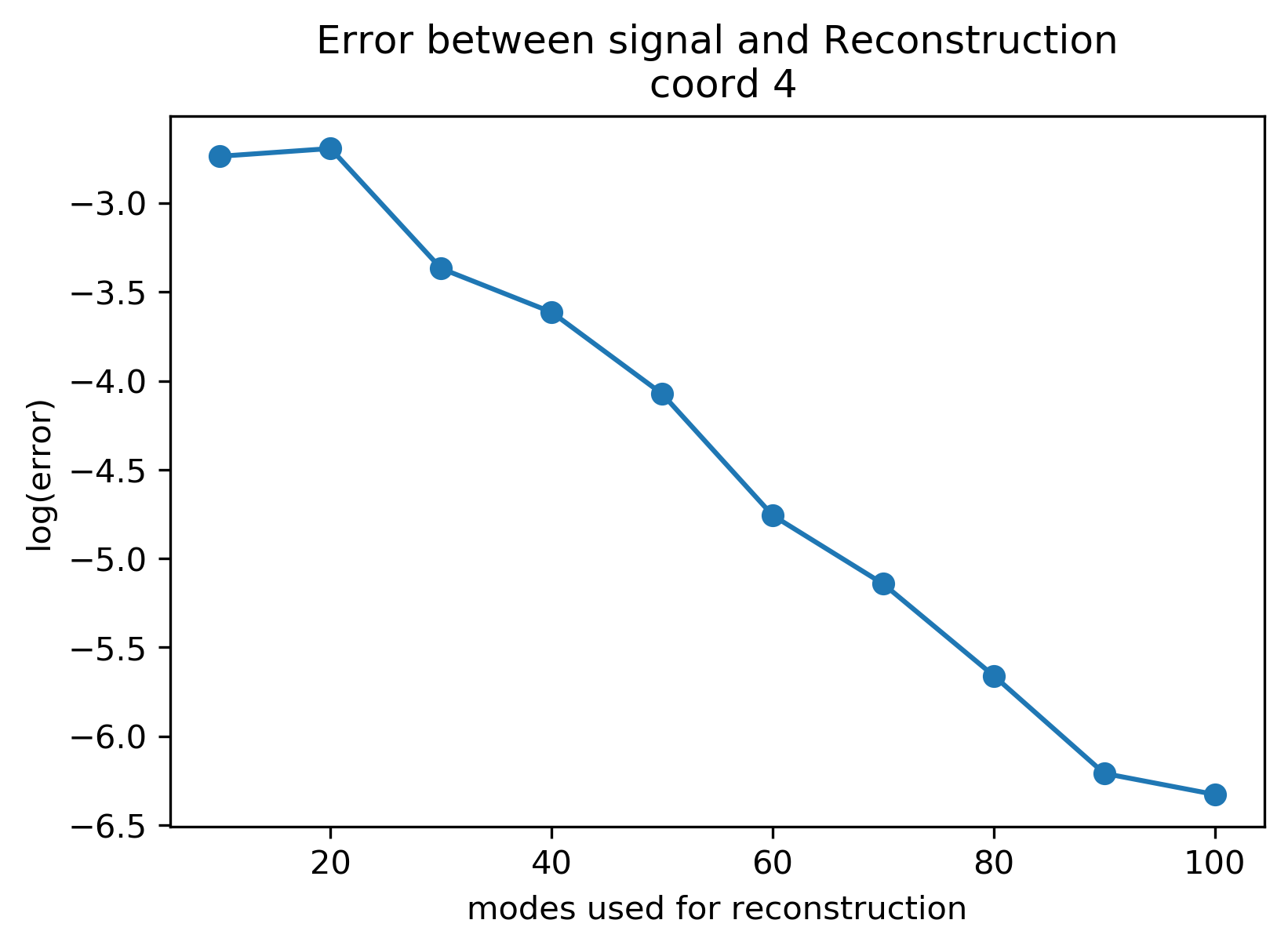}
\includegraphics[width = 0.45\textwidth, height=0.15\textheight]{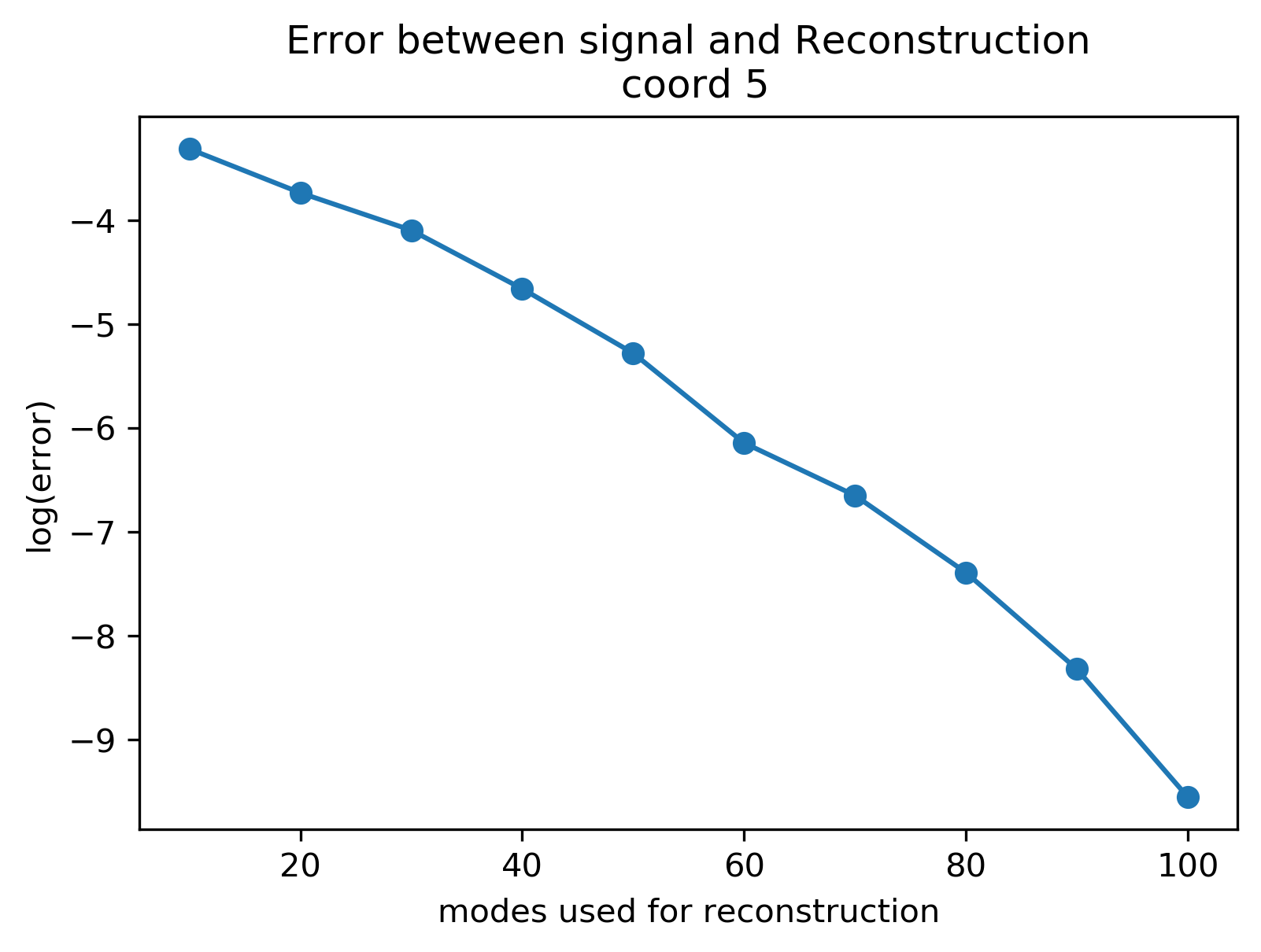}
\includegraphics[width = 0.45\textwidth, height=0.15\textheight]{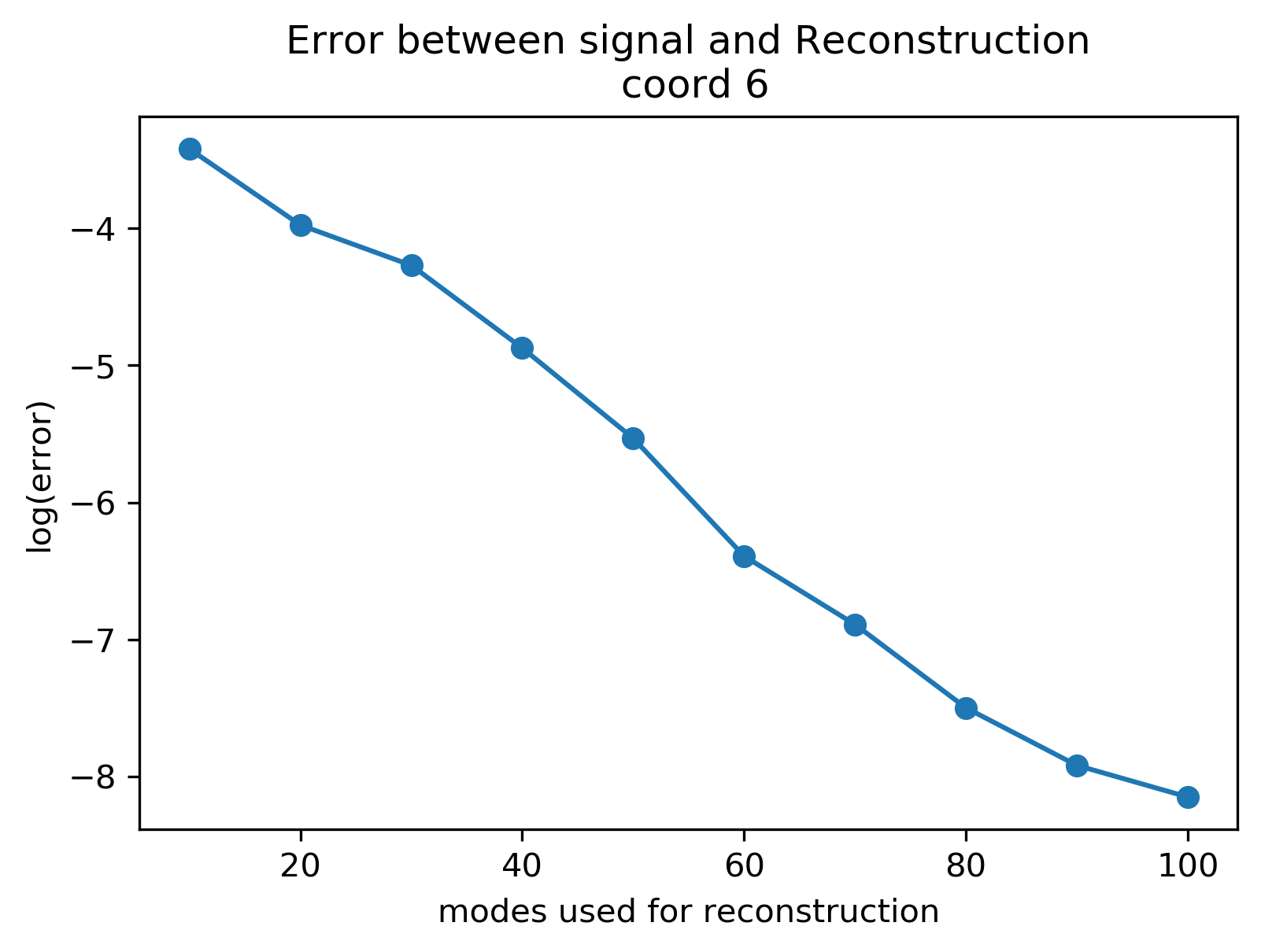}
\includegraphics[width = 0.45\textwidth, height=0.15\textheight]{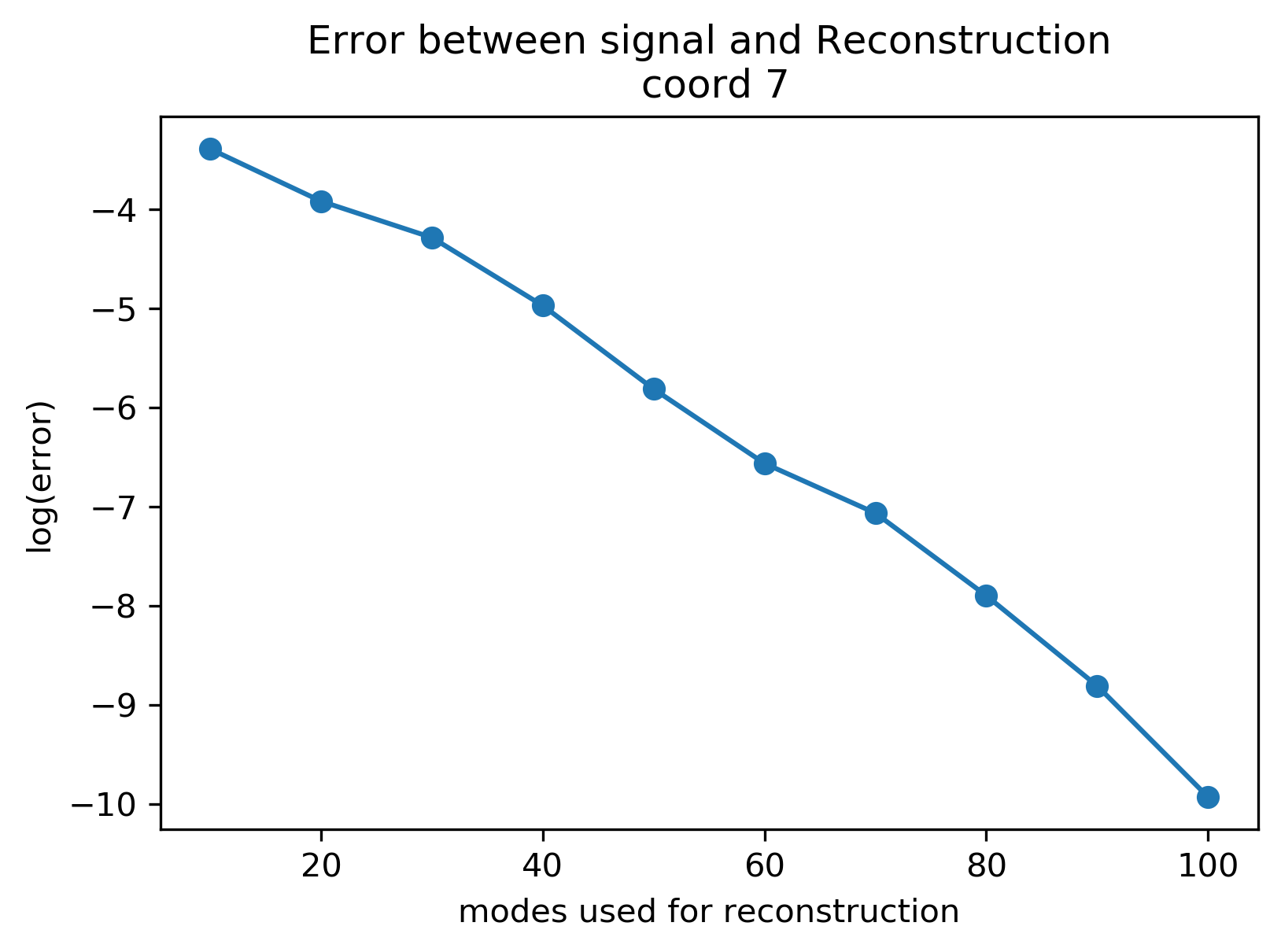}
\includegraphics[width = 0.45\textwidth, height=0.15\textheight]{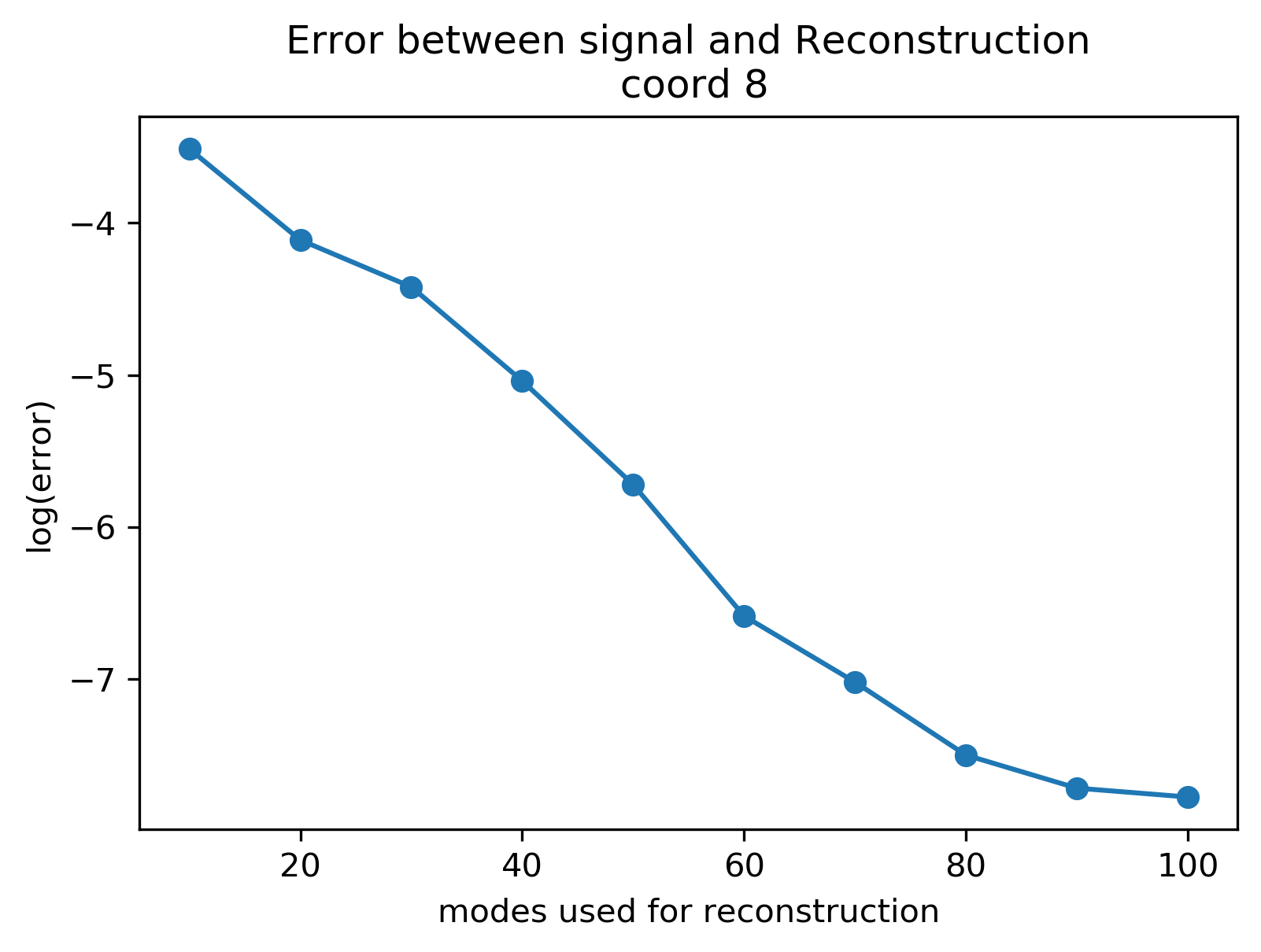}
\includegraphics[width = 0.45\textwidth, height=0.15\textheight]{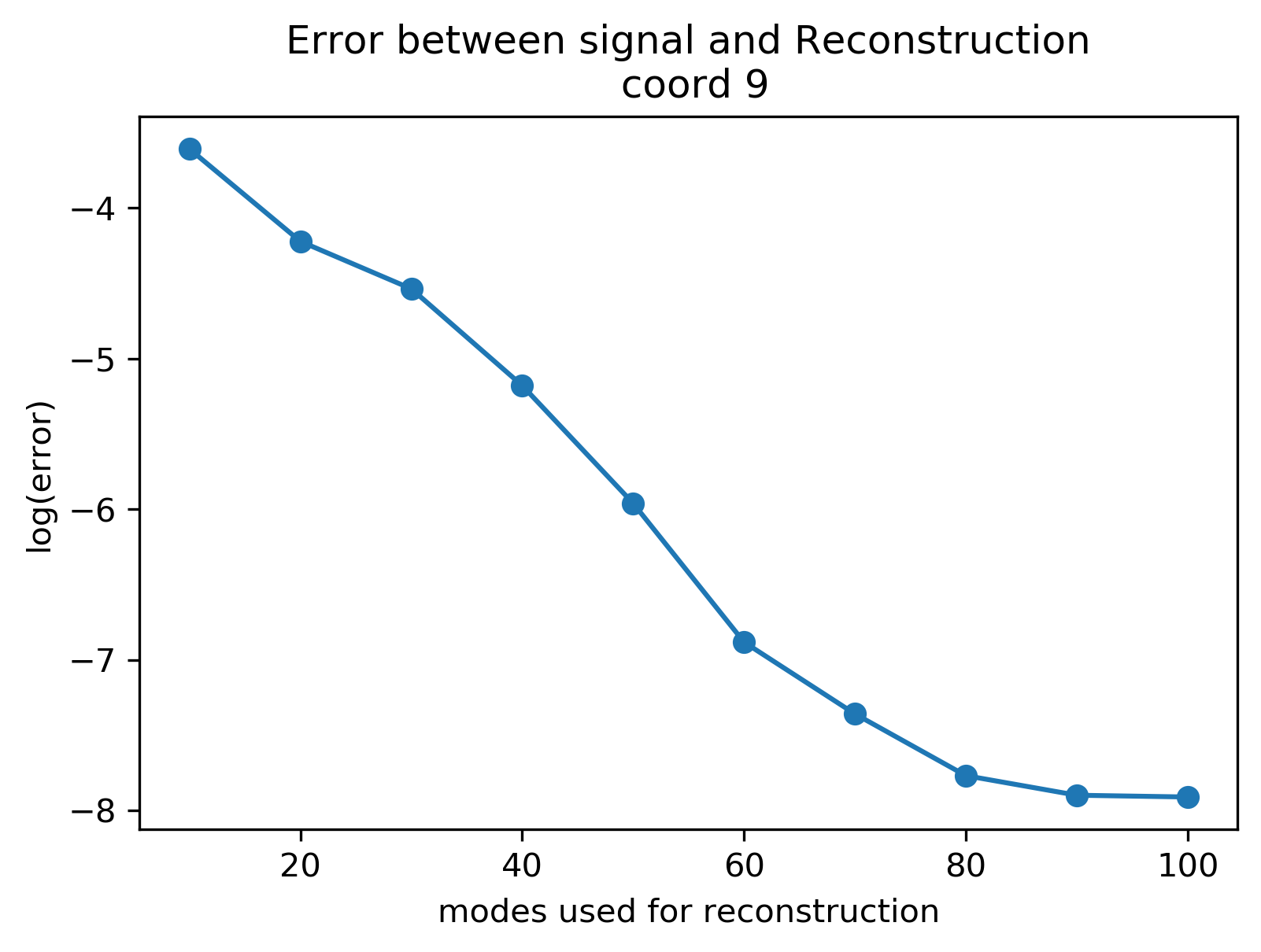}
\end{center}
\caption{\textbf{Kuramoto models, signal vs reconstruction error}: Comparison of error between the reconstruction signal and the true signal for each oscillator. The reconstruction error is computed using \eqref{eq:geodesic-dist}.}
\label{fig:kuramoto-reconstruction-error-complex-angles}
\end{figure}

\begin{figure}[htbp]
\begin{center}
\includegraphics[width = 0.45\textwidth, height=0.15\textheight]{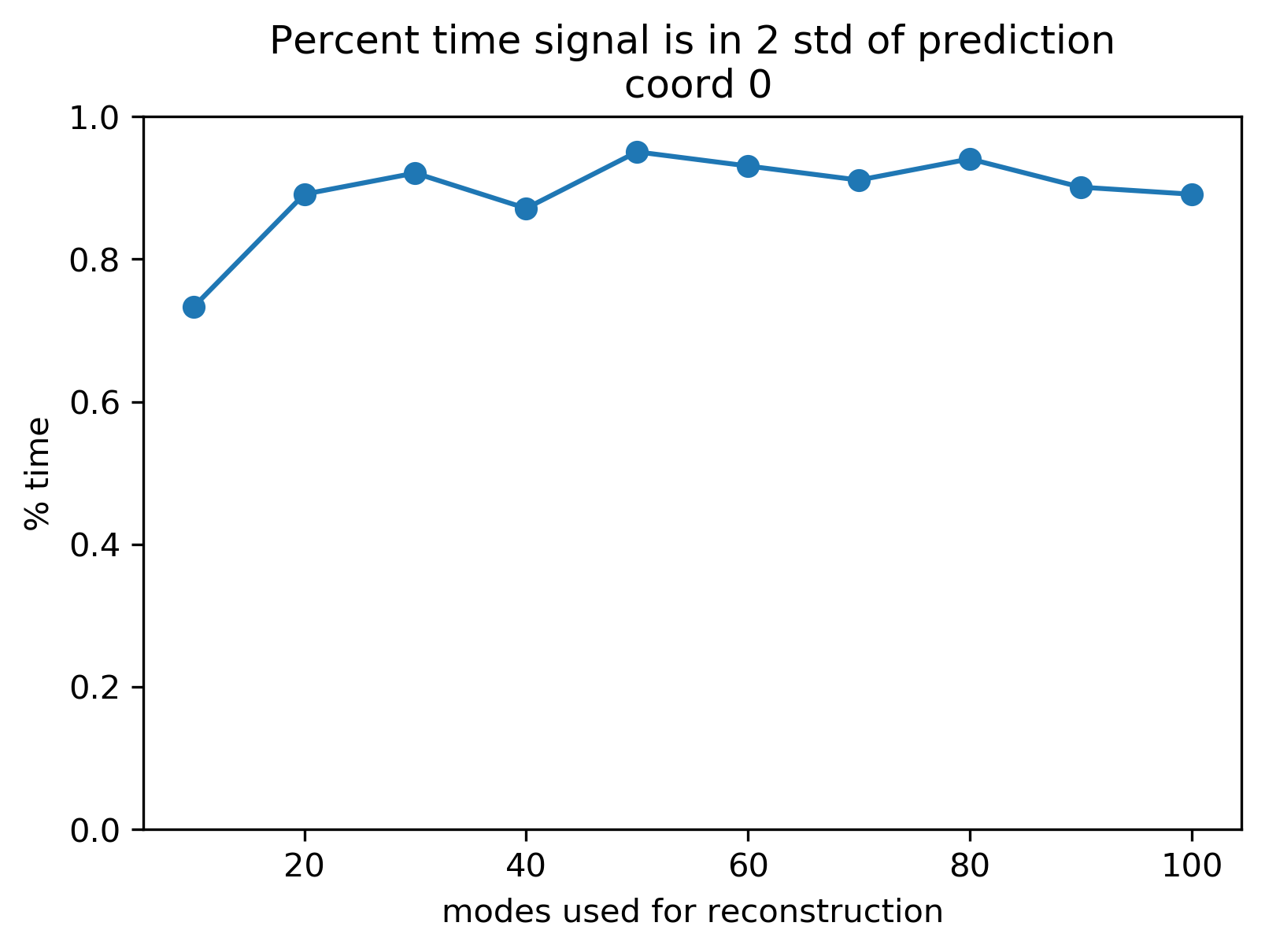}
\includegraphics[width = 0.45\textwidth, height=0.15\textheight]{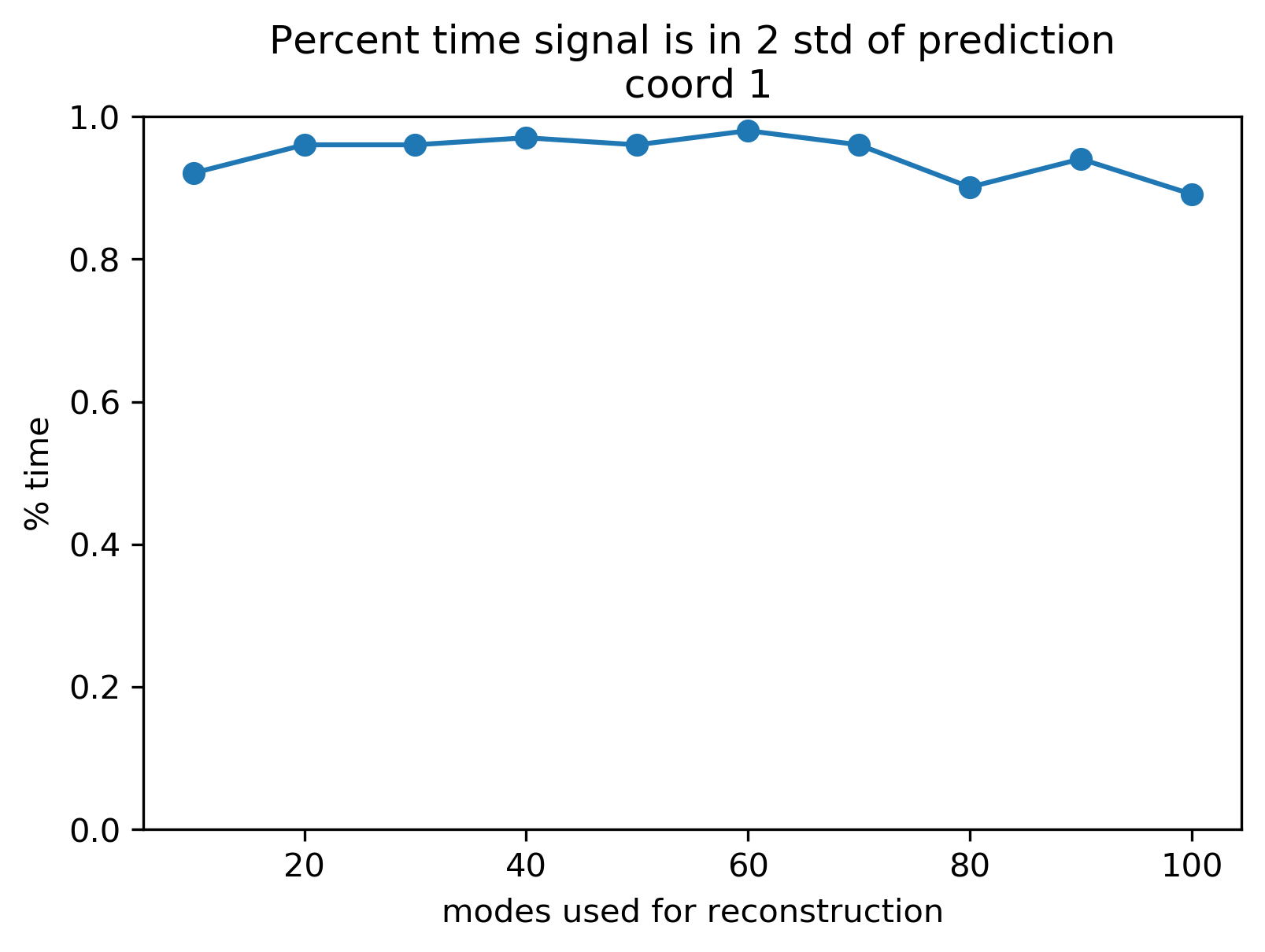}
\includegraphics[width = 0.45\textwidth, height=0.15\textheight]{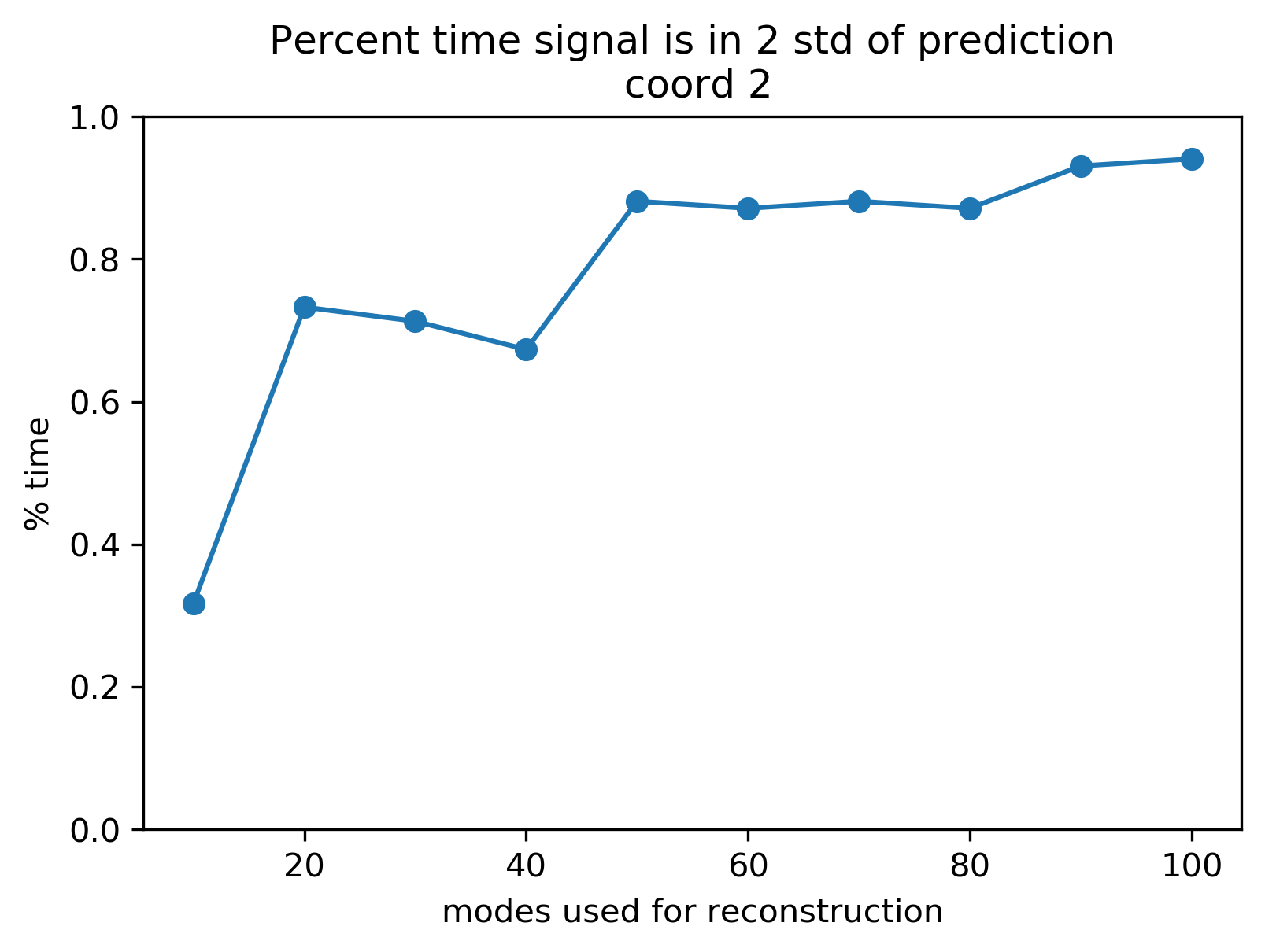}
\includegraphics[width = 0.45\textwidth, height=0.15\textheight]{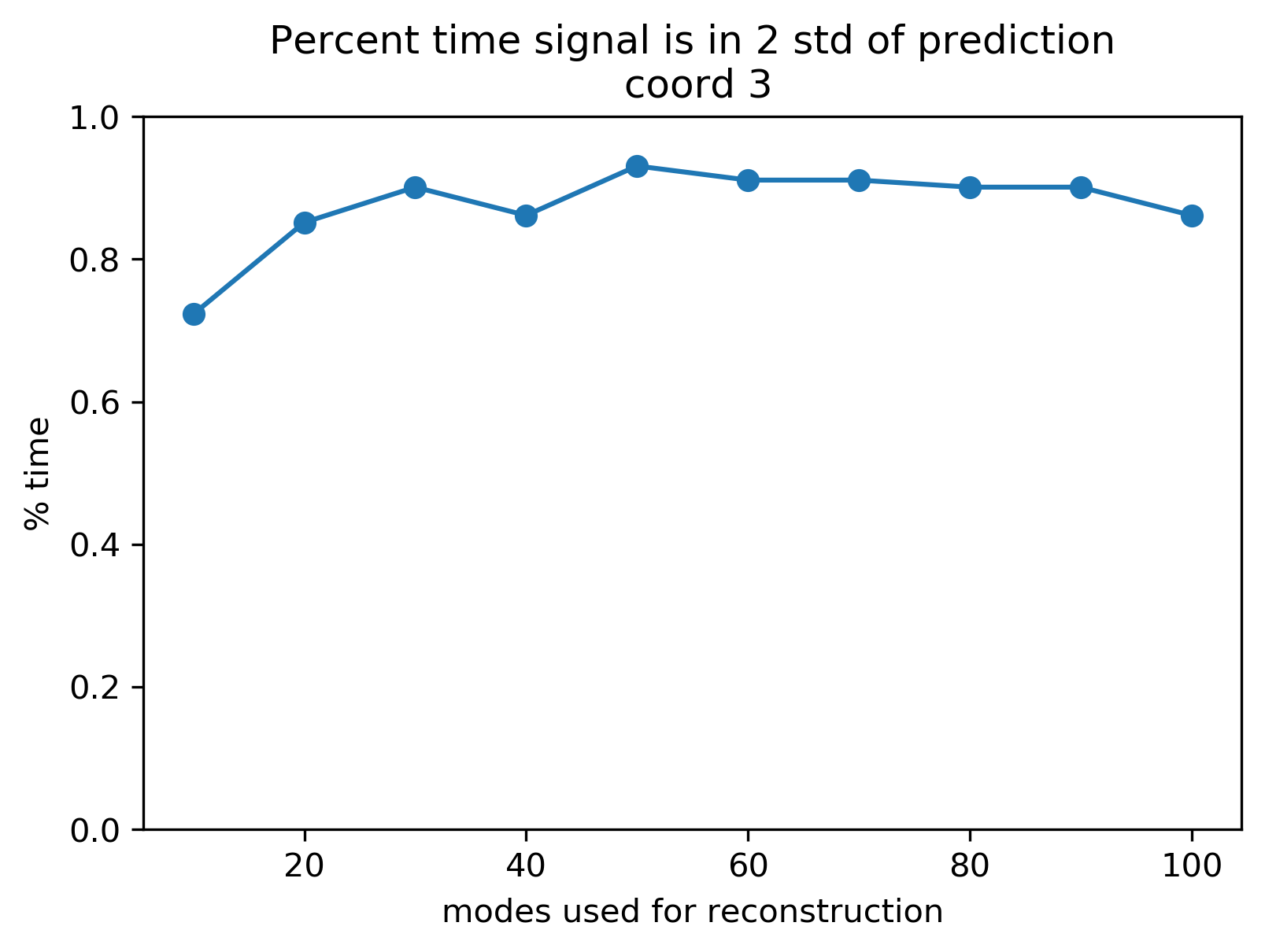}
\includegraphics[width = 0.45\textwidth, height=0.15\textheight]{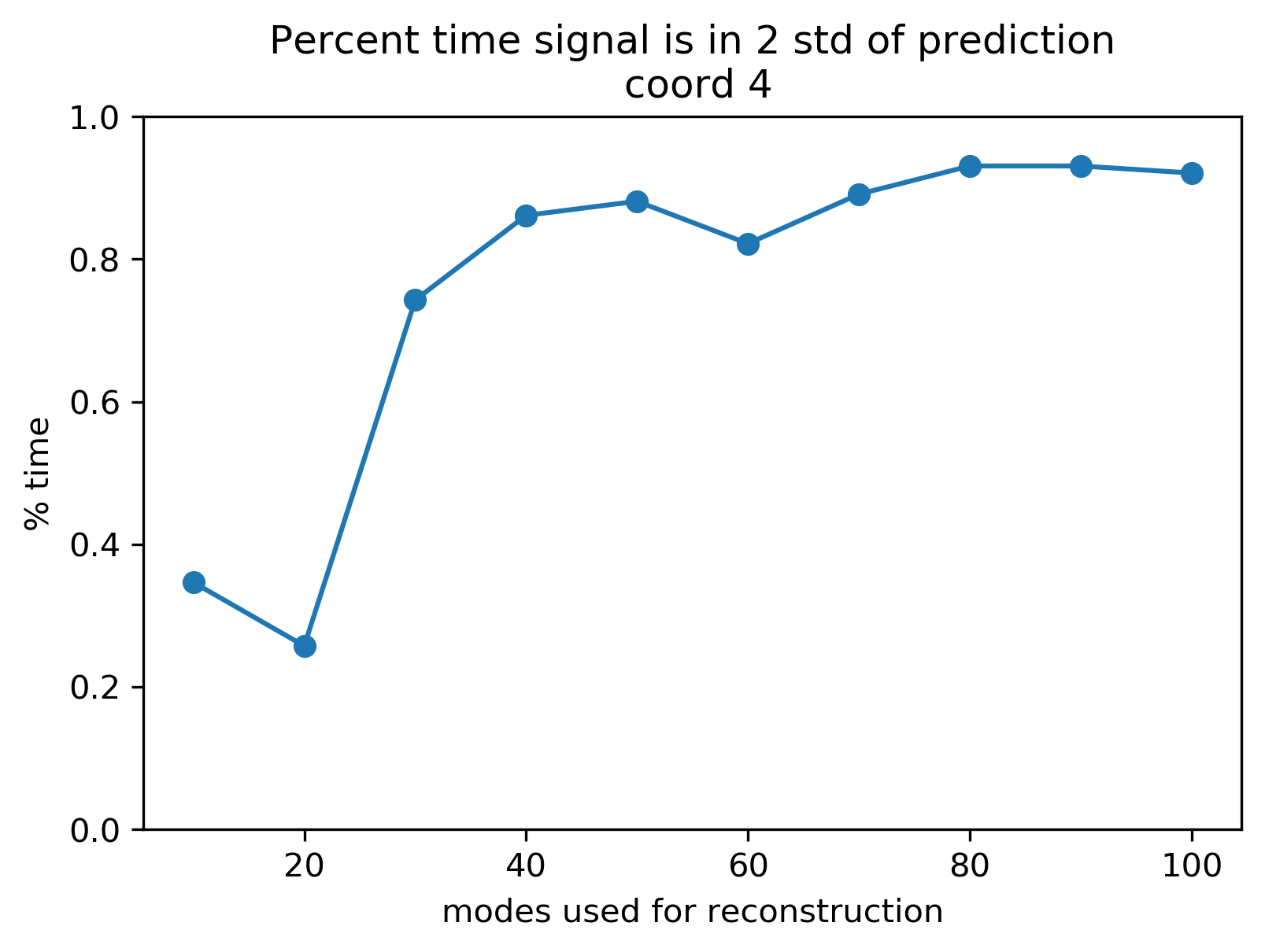}
\includegraphics[width = 0.45\textwidth, height=0.15\textheight]{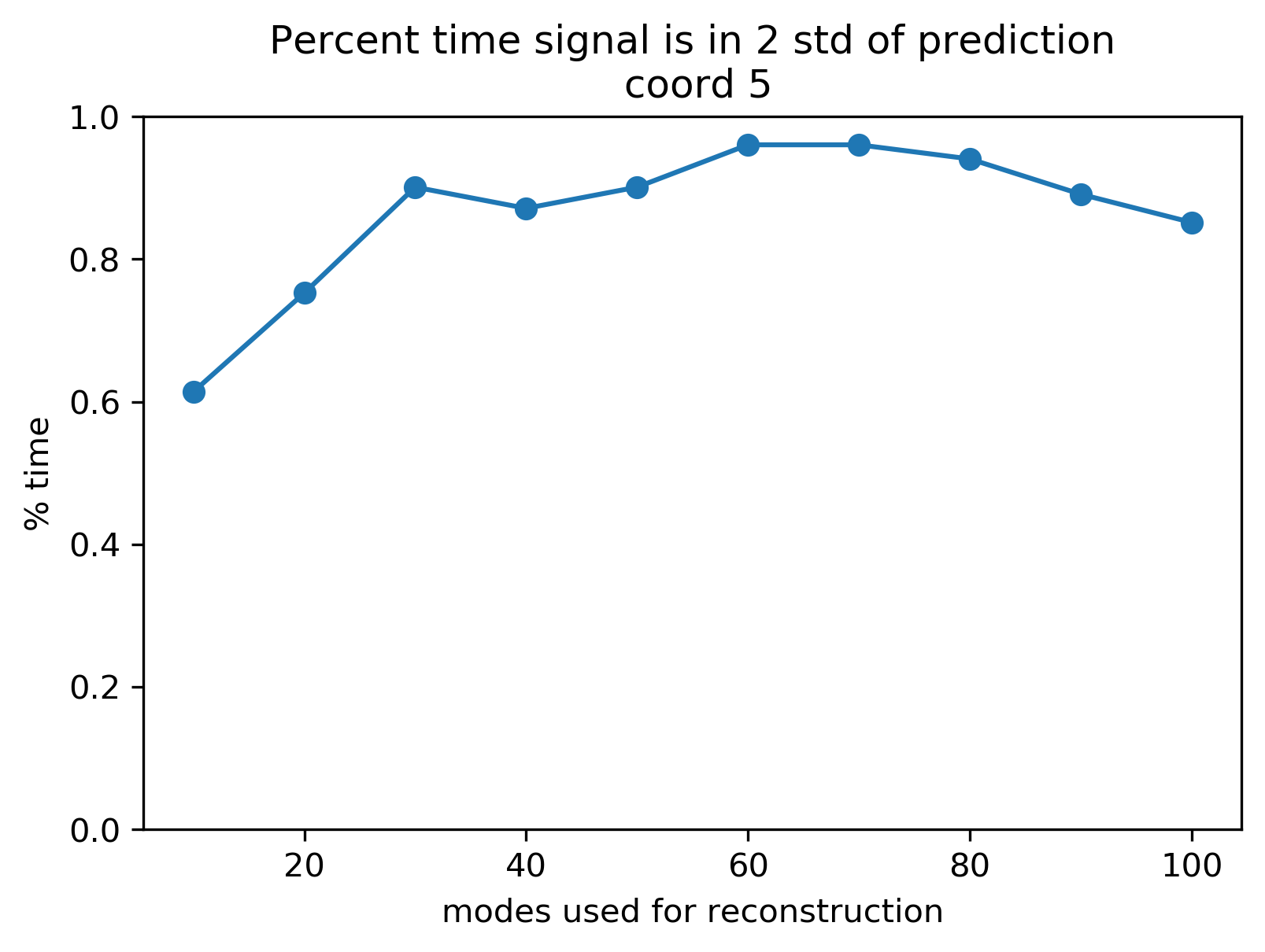}
\includegraphics[width = 0.45\textwidth, height=0.15\textheight]{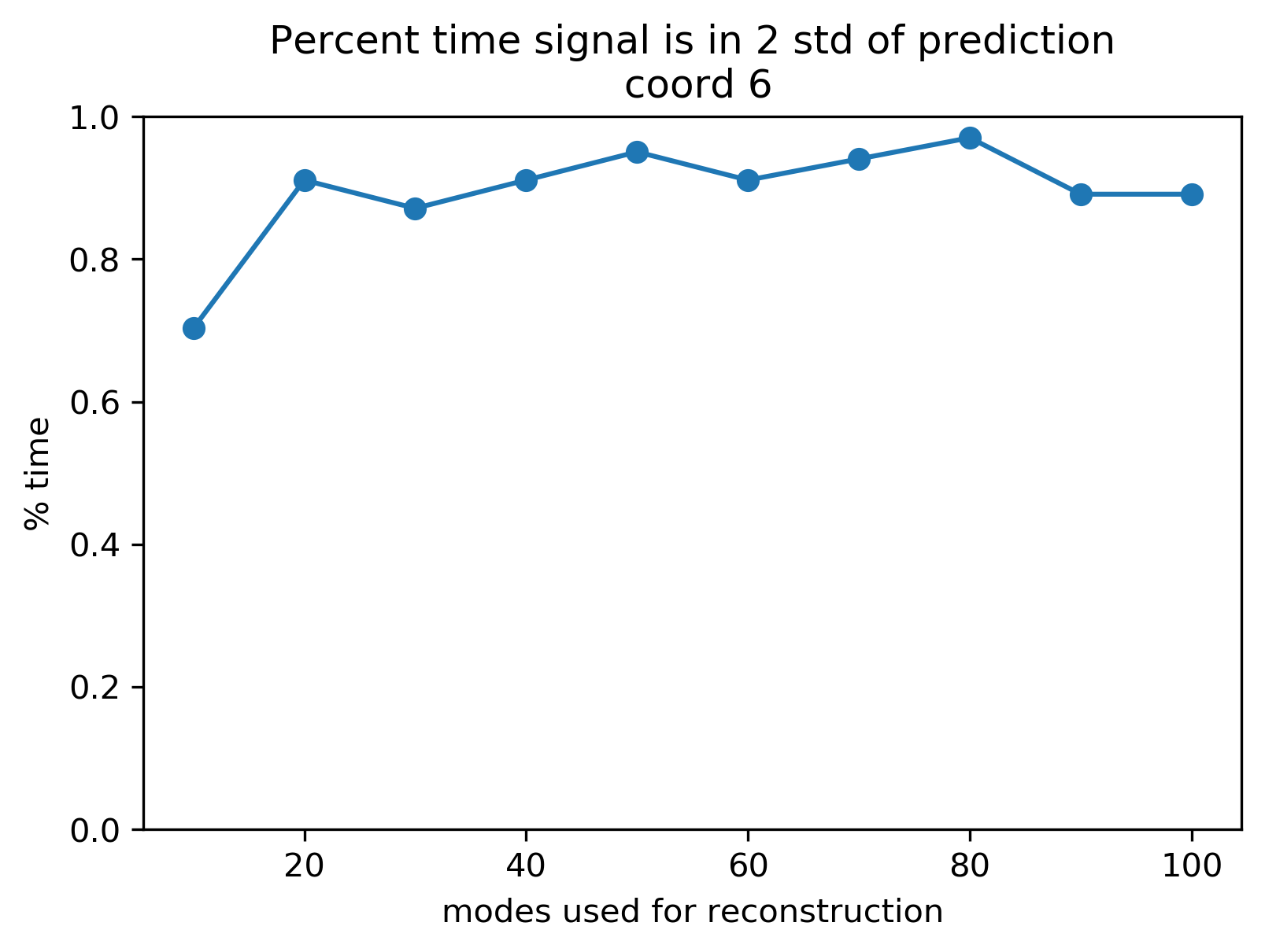}
\includegraphics[width = 0.45\textwidth, height=0.15\textheight]{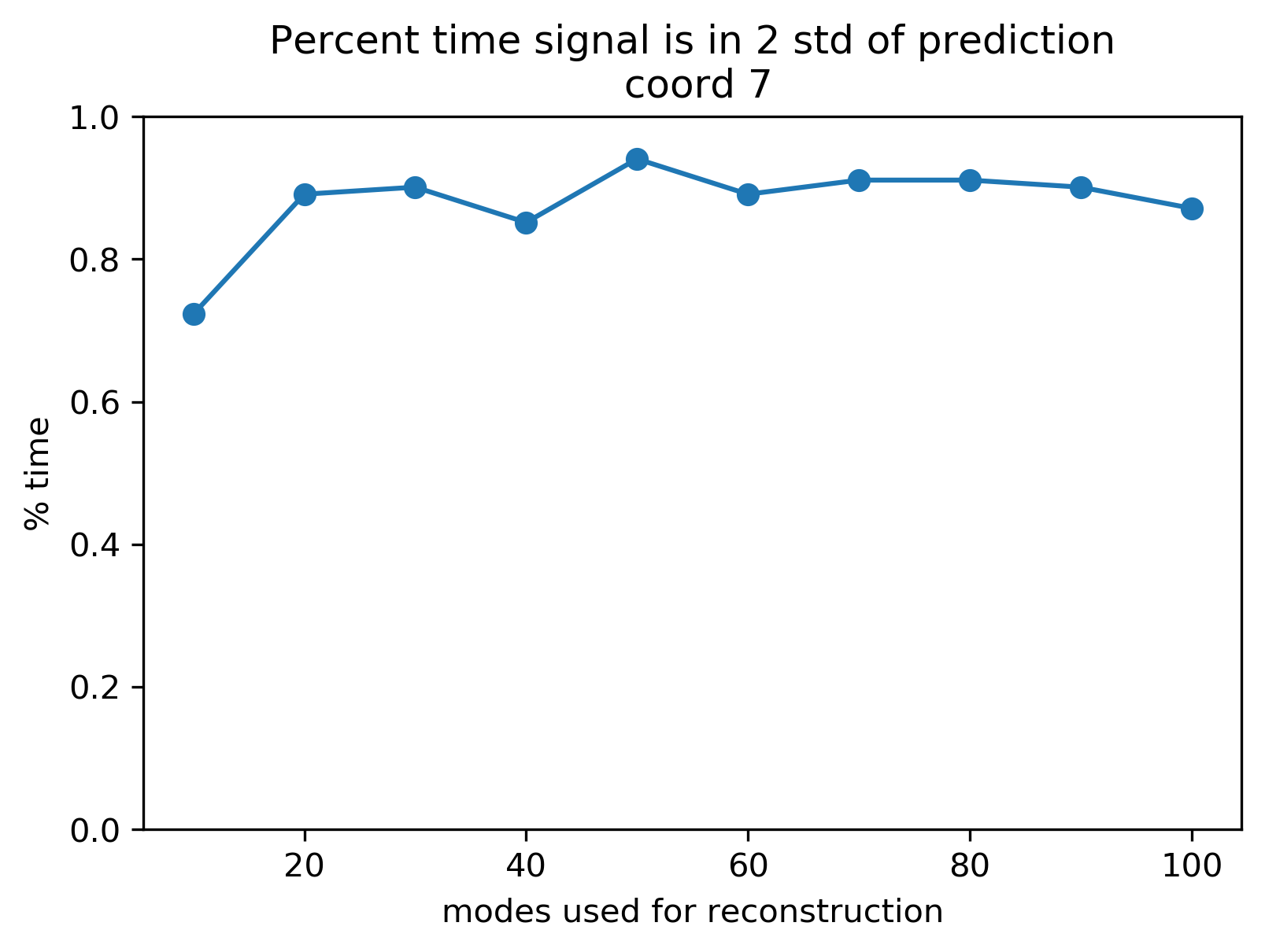}
\includegraphics[width = 0.45\textwidth, height=0.15\textheight]{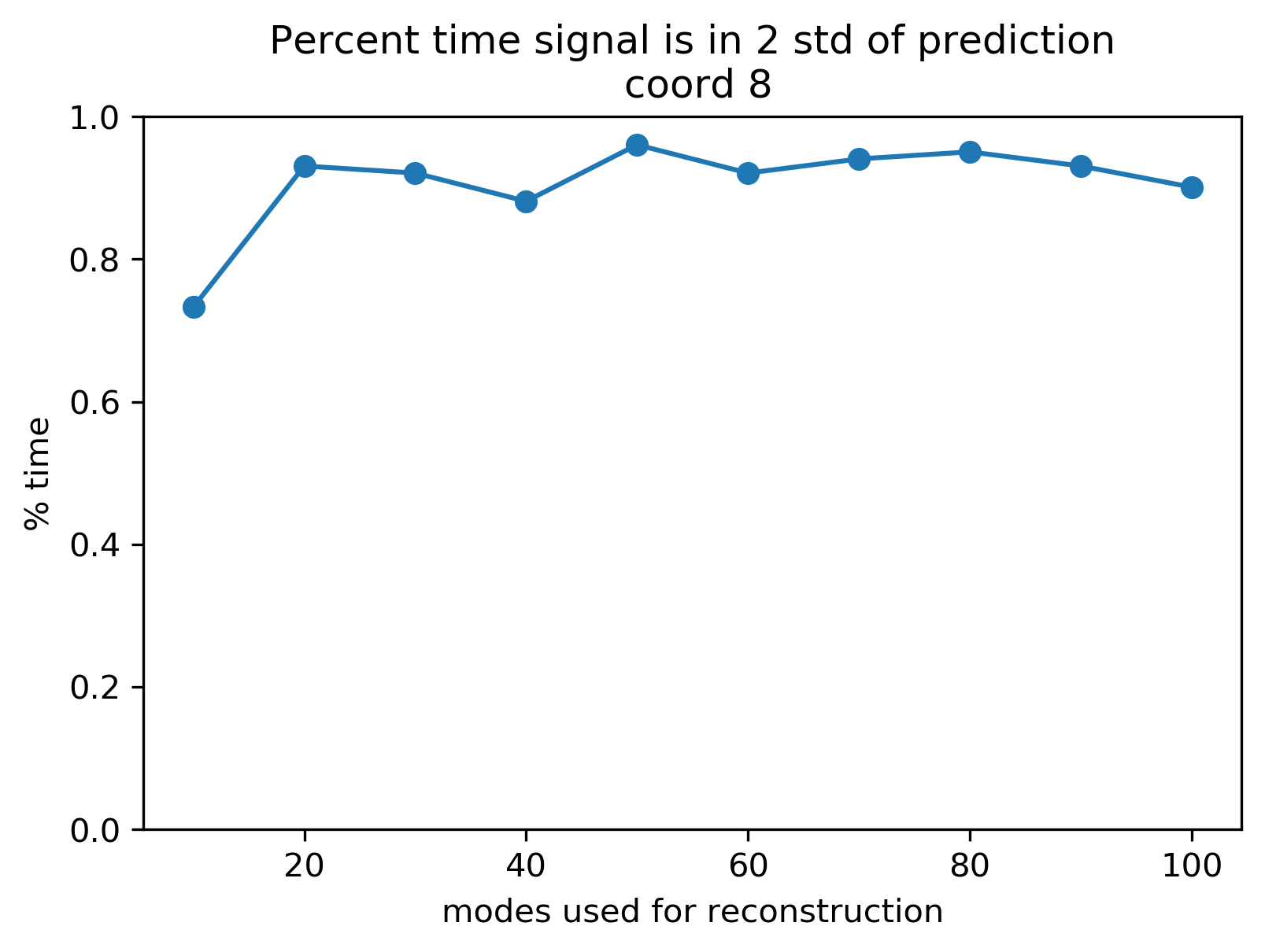}
\includegraphics[width = 0.45\textwidth, height=0.15\textheight]{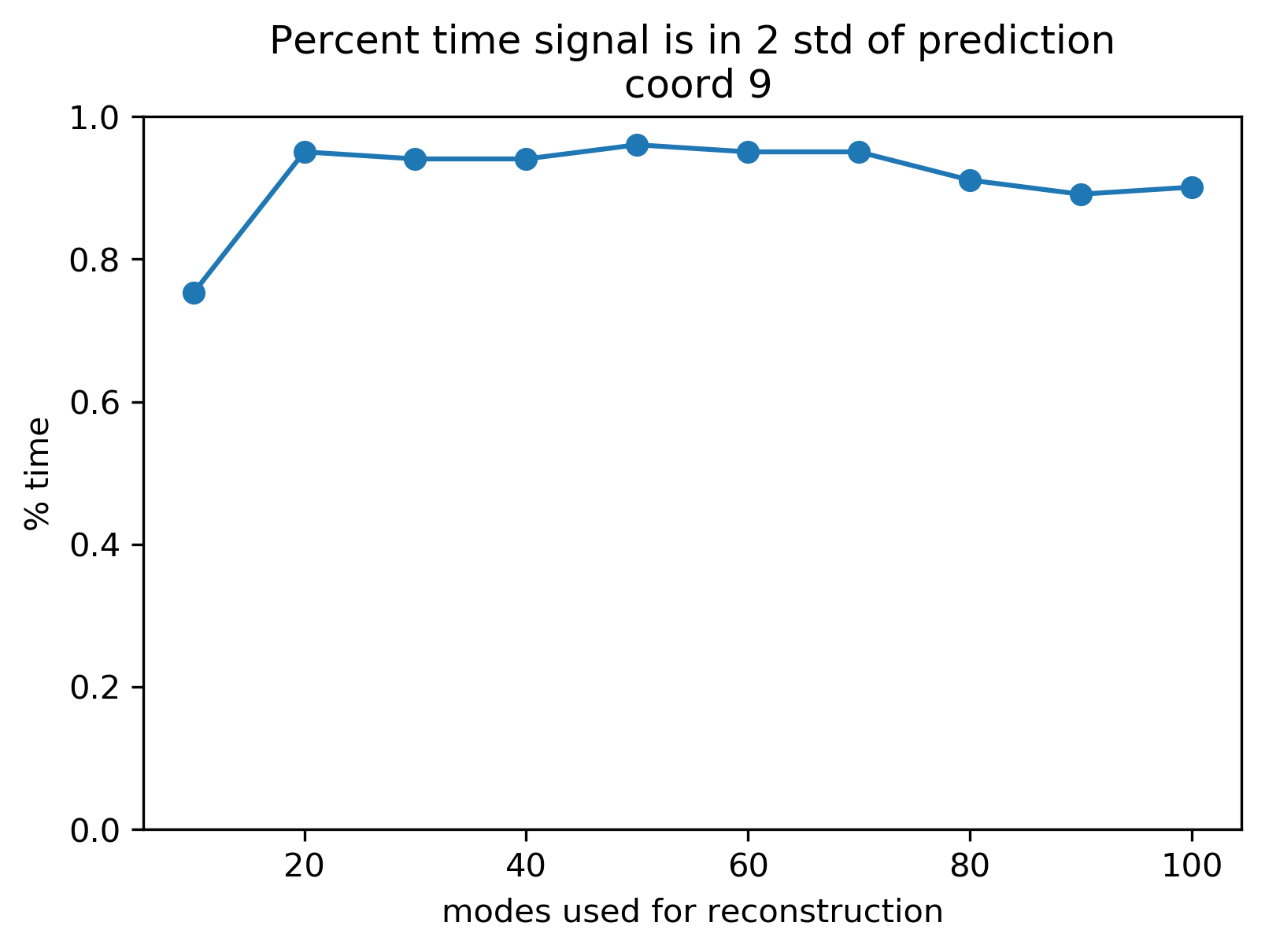}
\end{center}
\caption{\textbf{Kuramoto models, fraction of time true signal is within $\pm 2$ standard deviations of ROM signal}: Residence times in $\pm 2$ standard deviation band is computed via \eqref{eq:2std-residence-time}. A residence time of 1.0 corresponds to the true signal always being in the confidence bounds of the ROM.}
\label{fig:kuramoto-residence-time-complex-angles}
\end{figure}

\blue{
Figures \ref{fig:kuramoto-eigenvalues-complex-angles} - \ref{fig:kuramoto-residence-time-complex-angles} show an ablative study similar to the anharmonic oscillator example. Figure \ref{fig:kuramoto-reconstruction-complex-angles} shows the reconstruction of the signal of the first oscillator. The ROMs predictions converge rapidly to the true signal as the number of modes increases. Additionally, the computed variance rapidly converges to 0. In fact, one cannot distinguish a $\pm 2$ standard deviation band when using 60 modes or greater. 
}


We quantify the \blue{reduced order models in the ablation study} using two metrics. The first is the reconstruction error. This is simply the average ``geodesic'' distance between the trajectories. If the true trajectory is $\Theta = \set{\theta(t) : t \in [0, T], \theta(t) \in [0, 2\pi)}$ and the ROM trajectory is $\Phi = \set{\phi(t) : t \in [0, 1, \dots, T], \phi(t) \in [0, 2\pi)}$, then the \blue{average } ``geodesic'' distance between the trajectories is computed as 
	\begin{equation}\label{eq:geodesic-dist}
	d(\Theta, \Phi) = \frac{1}{T+1} \left(\sum_{t=0}^{T} \abs{e^{i\theta(t)} - e^{i\phi(t)}}^2 \right)^{1/2}.
	\end{equation}
\blue{Note that \eqref{eq:geodesic-dist} does not use the true geodesic distance on $[0, 2\pi) \mod 2\pi$ for the the averaging, where the true geodesic distance is} computed as
	\begin{equation}
	d_g(\theta(t), \phi(t)) = \min_{k\in \Z} \abs{\theta + (2\pi k + \phi)}.
	\end{equation}
However, it is a good approximation \blue{in the sense that $\abs{e^{i\theta(t)} - e^{i\phi(t)}}^2$ is a continuous, monotonically increasing function with respect to $d_g(\theta(t), \phi(t)))$.}

Figure \ref{fig:kuramoto-reconstruction-error-complex-angles} shows the reconstruction log error of the ROM's vs.\ the number of modes used for reconstruction. While the log error between the signal and the reduced order model is a good metric, we also have access to the modal noise which we use to put error bounds on the predictions given by the ROM. We compute the \blue{fraction} of time the true trajectory is within $\pm 2$ standard deviations of the ROM's prediction. This is computed as 
	\begin{equation}\label{eq:2std-residence-time}
	P = \frac{1}{T+1} \sum_{t=0}^{T} \ind{[0, 2\sigma)}( d_g(\theta(t), \phi(t)) )
	\end{equation}
where $d_g(a, b)$ is the \emph{true}\footnote{As opposed to the approximate geodesic distance \eqref{eq:geodesic-dist}} geodesic distance on $[0, 2\pi) \mod 2\pi$ and $\sigma$ is the modal noise standard deviation. Figure \ref{fig:kuramoto-residence-time-complex-angles} shows the fraction of time that the true signal is within 2 standard deviations of the ROM prediction. 

\blue{Again we can use the heuristic in section \ref{sec:heuristic} to choose how many modes to keep in the deterministic part of the ROM. According to the heuristic, 30 modes is the minimum number of modes to keep (see figure \ref{fig:heuristic-kuramoto}).}

\clearpage
\section{A heuristic for estimating the minimum number of modes to use for the ROM}\label{sec:heuristic}
Above, we have computed an ablation study of reduced order models for 3 examples, from a minimum of 10 modes up to a full reconstruction (100 modes in all examples). The question remains how to choose the final number of modes to be used in the reduced order model (ROM); too few and we get poor prediction accuracy because we are not accurately capturing the nominal dynamics, too many modes and we are trying to capture noise and the continuous spectrum with eigenvalues and eigenfunctions.

In this section, we present a heuristic method to determine the minimum number of modes that must be taken.\blue{This heuristic is motivated by plots that show that by increasing the number of modes, the the distributions seem to become Gaussian.} This method is based on testing the normality of the modal distribution that is computed for each ROM. The rationale is as follows. In our examples, we have structured dynamics plus added noise and a changing network topology. When we construct our model, we compute a deterministic model and then estimate the modal noise distribution. We assume that once we have adequately recovered the true structured dynamics with the computed deterministic part of the ROM the modal noise will have a Gaussian distribution. \blue{This can happen due to a variety of settings: the deterministic chaotic part has hyperbolic dynamics, or the stochastic elements are Gaussian distributed, or the number of dimensions is very large and projection yields the central limit theorem.} If too few modes are used for the deterministic model, much of the structured part of the dynamics will be modeled as noise, destroying the normality of the modal distribution. By using too many modes in the deterministic model, the noisy part of the dynamics are subsumed into a deterministic model, possibly destroying the normality of the distribution. \blue{Since we are fitting a Gaussian to the distributions, the non-normality of the distribution can lead to inaccurate estimates for the confidence bounds, thereby degrading the model's predictions.}

Algorithm \ref{alg:heuristic} specifies the heuristic. We note that this test can be computed automatically and does not require manually looking at histograms or quartile-quartile plots to determine this although we show such plots in Figures \ref{fig:heuristic-linear}, \ref{fig:heuristic-kuramoto}, and \ref{fig:heuristic-anharmonic} below. In each of the figures, panel (a) shows both the mean and median p-values vs. the number of modes. Panel (b) shows a box a whisker plot of the p-values to show their distribution. Panel (c) shows a quartile-quartile plot; the closer the blue dots follow the red line, the closer the distribution is to a normal distribution. Final, panel (d) shows the reduced order model reconstruction using the number of modes dictated by the heuristic algorithm.

\begin{algorithm}[H]\label{alg:heuristic}
\DontPrintSemicolon
  
  \KwData{The set modal noise distributions $\rho_i$ for model $i$.}
  \KwResult{Minimum number of modes to use for ROM.} 
  
  Let $\rho_i = \set{\rho_{i,j}}_{j=1}^{J}$ be the set of modal distributions for reduced order model $i$. $\rho_{i, j}$ is the modal distribution for the $j$-th coordinate of the system. The $\rho_i$'s should be ordered from smallest model to largest (e.g., $i = 10, 20, \dots, 100$ in our examples.)
  
  For $j=1,\dots, J$, perform a Shapiro-Wilk hypothesis test to compute $J$ p-values, $p_i = \set{p_{i,1}, \dots, p_{i,J}}$ (one for each modal distribution).
  
  For each $i$, compute the mean and median of the set $p_i$.
  
  \textbf{Heuristic:} The smallest ROM where both the mean and median are above 0.05 is the smallest number of modes one should take for the model.

\caption{Heuristic for determining the minimum number of modes used in the ROM.}
\end{algorithm}

\begin{figure}[htbp]
\centering
\begin{subfigure}[b]{0.45\textwidth}
    \centering
    \includegraphics[width = \textwidth]{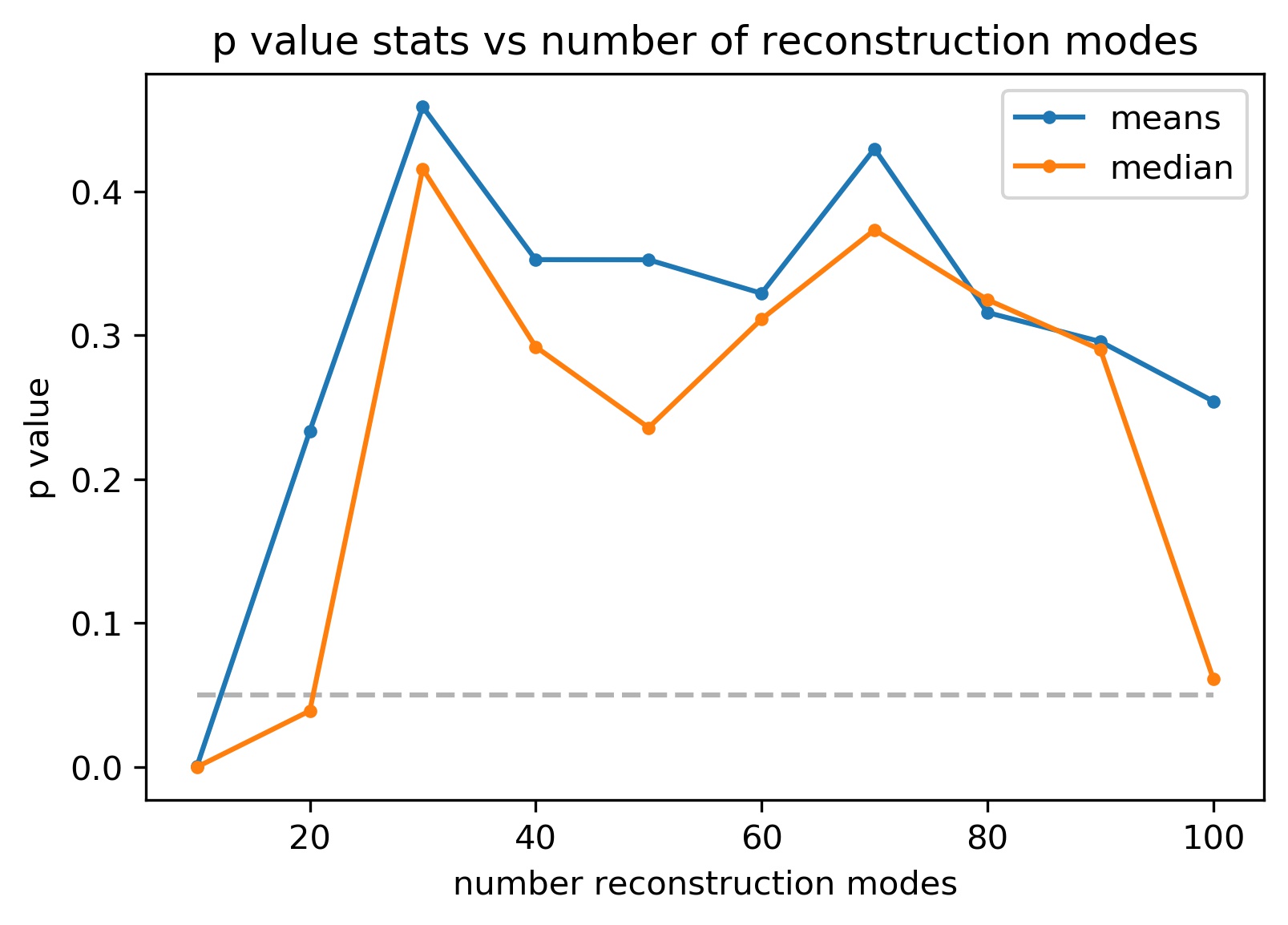}
    \caption{p values mean, median}
\end{subfigure}
\begin{subfigure}[b]{0.45\textwidth}
    \centering
    \includegraphics[width = \textwidth]{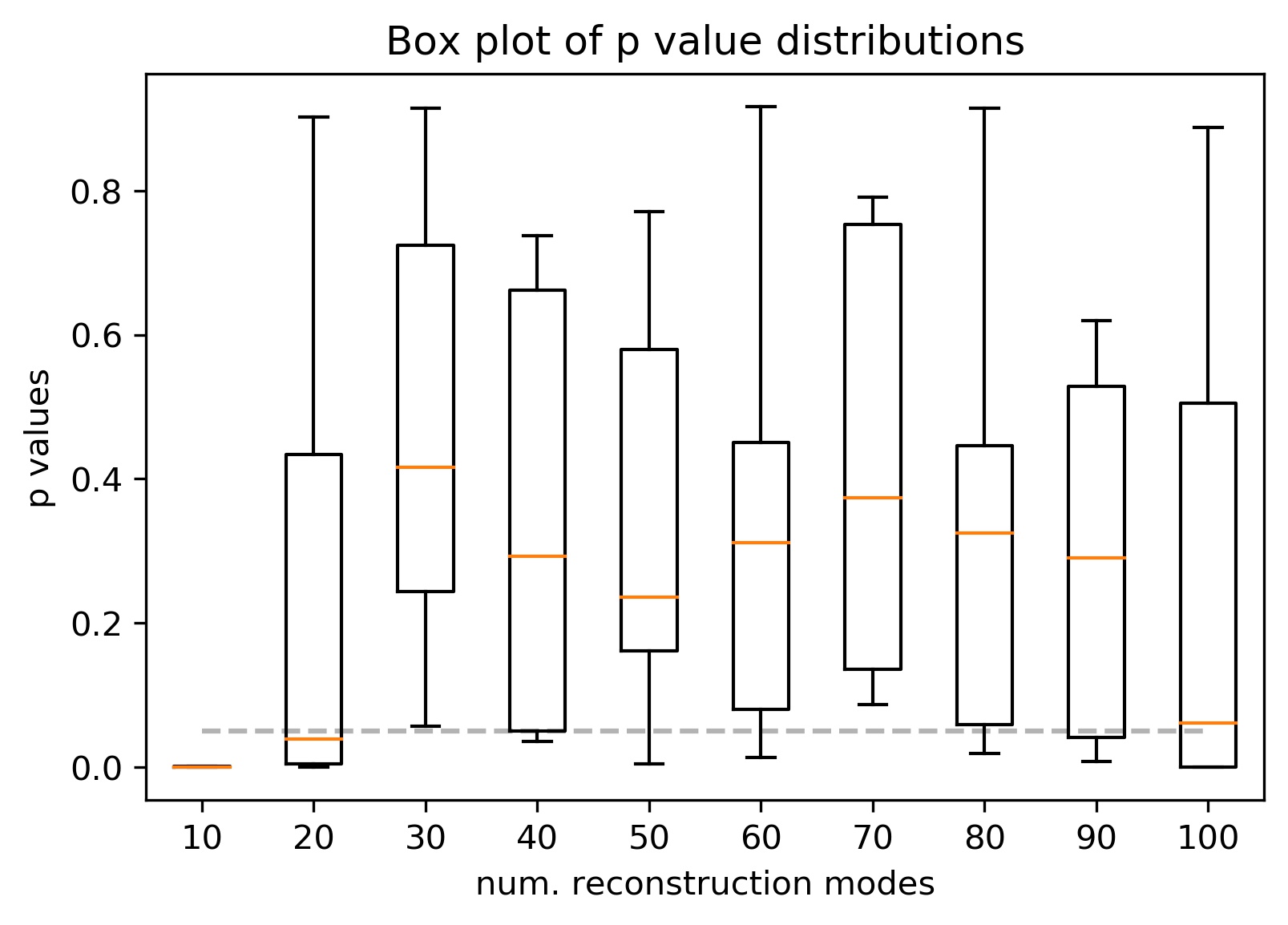}
    \caption{p values boxplot}
\end{subfigure}
\begin{subfigure}[b]{0.45\textwidth}
    \centering
    \includegraphics[width = \textwidth]{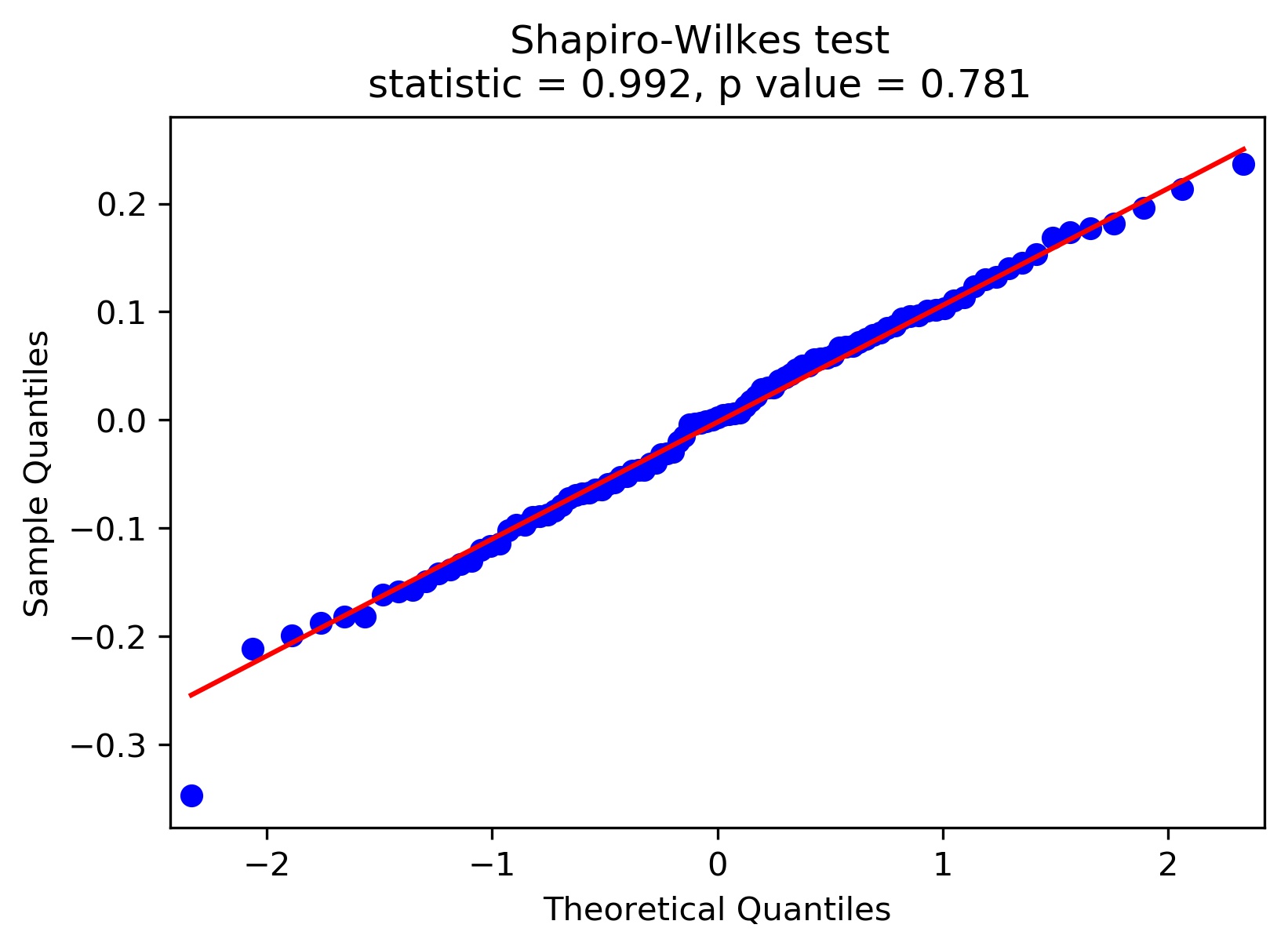}
    \caption{quartile-quartile plot of coordinate 0's modal distribution.}
\end{subfigure}
\begin{subfigure}[b]{0.45\textwidth}
    \centering
    \includegraphics[width = \textwidth]{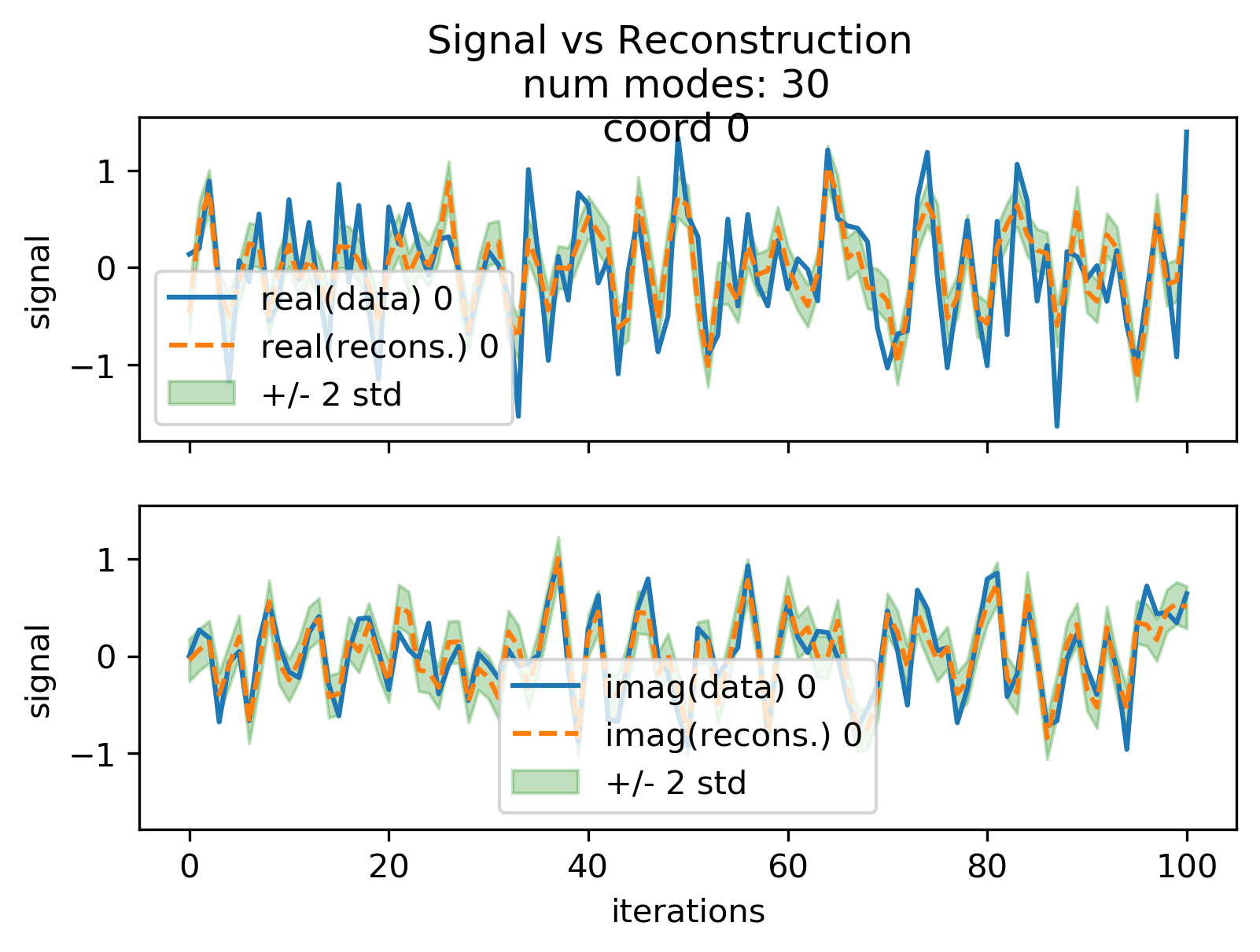}
    \caption{30 mode reconstruction of coordinate 0.}
\end{subfigure}
\caption{\textbf{Linear modal model's p-values:} Pane (a) shows a visual representation of the outputs of the heuristic, namely the mean and median p values for the different sizes of models. The heuristic says that 30 modes is the minimum number to use for the ROM. (b) is similar except it is a boxplot of the modal distributions' p values. (c) quartile-quartile plot for coordinate 0's modal distribution visually testing for normality. (d) coordinate 0's reconstruction with the 30 mode ROM.}
\label{fig:heuristic-linear}
\end{figure}

\begin{figure}[htbp]
\centering
\begin{subfigure}[b]{0.45\textwidth}
    \centering
    \includegraphics[width = \textwidth]{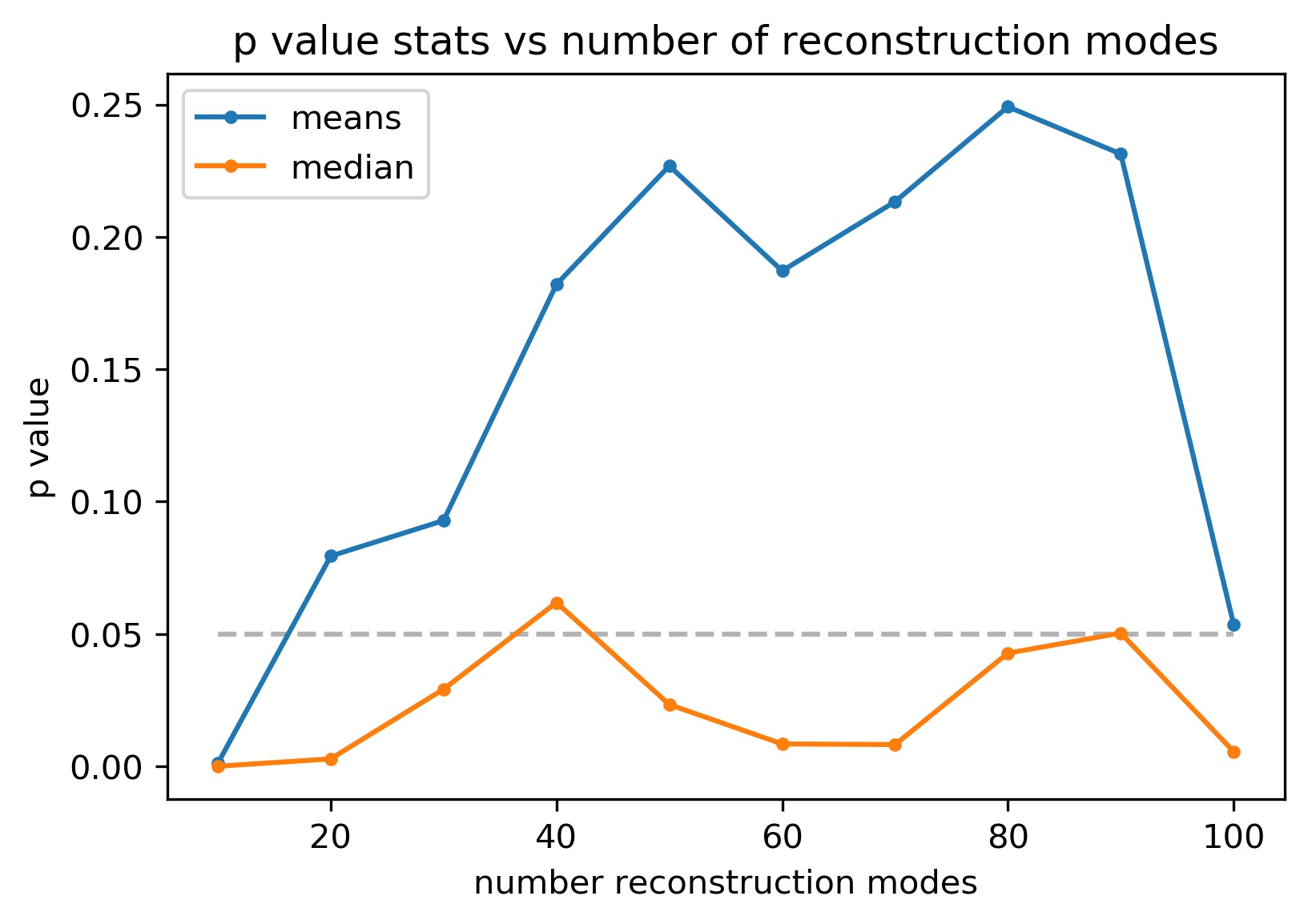}
    \caption{p values mean, median}
\end{subfigure}
\begin{subfigure}[b]{0.45\textwidth}
    \centering
    \includegraphics[width = \textwidth]{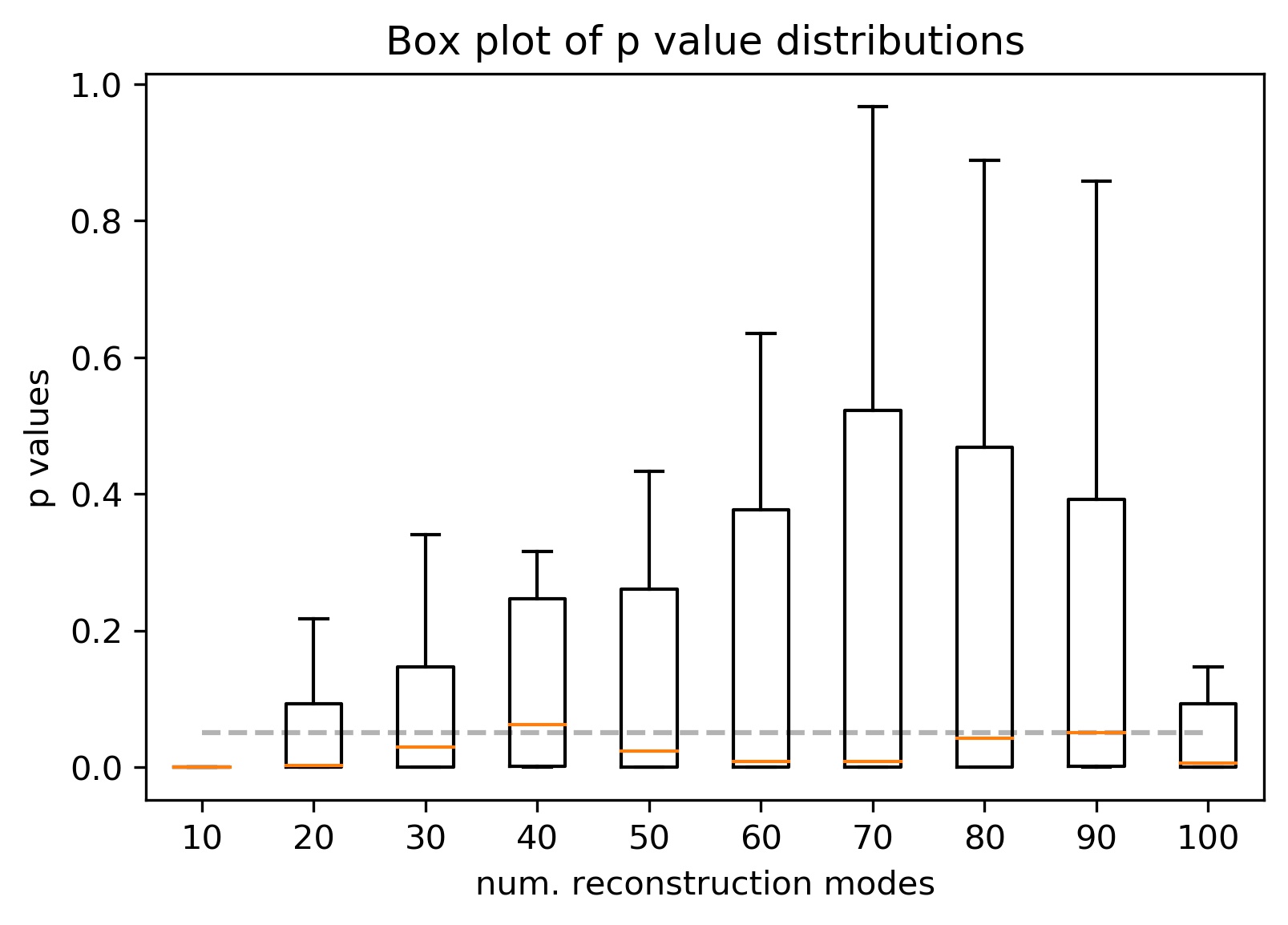}
    \caption{p values boxplot}
\end{subfigure}
\begin{subfigure}[b]{0.45\textwidth}
    \centering
    \includegraphics[width = \textwidth]{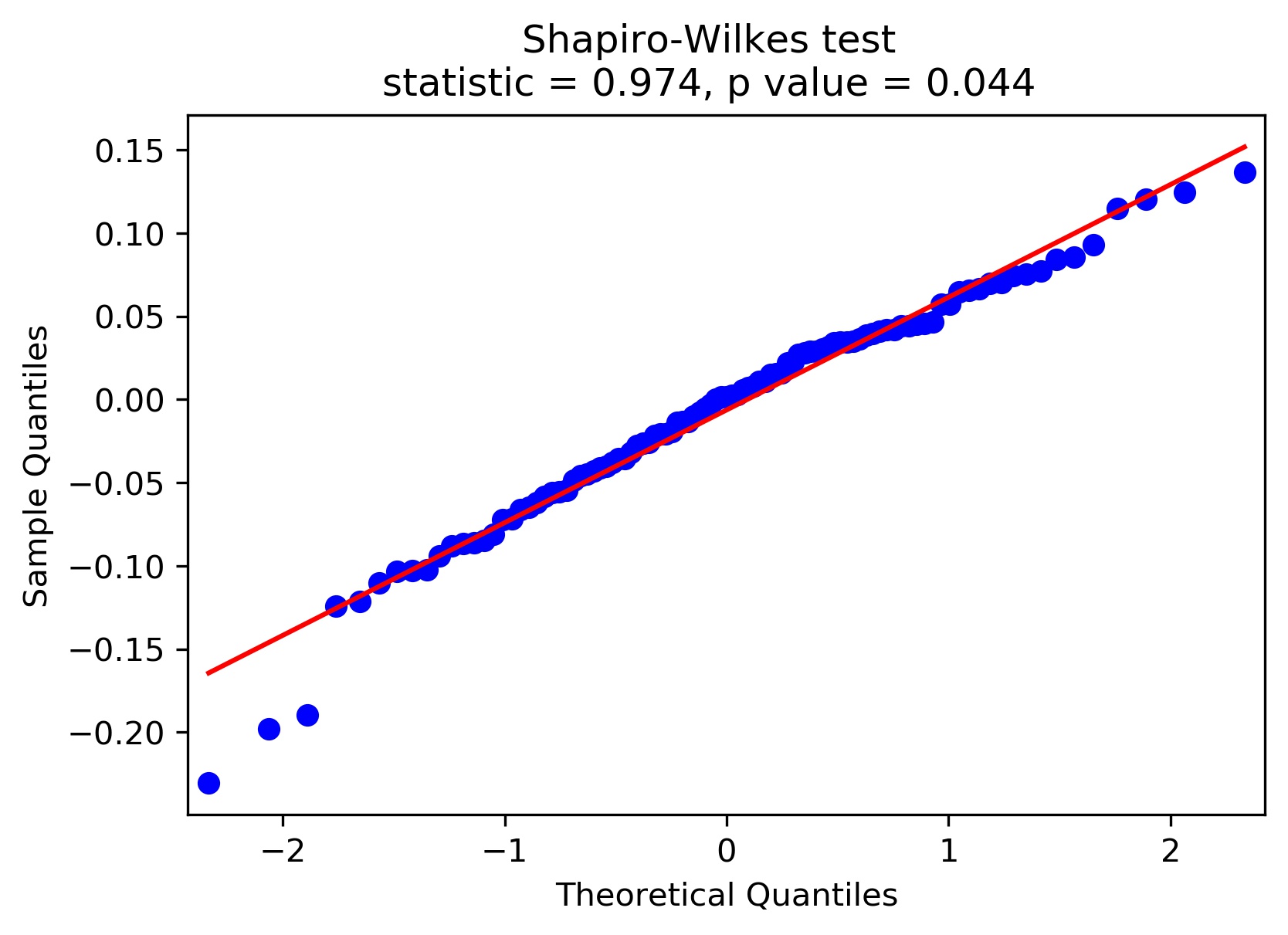}
    \caption{quartile-quartile plot of coordinate 0's modal distribution.}
\end{subfigure}
\begin{subfigure}[b]{0.45\textwidth}
    \centering
    \includegraphics[width = \textwidth]{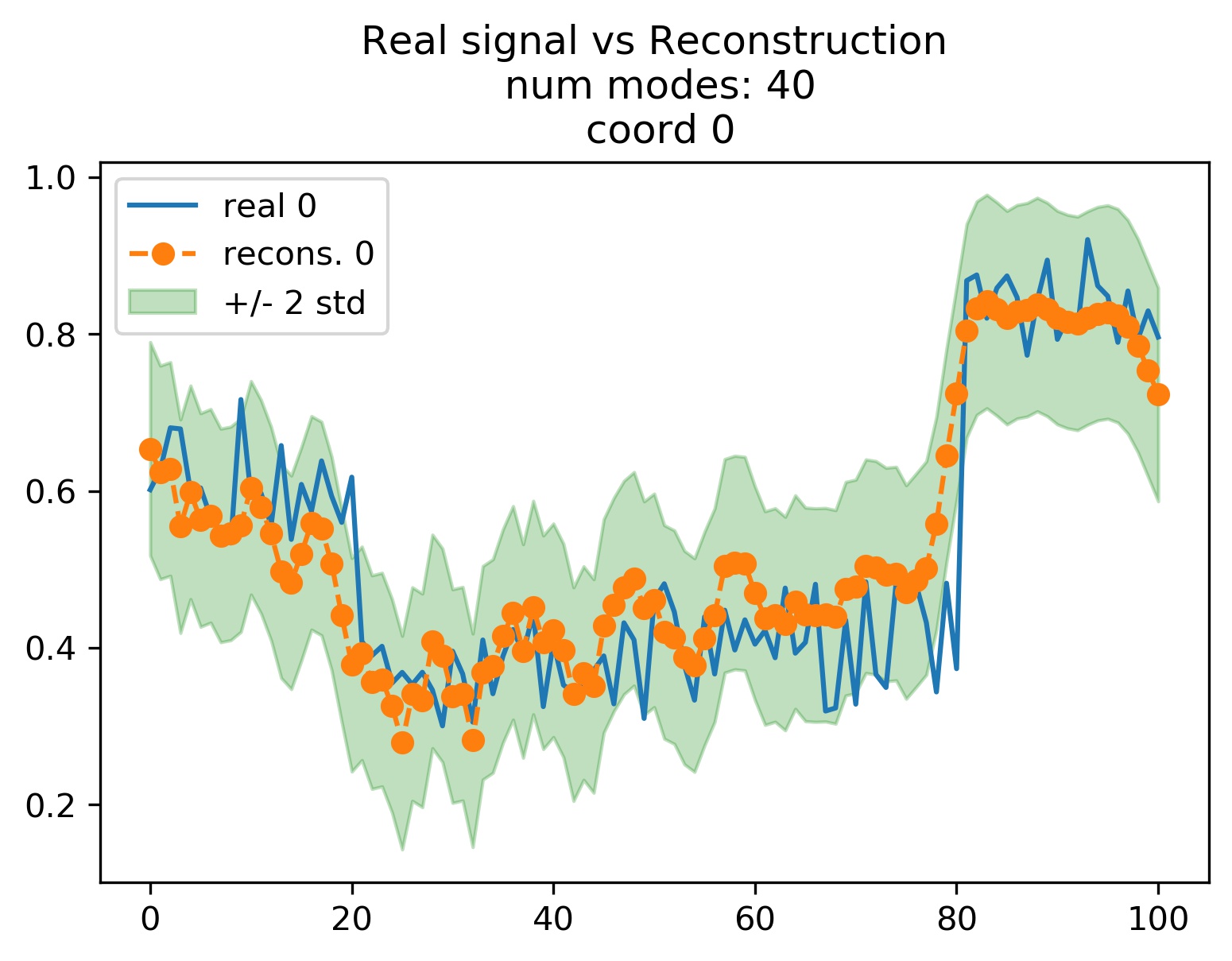}
    \caption{40 mode reconstruction of coordinate 0.}
\end{subfigure}
\caption{\textbf{Anharmonic model's p-values:} Pane (a) shows a visual representation of the outputs of the heuristic, namely the mean and median p values for the different sizes of models. The heuristic says that 30 modes is the minimum number to use for the ROM. (b) is similar except it is a boxplot of the modal distributions' p values. (c) quartile-quartile plot for coordinate 0's modal distribution visually testing for normality. (d) coordinate 0's reconstruction with the 40 mode ROM.}
\label{fig:heuristic-anharmonic}
\end{figure}

\begin{figure}[htbp]
\centering
\begin{subfigure}[b]{0.45\textwidth}
    \centering
    \includegraphics[width = \textwidth]{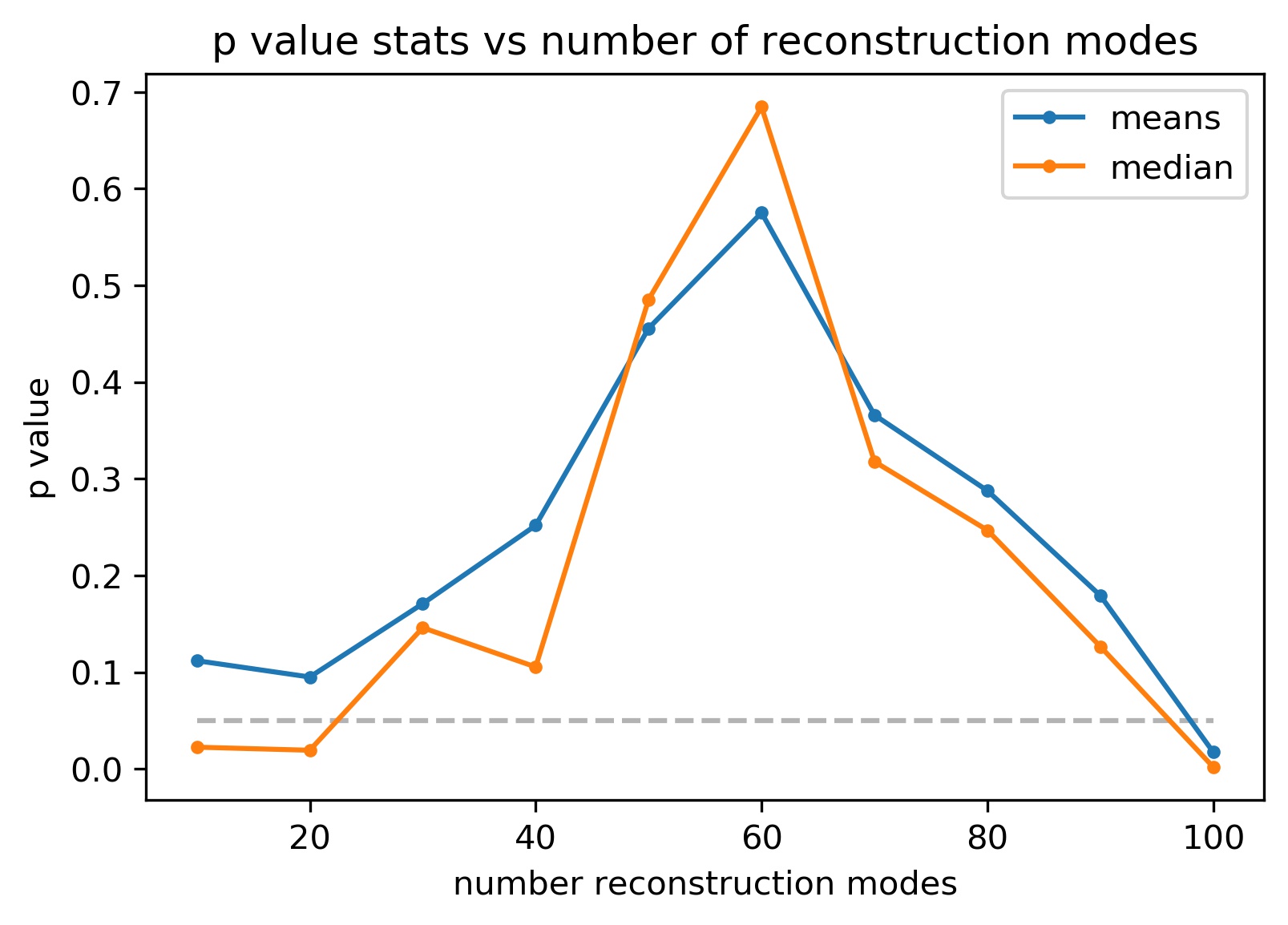}
    \caption{p values mean, median}
\end{subfigure}
\begin{subfigure}[b]{0.45\textwidth}
    \centering
    \includegraphics[width = \textwidth]{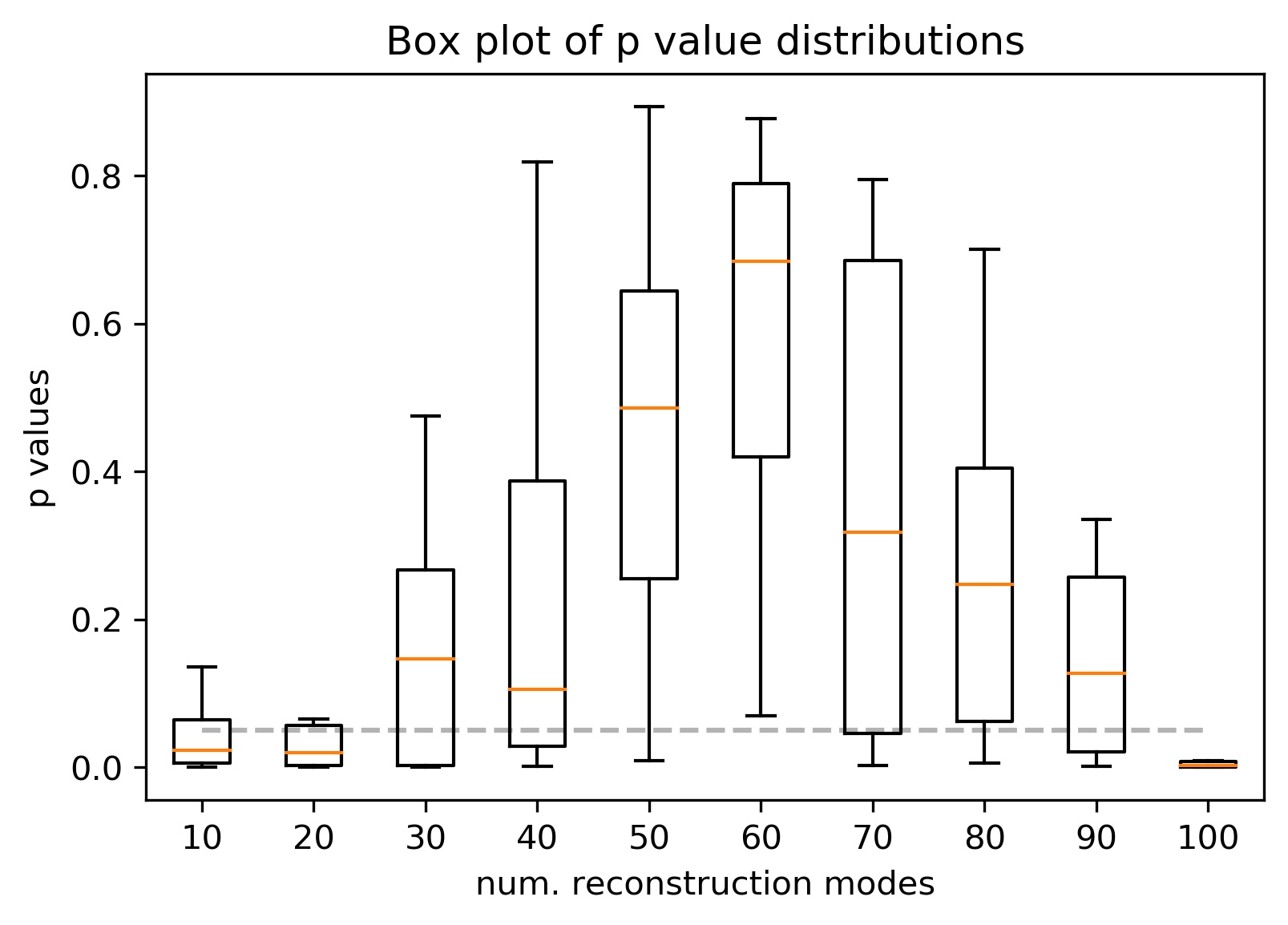}
    \caption{p values boxplot}
\end{subfigure}
\begin{subfigure}[b]{0.45\textwidth}
    \centering
    \includegraphics[width = \textwidth]{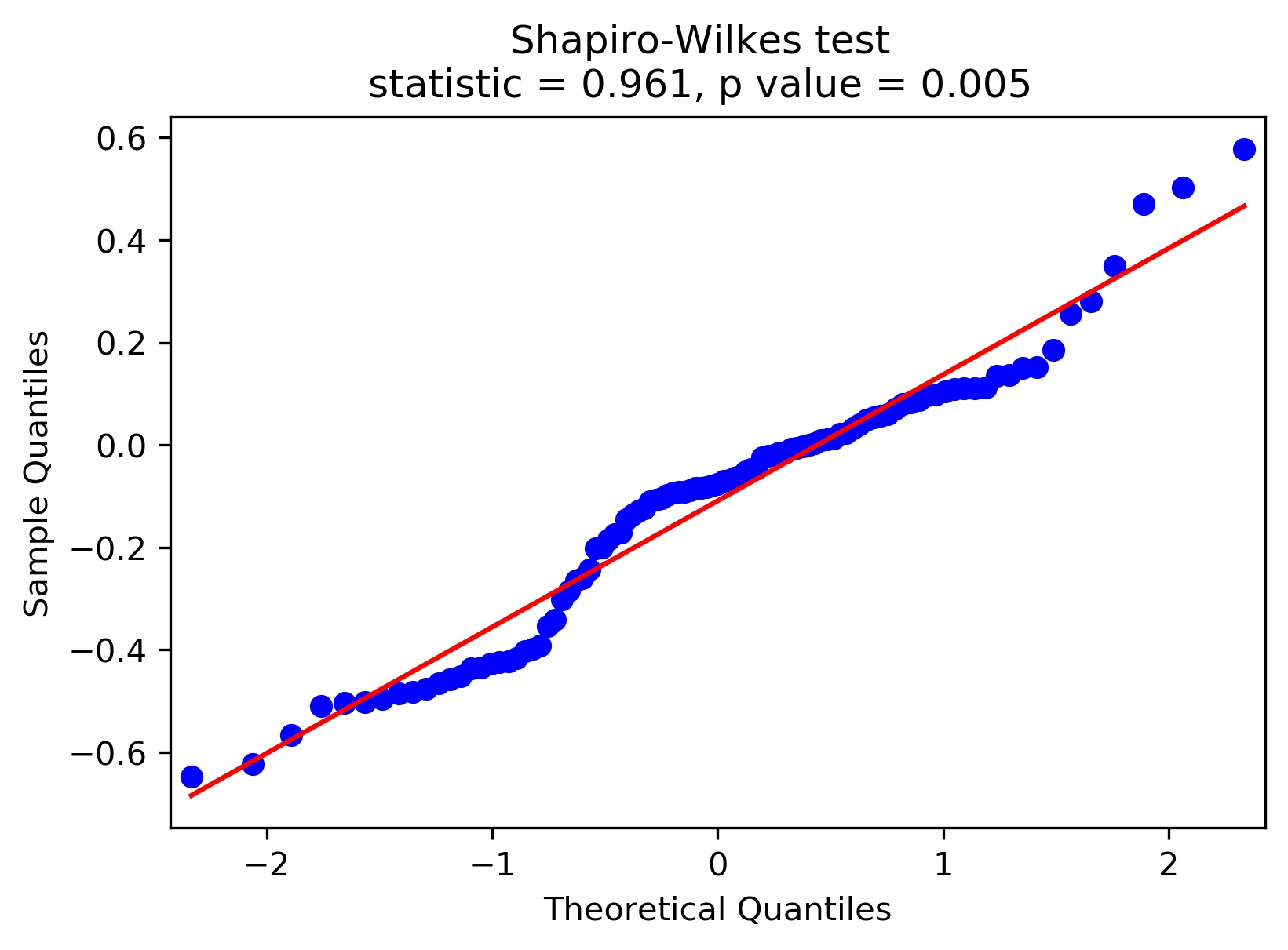}
    \caption{quartile-quartile plot of coordinate 0's modal distribution.}
\end{subfigure}
\begin{subfigure}[b]{0.45\textwidth}
    \centering
    \includegraphics[width = \textwidth]{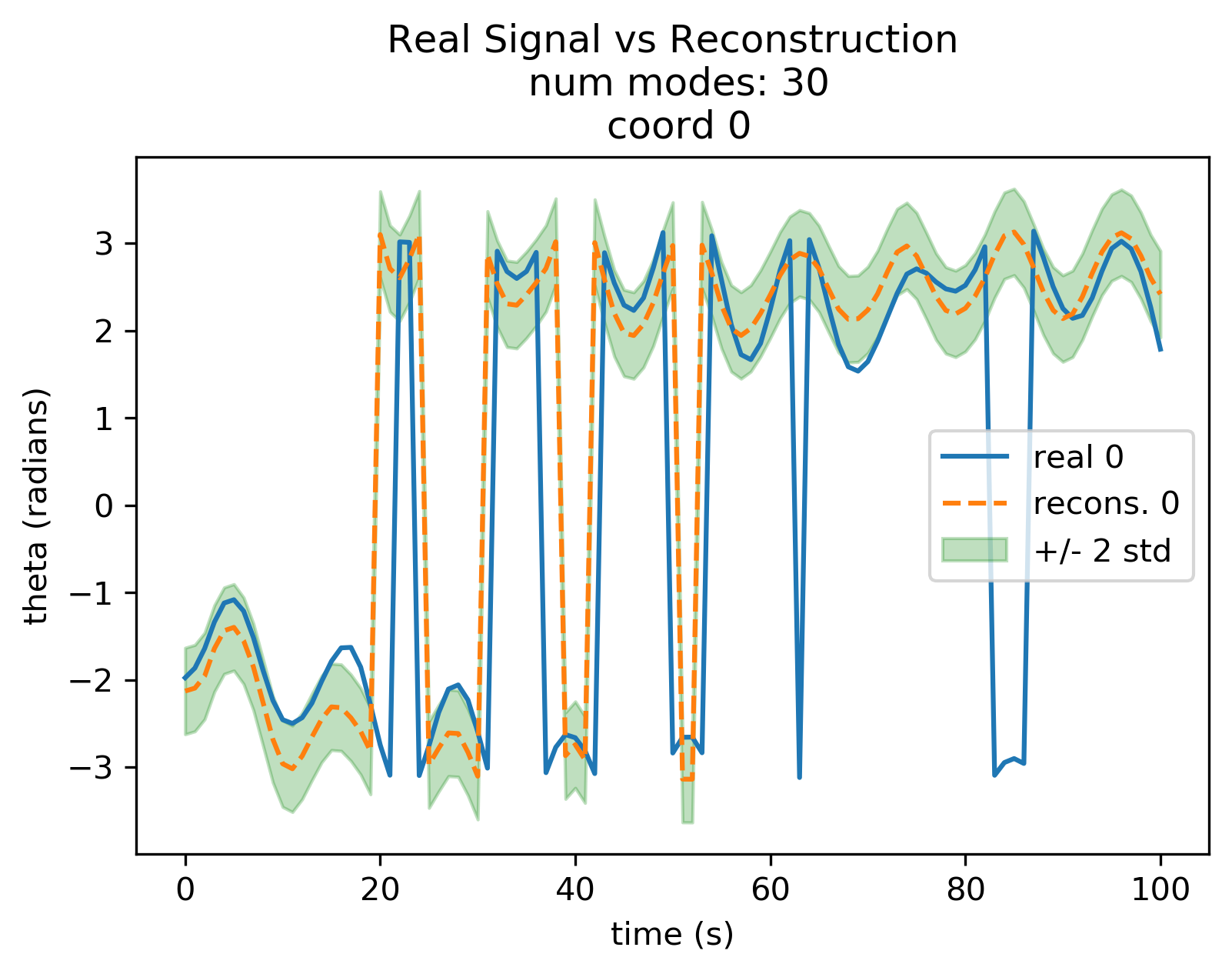}
    \caption{30 mode reconstruction of coordinate 0.}
\end{subfigure}
\caption{\textbf{Kuramoto model's p-values:} In (a) be show a visual representation of the outputs of the heuristic, namely the mean and median p values for the different sizes of models. The heuristic says that 30 modes is the minimum number to use for the ROM. (b) is similar except it is a boxplot of the modal distributions' p values. (c) quartile-quartile plot for coordinate 0's modal distribution visually testing for normality. While the mean and median of the p values are above the 0.05 threshold, this particular coordinate is not approximated well by a normal distribution. (d) coordinate 0's reconstruction with the 30 mode ROM.}
\label{fig:heuristic-kuramoto}
\end{figure}

\clearpage
\section{Discussion and Conclusions}\label{sec:conclusions}

We have developed a method for reduced order modeling using Koopman operator theory that gives confidence bounds on the predictions of the reduced order model.  While the reduced order model is by necessity a finite process, the spectral expansion of the Koopman operator is infinite. The reduced order model represents the process with a finite number of Koopman modes. The rest of the dynamics are modeled as a noise process. The part of the noise process that is in the subspace modeled by the ROM's modes is called the modal process.  This noise process is modeled as a Gaussian noise process and the estimated standard deviation of the noise is used to compute a confidence bound on the ROM's predictions.


We have applied this modeling to a sequence of examples. The first was a synthetic example where we could directly specify the Koopman eigenvalues and modes, so the computed Koopman spectrum could be compared to ground truth. The next two examples represented networked dynamical systems that had changing network topology. The first was a network of noisy, anharmonic oscillators whose connections were switched at each integer time. This system was chosen due to the fact that a single, non-noisy, anharmonic oscillator has a purely continuous spectrum, except the eigenvalue at 1, and thus does not have an expansion in terms of Koopman modes and eigenfunctions. Despite this, the methodology was able to model the system's evolution with a small number of computed modes plus the noise process. The actual trajectory of the system stayed within the confidence bounds of the ROM's prediction a large percentage of the time. The last example consisted of a noisy Kuramoto model with random coupling strength between the oscillators. As with the anharmonic oscillator example, the true signal stayed within the confidence bounds of the ROM's prediction the majority of the time.

As the number of modes increased, the ROM's got closer to a deterministic plus Gaussian noise model. As we pass the threshold, the variance of the Gaussian decreases. This allows us to propose a heuristic algorithm that enables choosing the number of deterministic modes that should be kept.

In this paper, we assumed a Gaussian distribution of the modal noise when constructing the confidence bounds for the ROM predictions. Future work could relax this assumption and use something like kernel density estimation to construct the bounds.


\section*{Acknowledgments}
This material is based upon work supported by the Air Force Office of Scientific Research under award numbers FA9550-17-C-0012 and FA9550-22-1-0531.

\bibliographystyle{siam}

\end{document}